\documentclass[11pt, a4paper]{amsart}

\usepackage{graphicx}
\usepackage{graphics}
\usepackage{enumerate}
\usepackage{latexsym}
\usepackage[all]{xy}
\usepackage{psfrag}
\usepackage[ansinew]{inputenc}
\usepackage[T1]{fontenc}
\usepackage{lmodern}
\xyoption{matrix} \xyoption{arrow}
\usepackage[all]{xy}
\input xy
\xyoption{all}

\usepackage{epsfig}
\usepackage[center,small,sc]{caption}

\usepackage{amscd,amsmath,amsthm, amssymb,fancyhdr,setspace}

\usepackage{pdfsync}
                      
\usepackage{tikz}
\usetikzlibrary{matrix,arrows,backgrounds,shapes.misc,shapes.geometric,patterns,calc,positioning}
\usetikzlibrary{calc,shapes,fadings}
\usetikzlibrary{decorations.pathmorphing,automata}
\usetikzlibrary{decorations.pathreplacing}
\tikzstyle arrowstyle=[scale=1]
\tikzstyle directed=[postaction={decorate,decoration={markings,
    mark=at position .65 with {\arrow[arrowstyle]{stealth}}}}]
\tikzstyle reverse directed=[postaction={decorate,decoration={markings,
    mark=at position .65 with {\arrowreversed[arrowstyle]{stealth};}}}]
\tikzstyle{box} = [rectangle, draw]
\tikzstyle{line} = [draw, -latex']

\newcommand*\circled[1]{\tikz[baseline=(char.base)]{
            \node[shape=circle,draw,inner sep=1pt, line width=.3, scale=.5] (char) {#1};}}

\usepackage{color}
\definecolor{teal}{rgb}{0.0, 0.5, 0.5}
\definecolor{tealblue}{rgb}{0.21, 0.46, 0.53}
\definecolor{tealgreen}{rgb}{0.0, 0.51, 0.5}
\definecolor{tuscanred}{rgb}{0.51, 0.21, 0.21}
\definecolor{sangria}{rgb}{0.57, 0.0, 0.04}
\definecolor{rufous}{rgb}{0.66, 0.11, 0.03}
\definecolor{pinegreen}{rgb}{0.0, 0.47, 0.44}
\definecolor{darkscarlet}{rgb}{0.34, 0.01, 0.1}
\definecolor{darkseagreen}{rgb}{0.56, 0.74, 0.56}
\definecolor{darkpastelred}{rgb}{0.76, 0.23, 0.13}
\definecolor{darkpink}{rgb}{0.91, 0.33, 0.5}
\definecolor{darkpastelblue}{rgb}{0.47, 0.62, 0.8}
\definecolor{alizarin}{rgb}{0.82, 0.1, 0.26}
\definecolor{candyapplered}{rgb}{1.0, 0.03, 0.0}

\newcommand{\calg}{\mathcal{G}}

\newcommand{\ocalg}{\overline{\mathcal{G}}}
\newcommand{\calgSW}{_{SW}\calg}
\newcommand{\calgNE}{\calg^{NE}}

\newcommand{\gi}[3]{G_{#1}, G_{#2}, \ldots, G_{#3}}
\newcommand{\gii}[3]{G'_{#1}, G'_{#2}, \ldots, G'_{#3}}
\newcommand{\re}[2] {\res_{\calg} ( \calg_{#1}, \calg_{#2})}
\newcommand{\gt}[4] {\graft_{{#1}, {#2}} ( \calg_{#3}, \calg_{#4})}
  
\newcommand{\e}{\delta}

\newcommand{\Hom}{\textup{Hom}}
\newcommand{\per}{\textup{per}}
\newcommand{\C}{\mathcal{C}}

\DeclareMathOperator{\res}{Res \,} 
\DeclareMathOperator{\graft}{Graft \,  }

\DeclareMathOperator{\Jac}{J}
\DeclareMathOperator{\End}{End}

\renewcommand{\S}{\mathcal{S}}
\newcommand{\M}{\mathcal{M}}

\renewcommand{\mod}{\mathrm{mod \;}}
\newcommand{\Ext}{\operatorname{Ext}}
\newcommand{\pred}{\operatorname{Pred}}
\newcommand{\suc}{{\operatorname{Succ}}}
\newcommand{\Pred}{\operatorname{Pred}}
\newcommand{\Succ}{{\operatorname{Succ}}}
\newcommand{\Int}{{\operatorname{Int}}}

\newtheorem{Theorem}{Theorem}[section]
\newtheorem{Lemma}[Theorem]{Lemma}
\newtheorem{Proposition}[Theorem]{Proposition}
\newtheorem{Corollary}[Theorem]{Corollary}
\newtheorem{Definition}[Theorem]{Definition}
\newtheorem{Question}[Theorem]{Question}
\theoremstyle{definition}\theoremstyle{plain}
\newtheorem{Remark}[Theorem]{Remark}
\begin{document}

\topmargin -.6cm 
\oddsidemargin .5cm \evensidemargin .5cm

\title[ Extensions in Jacobian Algebras and Cluster Categories]{Extensions in Jacobian Algebras and Cluster Categories of Marked Surfaces \\
}

\thanks{This work was supported through the Engineering and Physical Sciences Research Council, grant number EP/K026364/1, UK}

\author{Ilke Canakci}
\address{Department of Mathematics, University of Leicester, University Road, Leicester LE1 7RH, United Kingdom}
\email{ic74@le.ac.uk}
\author{Sibylle Schroll with an appendix by Claire Amiot} 
\address{Department of Mathematics, University of Leicester, University Road, Leicester LE1 7RH, United Kingdom}
\email{schroll@le.ac.uk}

\address{Institut Fourier-UMR 5582, 100 rue des maths,
38402 Saint Martin d'H\'eres, France}
\email{claire.amiot@univ-grenoble-alpes.fr}

\begin{abstract}
In the context of representation theory of finite dimensional algebras, 
string algebras have been extensively studied and most aspects of their representation theory are well-understood. One exception to this is the 
classification of extensions
between indecomposable modules. In this paper we explicitly describe such extensions for a  class of
string algebras, namely gentle algebras associated to surface triangulations. These algebras arise as Jacobian algebras of unpunctured surfaces. 
We relate the extension spaces of indecomposable modules to crossings of generalised arcs in the surface and give explicit bases of the extension spaces for indecomposable modules in  almost all cases. We  show that the dimensions of these extension spaces are given in terms of crossing arcs in the surface.

Our approach is new and consists  of interpreting 
snake graphs as indecomposable modules. 
In order to show that our basis is a spanning set,  we need to work in the 
associated cluster category where we explicitly
calculate the middle terms of  extensions and give bases of their extension spaces. We note that not all 
extensions in the cluster category give rise to extensions for the Jacobian 
algebra.

\end{abstract}

\date{\today}

\maketitle


\section{Introduction}

\parskip7pt
\parindent0pt

Cluster algebras were introduced by Fomin and Zelevinsky in 2002 in \cite{FZ1} in order to give an algebraic framework for the study of the (dual)  canonical bases in Lie theory. This work was further developed in \cite{BFZ, FZ2,FZ4}. Cluster algebras are commutative algebras 
given by generators, the \emph{cluster variables},  and relations. The construction of the generators is a recursive process 
from  an initial set of data. In general, even in 
small cases, this is a complex  process. However, there is a class of cluster algebras coming from surfaces \cite{FST, FT} (see also \cite{FG1,FG2}) where this process is 
encoded in the combinatorial geometry  of surface triangulations. Surface cluster algebras are an important part of the classification of 
(skew-symmetric) cluster algebras in terms of mutation
type, namely almost all cluster algebras of finite mutation type are surface cluster algebras  \cite{FeShTu}. 

Cluster algebras from surfaces have been widely studied via the combinatorial geometry of the corresponding surfaces \cite{FeShTu, FST, FT,  MSW, MSW2}. The 
same holds true for the associated cluster categories and Jacobian algebras. An important example of this is the crossing of two arcs in a surface. 
In the case of cluster algebras this gives rise to a multiplication formula for the corresponding cluster variables \cite{MW}. 
In the cluster category, the number 
of crossings of two arcs gives the dimension of the extension space between the associated indecomposable objects \cite{ZZZ, QZ}. For Jacobian algebras
of surfaces 
where all marked points lie in the boundary, in \cite{BZ}, building 
on \cite{BR} and \cite{ABCP}, Auslander-Reiten sequences  have been given in terms of arcs in the surface. 
In general, however, there has so far been no 
link between arbitrary crossings of arcs in the surface and the extensions between indecomposable modules in the Jacobian algebra. 

The Jacobian algebras under consideration are gentle algebras  
and their indecomposable modules, given by strings and bands, correspond to curves and closed curves in the surface \cite{ABCP} (see \cite{L1,Lab2,CL,L2} for a more general definition of Jacobian algebras via quiver with potential  and \cite{GLS} for classification of their representation type).  Gentle algebras form a special class of algebras, for example,  this class is closed under tilting and derived 
equivalence \cite{S} and \cite{SZ}.  They are part of the  larger family of string algebras which  are 
an important family of algebras of tame representation type whose  representation theory is well-understood. For example, their  Auslander-Reiten 
structure has been  determined \cite{BR} and in 
\cite{CB,K} the morphisms between indecomposable modules 
are completely described.
However, a complete description of the extensions between 
indecomposable modules is not known.

In the present paper, we describe extension spaces of  string modules over gentle
Jacobian algebras. Furthermore, we show that in analogy with 
the cluster category, in most cases, the number of crossings of two arcs still gives rise to a dimension formula of the extension space between the corresponding indecomposable modules in the Jacobian algebra by explicitly constructing a basis. 
However, not every crossing contributes to this dimension. We characterize exactly which crossings
contribute and which do not. We do this by introducing a new approach, consisting of using the snake graph calculus developed in \cite{CS} to explicitly construct extensions  resulting from crossing arcs. 
This gives a lower bound on the 
dimensions of the extensions spaces in the Jacobian algebra. 
In order to obtain an upper bound, we work in the cluster category. There we show explicitly how, in many cases, the four arcs in the surface, resulting from the smoothing 
of a crossing of two not necessarily distinct arcs, give rise to two extensions in the cluster category.

For a surface cluster algebra, the cluster variables are in bijection with arcs in the surface \cite{FST}. Moreover, in a given 
triangulation, each arc corresponds to a combinatorial object called a \emph{snake graph} \cite{MS, MSW,  Propp}.  Snake graphs have proven to be an 
important element in the understanding of surface cluster algebras, for example, in \cite{MSW2} snake graphs (and band graphs) were used to show that certain collections of 
loops and (generalised) arcs comprise  vector space bases  
for surface cluster algebras. Snake graphs have also been instrumental in the proof of the positivity conjecture for 
surface cluster algebras \cite{MSW}. 
Note that the conjecture has since been proved for all skew-symmetric cluster algebras \cite{LS}.  

If two (generalised) arcs $\gamma_1$ and $\gamma_2$ in a marked surface $(S,\M)$ cross then the geometric operation of smoothing the crossing is given by 
locally replacing the crossing {\Large $\times$} 
with the pair of segments   $\genfrac{}{}{0pt}{5pt}{\displaystyle\smile}{\displaystyle\frown}$  or with the pair of segments { $\supset \subset$}. 
This gives rise to four  new arcs $\gamma_3, \gamma_4$ and $\gamma_5, \gamma_6$ corresponding to the two different ways of smoothing the crossing. 
The corresponding elements $x_{\gamma_1},  \ldots, x_{\gamma_6}$ in the cluster algebra satisfy the so-called \emph{skein relations} given by  
$x_{\gamma_1} x_{\gamma_2} = y_- x_{\gamma_3} x_{\gamma_4} + y_+ x_{\gamma_5} x_{\gamma_6}$ where $y_-, y_+$ are some coefficients  \cite{MW}.

Suppose from now on that $(S,\M)$ is a marked surface such that  all marked points are in the boundary of $S$ and 
let $T$ be a triangulation of $(S,\M)$. All arcs are considered to be generalised, that is they  might have  self-crossings. 
The abstract snake graph calculus developed in \cite{CS} applies in this setting and 
gives a combinatorial interpretation in terms of snake graphs of the  arcs resulting from the smoothing of two crossing arcs. 
We remark that we never actually smooth self-crossings as in \cite{CS2}, but instead in the case of a self-crossing
we  consider two copies of the same arc.  We then
use the combinatorial description in \cite{CS} to study the extension space over the associated Jacobian algebra $J(Q,W)$ and over the cluster category $\C(S,\M)$
defined in 
\cite{A} by giving explicit bases for these spaces in almost all cases. As mentioned above, the string modules over  $J(Q,W)$ are in bijection with the arcs
in the surface \cite{ABCP} not contained in $T$, and the arcs in the surface correspond in turn to snake graphs \cite{MSW}.  Therefore there is a  correspondence associating a
snake graph corresponding to an arc in $(S, \M, T)$ to the string module corresponding to the same arc and this defines a sign function on the snake graph, see  Proposition~\ref{SnakeStrings}. 

Based on the snake graph calculus developed in \cite{CS}, given two string modules, we define three types of crossings of modules corresponding to 
the three different types of crossings of the associated  snake graphs. 
Namely, if, in the language of \cite{CS2}, the snake graphs cross with an overlap then we say that the corresponding string modules \emph{cross in a module}. 
If the snake graphs cross with grafting and  $s=d$ (see Section 2.4 for the definition of $s$ and $d$) then we say that the corresponding string modules \emph{cross in an arrow} and 
finally if the snake graphs cross with grafting and $s \neq d$ where $s$ and $d$ are parameters associated to the snake graphs  then the corresponding string modules \emph{cross in a 3-cycle}. 
Our first result is then to determine when two  
crossing string modules $M_1$ and $M_2$  give rise to a non-zero element in $\Ext^1_{J(Q,W)}(M_1, M_2)$.

{\bf Theorem~\ref{ShortExactSequences}}\label{ShortExactSequencesIntro} \emph{Let $M_1$ and $M_2$ be two  string modules (not necessarily distinct) over $J(Q,W)$ corresponding to arcs
    $\gamma_1$ and $\gamma_2$ in $(S,\M)$. Then for a given crossing 
    of $M_1$ and $M_2$ corresponding to a crossing of $\gamma_1$ and $\gamma_2$,
    there are string modules $M_3$ and $M_4$ obtained by `smoothing the crossing' of $\gamma_1$ and $\gamma_2$ such that there exists an extension of $M_1$ by $M_2$
    \begin{itemize}
    \item[(1)]  with two non-zero middle terms given by $M_3$ and $M_4$ if and only if $M_1$ crosses $M_2$ in a module at this crossing,
    \item[(2)] with one non-zero middle term given by $M_3$ if and only if $M_1$ crosses $M_2$ in an arrow at this crossing. 
    \end{itemize}
    When  $M_1$ crosses  $M_2$ in a 3-cycle, $M_3$ and $M_4$ do not give an extension of $M_1$ by $M_2$ corresponding to this crossing.}

A geometric interpretation of the three different types of crossings in Theorem~\ref{ShortExactSequences} is given in Remark~\ref{Geometric interpretation of crossings}.

Note that there is a direction in the crossing of modules, that is `$M_1$ crosses $M_2$' is different from `$M_2$ crosses $M_1$' and  that this 
distinction does not appear in terms of crossings of the corresponding  arcs, see Section~\ref{CrossingModulesSection}. 

We remark that Theorem~\ref{ShortExactSequences} can be interpreted as skein relations for string modules  and such skein relations have been announced in \cite{GLS2} in the setting of Caldero-Chapoton algebras.

In the cluster category $\C(S,\M),$ the indecomposable objects correspond to arcs and (non-contractible)  closed loops in $(S,\M)$ and therefore they are referred to as string and band objects, respectively, see \cite{BZ}. It  follows from the results on AR-triangles in \cite{BZ} that all triangles between indecomposable objects in $\C (S,\M)$ have at most two middle terms. In  \cite{ZZZ} it
is shown  that the dimension of the extension space of  two string  objects in the cluster category is equal to the number of crossings of the corresponding arcs (see also \cite{BM} for a similar result for tube categories).  This suggests a close connection between the geometric crossing of arcs and the extension spaces.

Indeed, we show in Theorem~\ref{Triangles} (see below) that the middle terms of triangles in the cluster category arise from smoothing the crossings of  arcs in the surface in almost all cases. More precisely, consider the two pairs of arcs $\gamma_3, \gamma_4$ and $\gamma_5, \gamma_6$ obtained from smoothing a crossing of two arcs
$\gamma_1$ and $\gamma_2$ (with a suitable orientation). We show that the 
pair $\gamma_3, \gamma_4$ always gives rise to an element in $\Ext_\C(\gamma_1, \gamma_2)$. We show that the other pair of arcs $\gamma_5, \gamma_6$ gives rise
to an element in $\Ext_\C(\gamma_2, \gamma_1)$ if the crossing of $\gamma_1$ and $\gamma_2$  itself does not have a {\it self-crossing } (such a self-crossing is given by a {\it self-crossing overlap} in terms of the corresponding snake graphs. See Section 2, Definition 2.11 and Theorem 2.12 for the definition of  a self-crossing overlap of snake graphs and the correspondence with crossings of arcs). 
We remark that an important factor in the proof of Theorem~\ref{Triangles} is  
the geometric interpretation of Iyama-Yoshino 
reduction \cite{IY} given by Marsh-Palu \cite{MP}.

{\bf Theorem~\ref{Triangles}} \label{TrianglesIntro} {\it Let $\gamma_1$ and $\gamma_2$ be two  string objects (not necessarily distinct) in $\C(S,\M)$
    such that their corresponding arcs cross in $(S,\M)$.
    Let $\gamma_3, \gamma_4, \gamma_5, \gamma_6$ be the string objects corresponding to 
    the smoothing of a particular crossing of a suitable orientation of the corresponding arcs 
    $\gamma_1$ and $\gamma_2$. Then there is a non-split triangle in $\C(S,\M)$ given by
    \begin{equation}
    \gamma_2 \longrightarrow \gamma_3 \oplus \gamma_4 \longrightarrow \gamma_1 \longrightarrow \gamma_2[1]
    \end{equation} 
    and  if the crossing of $\gamma_1$ and $\gamma_2$ is not in a self-crossing  overlap in some triangulation of $(S,\M)$  then we obtain a non-split triangle given by 
    \begin{equation} \gamma_1 \longrightarrow \gamma_5 \oplus \gamma_6 \longrightarrow \gamma_2 \longrightarrow \gamma_1[1]. \end{equation}
    If any of $\gamma_3, \gamma_4, \gamma_5, \gamma_6$ are  boundary arcs, the corresponding objects $\C(S,\M)$  are zero objects. }

We raise the question (Question~\ref{Question}) of what the middle terms are of those triangles in (2) above  corresponding to the crossings   of $\gamma_1$ and $\gamma_2$  that is a self-crossing overlap. 
That there must also be a non-zero extension $\Ext_\C(\gamma_2, \gamma_1)$ corresponding in some way to this crossing follows from the 2-Calabi Yau property of $\C(S, \M)$ and from the result in \cite{ZZZ} giving the dimension of the extension space in terms of the number of crossing of the arcs.

Theorem~\ref{Geometric interpretation of crossings}  together with  Theorem~\ref{Triangles} and the dimension formula in \cite{ZZZ}, give the following result for extensions in the Jacobian algebra. Here $\Int(\gamma, \delta)$ is the minimal number of 
intersections of two arcs $\gamma$ and $\delta$. We remark that we use the following convention:
if $\gamma = \delta$ then $\Int(\gamma, \gamma)= 2 m$ where $m$ is the minimal number of self-crossings of $\gamma$.

{\bf Corollary 4.2} {\it Let $M, N$ be two   string modules over $J(Q,W)$ and let $\gamma_M$ and $\gamma_N$ be the corresponding arcs in $(S,\M)$ such that $\gamma_M$ and $\gamma_N$ have no crossing with self-crossing overlap.  
    
    (1) A basis of $\Ext^1_{J(Q,W)}(M,N)$ is given 
    by all short exact sequences arising from $M$ crossing $N$ 
    in a module or an arrow and where the middle terms are as described in Theorem~\ref{ShortExactSequences};
    
    (2) We have 
    $$\dim \Ext^1_{J(Q,W)}(M,N) + \dim \Ext^1_{J(Q,W)}(N, M) = \Int(\gamma_M, \gamma_N) - k-k'$$
    where $k$ (resp. $k'$) is the number of times that $M$ crosses $N$ (resp.  $N$ crosses $M$)
    in a 3-cycle. In particular, if $M=N$ we have 
    $$2 \dim \Ext^1_{J(Q,W)}(M, M) = \Int(\gamma_M, \gamma_M) - 2k.$$
}

{\bf Acknowledgments:} We  would like to thank Ralf Schiffler for helpful conversations. We also would like to thank the referee for pointing out a gap in an earlier version of the proof of Theorem 4.1 and for the helpful comments and suggestions that they make.

\section{Background on Jacobian algebras, cluster categories and snake graphs}

Throughout let  $k$ be an algebraically closed field. 

\subsection{Bordered marked surfaces}

In this section we follow \cite{CS, MSW} in our exposition. Let $S$ be a connected oriented 2-dimensional Riemann surface with non-empty boundary. Let $\M$ be a finite set of marked points on $S$ such that all marked points lie in the boundary of $S$ 
and each boundary component contains at least one marked point. 
Call the pair $(S,\M)$ a \emph{(bordered) marked surface}. If $S$ is a disc then let $|\M| \geq 4$. 

\begin{Definition}
A \emph{generalised arc} in $(S,\M)$ is a curve $\gamma$ in $S$, considered up to homotopy, such that

(1) the endpoints of $\gamma$ are in $\M,$

(2) except for the endpoints $\gamma$ is disjoint from the boundary of $S,$

(3) $\gamma$ does not cut out a monogon or a bigon. 

The curve  $\gamma$ is called an \emph{arc}, considered up to homotopy, if it satisfies (1), (2), (3), and 

(4) $\gamma$ does not cross itself, except that its endpoints might coincide. 
\end{Definition} 

A generalised arc may cross itself a finite number of times. 

A \emph{boundary segment} is the homotopy class of a curve that lies in the boundary and connects two (not necessarily distinct) neighbouring marked points on the same boundary component.
Note that a boundary segment is not considered to be an arc. However, we sometimes do refer to it as a boundary arc.

\begin{Definition} For two arcs $\gamma$, $\gamma'$ in $(S,\M)$, let $\rm{Int}(\gamma, \gamma')$ be the minimal number of crossings of curves $\alpha$, $\alpha'$ where $\alpha$ and 
$\alpha'$ range over the homotopy classes of $\gamma$ and $\gamma'$, respectively. We say that arcs  $\gamma$, $\gamma'$ are compatible if $\rm{Int}(\gamma, \gamma')=0$.
\end{Definition}

\begin{Definition}
A \emph{triangulation} of $(S,\M)$ is a maximal collection of pairwise compatible arcs. 
A \emph{flip} of an arc $\gamma$ in a triangulation $T$ of $(S,\M)$ replaces the arc $\gamma$ with the unique arc $\gamma'$ such that $T\setminus \{\gamma\} \cup \{\gamma'\}$ is a triangulation
of $(S,\M)$. 
\end{Definition}

All triangulations of $(S,\M)$ are connected by a series of flips.

\begin{Definition}\label{SmoothingCrossingArcs}
Let $\gamma_1$ and $\gamma_2$ be generalised arcs such that $\gamma_1$ and $\gamma_2$ cross at a point $x$. We define the \emph{smooting of the crossing of $\gamma_1$ and $\gamma_2$ at the point
    $x$} to be the pairs of arcs $\{\alpha_1, \alpha_2\}$ and $\{\beta_1, \beta_2\}$ where 

- $\{\alpha_1, \alpha_2\}$ is the same as $\{\gamma_1, \gamma_2\}$ except locally where the crossing {\Large $\times$} is replaced with the pair of segments 
$\genfrac{}{}{0pt}{5pt}{\displaystyle\smile}{\displaystyle\frown},$ 

- $\{\beta_1, \beta_2\}$ is the same as $\{\gamma_1, \gamma_2\}$ except locally where the crossing {\Large $\times$} is replaced with the pair of segments { $\supset \subset$}.    
\end{Definition}

We remark that if we consider oriented arcs then the orientation of $\gamma_1$ and $\gamma_2$ permits to distinguish the set of arcs $\{\alpha_1, \alpha_2\}$ from the set of arcs $\{  \beta_1$, $\beta_2\}$.

From now on we will not make a distinction between arcs and generalised arcs and we will simply call them arcs unless otherwise specified.

\subsection{Gentle algebras from surface triangulations}

In this section we recall the definition of gentle algebras and introduce some related notation which we will be using throughout the paper.

Let $Q = (Q_0, Q_1)$ be a quiver, denote by $kQ$ its path algebra and for an admissible ideal $I$, let $(Q,I)$ be the associated bound quiver. Denote by 
$\mod A$ the module category of finitely generated right $A$-modules of an algebra $A$.

\begin{Definition} An algebra $A$ is \emph{gentle} if it is Morita equivalent to an algebra $kQ/I$ such that

(S1) each vertex of $Q$ is the starting point of at most two arrows and is the end point of at most two arrows; 

(S2) for each arrow $\alpha$ in $Q_1$ there is at most one arrow $\beta$ in $Q_1$ such that $\alpha\beta$ is not in $I$ and there is at most one arrow $\gamma$ in $Q_1$
such that $\gamma\alpha$ is not in $I$; 

(S3) $I$ is generated by paths of length 2;

(S4) for each arrow $\alpha$ in $Q_1$ there is at most one arrow $\delta$ in $Q_1$ such that $\alpha \delta$ is in $I$ and there is at most one   arrow $\varepsilon$ in $Q_1$ such that
$\varepsilon \alpha$ is in $I$. 
\end{Definition}

For $\alpha \in Q_1$, let $s(\alpha)$ be the start of $\alpha$ and $t(\alpha)$ be the end of $\alpha$.

For each arrow $\alpha$ in $Q_1$ we define the formal inverse $\alpha^{-1}$ such that $s(\alpha^{-1}) = t(\alpha)$ and $t(\alpha^{-1}) = s(\alpha)$. A word $w  =
\varepsilon_1 \varepsilon_2 \ldots \varepsilon_n$ is a \emph{string} if either $\varepsilon_i$ or
$\varepsilon_i^{-1}$ is an arrow in $Q_1$, if $s(\varepsilon_{i+1}) = t(\varepsilon_{i}),$  if
$\varepsilon_{i+1} \neq \varepsilon_{i}^{-1}$ for all $ 1 \leq i \leq n-1$ and if no subword of $w$ or its inverse is in $I$. Let $s(w) = 
s(\varepsilon_1)$ and $t(w) = t(\varepsilon_n)$.
Denote by $\mathcal{S}$ the set of strings modulo the equivalence relation $w \sim w^{-1}$, where $w$ is a string.

A string $w$ is a direct
string if $w = \alpha_1 \alpha _2 \ldots \alpha_n$ and $\alpha_i \in Q_1$ for all $ 1 \leq i \leq n$ and $w$ is an inverse string if $w^{-1}$ is a direct string.

The terminology of string modules, in particular, the notions of hooks and cohooks were defined in  \cite{BR}.  However, the definitions
of hooks and cohooks we give here differ slightly from the usual definitions. More precisely, our hooks and cohooks do not necessarily satisfy the maximality
conditions on direct and inverse strings appearing in the standard literature. 

Given a string $w$, define  four substrings ${_hw}, {w_h}, {_cw}, {w_c}$ of $w$  as follows:

We say ${_hw}$ is obtained from $w$ by deleting a hook on $s(w)$ where
$$ {_hw} = \left\{ \begin{array}{ll}
0 & \mbox{if $w$ is an inverse string,} \\
{_hw} & \mbox{where ${_hw}$ is obtained from $w$ by deleting the first direct arrow in $w$} \\
&  \mbox{and the inverse string preceding it.}
\end{array}\right.$$ 
We say ${_cw}$ is obtained from $w$ by deleting a cohook on $s(w)$ where
$${_cw} = \left\{ \begin{array}{ll}
0 & \mbox{if $w$ is a direct string,} \\
{_cw} & \mbox{where ${_cw}$ is obtained from $w$ by deleting the first inverse arrow in $w$} \\
&  \mbox{and the direct string preceding it.}
\end{array}\right.$$
We say ${w_h}$ is obtained from $w$ by deleting a hook on $t(w)$ where
$${w_h} = \left\{ \begin{array}{ll}
0 & \mbox{if $w$ is a  direct string,} \\
{w_h} & \mbox{where ${w_h}$ is obtained from $w$ by deleting the last inverse arrow in $w$} \\
&  \mbox{and the direct string succeeding  it.}
\end{array}\right.$$
We say ${w_c}$ is obtained from $w$ by deleting a cohook on $t(w)$ where
$${w_c} = \left\{ \begin{array}{ll}
0 & \mbox{if $w$ is an inverse string,} \\
{w_c} & \mbox{where ${w_c}$ is obtained from $w$ by deleting the last direct arrow in $w$} \\
&  \mbox{and the inverse string succeeding  it.}
\end{array}\right.$$

Let $T$ be a triangulation of $(S,\M)$ with associated quiver with potential $(Q,W)$ and let $J(Q,W)$ be the associated Jacobian algebra as defined in \cite{L1}. As recalled in the introduction, given a marked surface
where all marked points lie in the boundary, this algebra coincides with 
the gentle algebra defined in \cite{ABCP}.
Let  $\S$ be the set of all strings in $J(Q,W)$. 
Given a string $w \in \S$ we denote by $M(w)$ the corresponding string module in $J(Q,W)$. Note that $M(w) \simeq M(w^{-1})$, for a string $w$.
Conversely, given a string module $M$ we can associate to it a string (or its inverse).  That is, there exists a string  $w_M$ such that $M \simeq M(w_M)$ in which case we also have $M \simeq M(w_M^{-1})$.
The string corresponding to a simple module at vertex $i$ of $Q$ is denoted by $i$, that is, it is given by the single vertex $i$.
Given an arc $\gamma$ in the surface  in \cite{ABCP} a  string $w_\gamma$ is associated to an orientation of $\gamma$. We denote  $ M(w_\gamma)$  the associated string module. The opposite orientation of $\gamma$ gives rise to the inverse string $(w_\gamma)^{-1}$ and we have that  $ M(w_\gamma) \simeq   M((w_\gamma)^{-1})$. 
Conversely,  by \cite{ABCP} any string
module $M $ is associated to an arc  $\gamma_M$ in $( S, \M, T)$. 

\begin{figure}[ht]
\begin{tikzpicture}[trans/.style={thick,->,shorten >=2pt,shorten <=2pt,>=stealth}]

\node[scale=.4](a) at (0,2) {$\bullet$};
\node[scale=.4](b) at (0,0)  {$\bullet$};
\node[scale=.4](a') at (5,2) {$\bullet$};
\node[scale=.4](b') at (5,0) {$\bullet$};

\node[scale=.4](c) at (-1,1) [label={[label distance=-.05cm, scale=.7, color=red]180:$s(\gamma)$}, scale=.7] {$\bullet$};
\node[scale=.4, color=red](c') at (6,1) [label={[label distance=-.05cm, scale=.7, color=red]0:$e(\gamma)$}, scale=.7] {$\bullet$};

\node[scale=.4] (a1) at ( $(a)!0.25!(a') $){$\bullet$};
\node[scale=.4] (a2) at ( $(a)!0.5!(a') $){$\bullet$};
\node[scale=.4] (a3) at ( $(a)!0.8!(a') $){$\bullet$};

\node[scale=.4] (b1) at ( $(b)!0.3!(b') $){$\bullet$};
\node[scale=.4] (b2) at ( $(b)!0.5!(b') $){$\bullet$};
\node[scale=.4] (b3) at ( $(b)!0.75!(b') $){$\bullet$};

\node[scale=.6] at (-.1,1.2){$\tau_1$};
\node[scale=.6] at (.8,1.2){$\tau_2$};
\node[scale=.6] at (1.45,1.2){$\tau_3$};
\node[scale=.6] at (2,1.2){$\tau_4$};
\node[scale=.6] at (2.6,1.2){$\tau_5$};

\node[scale=.6] at (3.8,1.2){$\tau_{d-2}$};
\node[scale=.6] at (4.4,1.25){$\tau_{d-1}$};
\node[scale=.6] at (5.1,1.2){$\tau_{d}$};

\draw (a.center) -- (b.center) (a'.center) -- (b'.center) (a.center) -- (c.center) (c.center) -- (b.center) (a'.center) -- (c'.center) (c'.center) -- (b'.center)
(a.center)--(b1.center)--(a1.center) (b1.center)--(a2.center) (b2.center)--(a2.center) (b3.center)--(a3.center)--(b3.center)
(a.center)--(a2.center) (a3.center)--(a'.center) (b.center)--(b2.center) (b3.center)--(b'.center) (a'.center)--(b3.center);
\draw [dotted, line width=.8] (a2.center)--(a3.center) (b2.center)--(b3.center);

\path[line width=.75, color=red] (c.center) edge node[pos=.6, fill=white,outer sep=1mm, scale=.6]{$\gamma$} (c'.center);

\draw[trans, color=brown] ($(a)!0.25!(b1)$)--($(a)!0.25!(b)$)node [midway,yshift=-5pt,scale=.6] {$\alpha_1$};
\draw[trans, color=brown] ($(a)!0.75!(b1)$)--($(a1)!0.75!(b1)$)node [midway,yshift=-5pt,scale=.6] {$\alpha_2$};
\draw[trans, color=brown] ($(a1)!0.75!(b1)$)--($(a2)!0.75!(b1)$)node [midway,yshift=-5pt,scale=.6] {$\alpha_3$};
\draw[trans, color=brown] ($(a2)!0.3!(b2)$)--($(a2)!0.3!(b1)$)node [midway,yshift=-5pt,scale=.6] {$\alpha_4$};
\draw[trans, color=brown] ($(a3)!0.7!(b3)$)--($(a')!0.7!(b3)$)node [midway,yshift=-5pt,xshift=3pt,scale=.6] {$\alpha_{d-2}$};
\draw[trans, color=brown] ($(a')!0.3!(b')$)--($(a')!0.3!(b3)$)node [midway,yshift=3pt,xshift=-.3pt,scale=.6] {$\alpha_{d-1}$};

\end{tikzpicture}

\caption{Denoting the vertices in the quiver corresponding to an arc $\tau_i$ of the triangulation also by $\tau_i$ {z}, the string corresponding to the arc $\gamma$ is given by  }
\label{FigStringConstruction}
\end{figure}
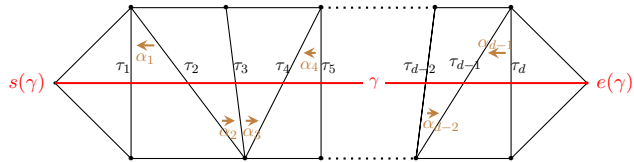

For the convenience of the reader we briefly recall in Figure~\ref{FigStringConstruction} the construction of a string given an orientated arc in a triangulated surface as defined in \cite{ABCP}, see also \cite{BZ}.

\vspace{-.1cm}

\begin{center}

\begin{tikzpicture}[auto]
\node at (-.8,0){$w(\gamma)=$};
\node (a0) at (0,0) {$\tau_1$};
\node (a1) at (1,0) {$\tau_2$};
\node (a2) at (2,0) {$\tau_3$};
\node (a3) at (3,0) {$\tau_4$};
\node (a4) at (4,0) {$\tau_5$};
\node at (4.5,0){$\cdots$};
\node (ad-2) at (5.2,0){$\tau_{d-2}$};
\node (ad-1) at (6.5,0){$\tau_{d-1}$};
\node (ad) at (7.7,0){$\tau_{d}$};
\draw[<-] (a0) to  node[scale=.7] {$\alpha_1$} (a1);
\draw[->] (a1) to  node[scale=.7] {$\alpha_2$} (a2);
\draw[->] (a2) to  node[scale=.7] {$\alpha_3$} (a3);
\draw[<-] (a3) to  node[scale=.7] {$\alpha_4$} (a4);
\draw[->] (ad-2) to  node[scale=.7] {$\alpha_{d-2}$} (ad-1);
\draw[<-] (ad-1) to  node[scale=.7] {$\alpha_{d-1}$} (ad);
\end{tikzpicture}

\end{center}


\subsection{Cluster categories of marked surfaces}


Cluster categories were first introduced in \cite{BMRRT} for acyclic quivers and independently in \cite{CCS} for type $A$. 
Generalised cluster categories were defined in \cite{A}.
Namely, given a quiver with potential $(Q,W)$ such that the Jacobian algebra $J(Q,W)$ is finite dimensional, denote by $\Gamma := \Gamma(Q,W)$ the associated 
Ginzburg dg-algebra.   Consider the perfect derived category  $\per \; \Gamma$ which is the smallest triangulated  subcategory of the derived category 
$\mathcal{D} (\Gamma)$ containing $\Gamma$ which is
stable under taking direct summands
and consider the bounded derived category $\mathcal{D}^\flat(\Gamma)$ of $\Gamma$.  The generalized cluster category $\C(Q,W)$ is the 
quotient $\per \; \Gamma / \mathcal{D}^\flat (\Gamma)$.
It is shown in \cite{A} that $\C(Q,W)$ 
is Hom-finite, 2-Calabi-Yau,  the image of $\Gamma$ in $\C(Q,W)$ is a cluster tilting object $T_\Gamma$, and the  
endomorphism algebra of $T_\Gamma$ is isomorphic to the Jacobian algebra $J(Q,W)$.
Furthermore, the categories $\C(Q,W)/T_\Gamma$ and $\mod  J(Q,W)$ are equivalent and the functor $\Hom_{\C(Q,W)}(T_\Gamma [-1], -)$ is the projection functor from $\C(Q,W)$
to $\mod  J(Q,W)$ \cite{KZ}.

Now let $(S,\M)$ be a marked surface,  $T, T'$  triangulations of $(S,\M)$, and let $(Q,W)$ and $(Q',W')$ be the quivers with potential 
associated with $T$ and $T',$ respectively. 
It follows from \cite{FST, KY, L1} that $\C(Q,W)$ and $\C(Q',W')$ are triangle equivalent and hence the cluster category is independent of the 
triangulation of $(S,\M)$. We will thus denote the cluster category by $\C(S,\M)$ or $\C$.

In \cite{BZ} the cluster category $\C{(S,\M)}$ associated to a surface with marked points on the boundary is explicitly described. 
In particular, a parametrization of the 
indecomposable objects of $\C(S,\M)$ is
given in terms of string objects and band objects. The string objects correspond bijectively to the homotopy classes of 
non-contractible curves in $(S,\M)$ that are not 
homotopic to a boundary segment of $(S,\M)$ and subject to the equivalence relations $\gamma \sim \gamma^{-1}$. The band objects correspond bijectively to 
the elements of $k^* \times \Pi_1^*(S,\M)/\sim$ where $\Pi_1^*(S,\M)$ are the invertible elements of the fundamental group of $(S,\M)$ 
and where $\sim$ is the equivalence
relation generated by $\gamma \sim \gamma^{-1} $ and cyclic permutation of $\gamma$.

Furthermore, it is shown in \cite{BZ} that the $AR$-translation of an indecomposable object $\gamma$ corresponds to simultaneously rotating
the start and end points of $\gamma$ in the orientation of $(S,\M)$.

Unless otherwise stated we will not distinguish between arcs and the corresponding indecomposable objects in $\C(S,\M)$.

The following theorem plays a crucial role in our results. 

\begin{Theorem}\cite[Theorem 3.4]{ZZZ}\label{ZZZ}
Let $\gamma$ and $\delta$ be two (not necessarily distinct) arcs in $(S,\M)$. Then $$\dim_k \Ext^1_\C(\gamma, \delta) = \Int( \gamma, \delta).$$
\end{Theorem}

\subsection{Snake Graphs}\label{SnakeSection}
In this section we state and prove two results relating to snake graphs, namely Propositions~\ref{self-crossing equals crossing overlap} and~\ref{SnakeStrings}.
These  results form the basis of many of the proofs in sections 3 and 4. 
For the convenience of the reader we recall in this section all relevant results on snake graphs that we refer to in later sections. 
We define snake graphs associated to triangulations of surfaces as in \cite{MSW} and \cite{CS, CS2}. Below, we closely follow the 
exposition in \cite{CS2} adapting it to snake graphs associated to surface triangulations. 

Let $T$ be a triangulation of $(S,\M)$ and $\gamma$ be an arc in $(S,\M)$ which is not in $T$.
Choose an orientation of $\gamma$.  We choose $\gamma$ to be  a representative in its homotopy class which transversally intersects the arcs of $T$ such that no arc of $T$ is crossed twice in succession.
Let $\tau_{i_1} , \ldots, \tau_{i_d}$ be  the arcs of $T$ crossed by $\gamma$ in the order 
given by the orientation of $\gamma$. Note that we choose $\gamma$ to be a representative in its homotopy 
class such that $\gamma$  has a minimal number of intersections with 
$\tau_{i_1} , \ldots, \tau_{i_d}$. It is however possible that $\tau_{i_j} = \tau_{i_k}$ for $j \neq k$.
For an arc
$\tau_{i_j}$, let $\Delta_{j-1}$ and $\Delta_j$ be the two triangles in $(S,\M,T)$ that share the arc $\tau_{i_j}$ and such that $\gamma$ first crosses $\Delta_{j-1}$ and then $\Delta_j$.
Note that each $\Delta_j$ always has three distinct sides, but that two or all three of the vertices of $\Delta_j$ might be identified. Let $G_j$ be the graph with 4 vertices and 
5 edges, having the shape of a square (with a fixed side length) with a diagonal that satisfies the following property: 
there is a bijection of the edges of $G_j$ and the 5 distinct arcs in the triangles $\Delta_{j-1}$ and $\Delta_j$ 
and such that the diagonal in $G_j$ corresponds to the arc $\tau_{i_j}$. That is, $G_j$ corresponds to the quadrilateral with diagonal $\tau_{i_j}$ formed by $\Delta_{j-1}$ and $\Delta_j$
in $(S, \M, T)$.

Given a planar embedding $\tilde{G}_j$ of $G_j$, we define the \emph{relative orientation} ${\rm{rel}}(\tilde{G}_j, T)$ of $\tilde{G}_j$ with respect to $T$ to be $1$ or $-1$
depending on 
whether the triangles in $\tilde{G}_j$ agree or disagree with the (common) orientation of the triangles $\Delta_{j-1}$ and $\Delta_j$ in $(S, \M, T)$. 

Using the notation above, the arcs $\tau_{i_j}$ and $\tau_{i_{j+1}}$ form two edges of the triangle $\Delta_j$. Let $\sigma_j$ be the third arc in this triangle. 
We now recursively glue together the tiles $G_1, \ldots, G_d$ one by one from 1 to $d$ in the following way: choose planar embeddings of the $G_j$ such that  ${\rm{rel}}(\tilde{G}_j, T)
\neq  {\rm{rel}}(\tilde{G}_{j+1}, T)$. Then glue $\tilde{G}_{j+1}$ to $\tilde{G}_j$ along the edge labelled $\sigma_j$. 

After gluing together the $d$ tiles $G_1, \ldots, G_d$, we obtain a graph (embedded in the plane) which we denote by $\mathcal{G}_\gamma^\Delta$. 

\begin{Definition}
The \emph{snake graph} $\mathcal{G}_\gamma$ associated to $\gamma$ is obtained from $\mathcal{G}_\gamma^\Delta$ by removing the diagonal in each tile. If $\tau \in T$ then 
we define the associated snake graph $\mathcal{G}_\tau$ to be the graph consisting of one single edge with two distinct vertices (regardless of whether the endpoints of $\tau$ are
distinct or not). 
\end{Definition}

\begin{figure}
\begin{tikzpicture}
\node at (-.5,-.5){\begin{tikzpicture}[scale=2.1, trans/.style={thick,->,shorten >=2pt,shorten <=2pt,>=stealth}]
    \draw (0,0) circle (1cm) ;
    \draw (0,0) circle  (.3cm) ;
    \node (a) at (90:1cm) [scale=.4] {$\bullet$};
    \node (b)at (270:1cm) [scale=.4] {$\bullet$};
    \node (c) at (90:.3cm) [scale=.4] {$\bullet$};
    \path (a.center) edge node [pos=.5, fill=white,outer sep=1mm,scale=.5]{$3$} (c.center);
    \draw (a.center) .. controls +(-.4,-.2) and +(-1.5,-.1) .. ($(a.center)!.8!(b.center)$) .. controls +(.5,0) and +(.7,.2) .. node [pos=.1, fill=white,outer sep=1mm,scale=.5]{$1$} (c.center);
    \draw (c.center) .. controls +(.4,.8) and +(1,.5) .. node [pos=.7, fill=white,outer sep=1mm,scale=.5]{$2$} (b.center);
    \draw[color=red, line width=.8] (a.center) .. controls +(-.2,0) and +(-1.2,0) .. node [pos=.7, fill=white,outer sep=1mm,scale=.5]{$\gamma$} ($(a.center)!.8!(270:.85cm)$) .. controls +(.4,0) and +(.7,0) .. ($(a.center)!.5!(90:.1cm)$) .. controls +(-.3,0) and +(-.05,.05)..  ($(a.center)!.8!(210:.6cm)$) .. controls +(.05,-.8) and +(0,-.05)..  ($(a.center)!.8!(-20:.8cm)$) .. controls (.7,.2)  and (.5,.5).. (a.center);
    \end{tikzpicture}};

\node at (4,1){\begin{tikzpicture}[scale=.5]
    \draw (0,0)--(1,0)--(1,1)--(0,1)--(0,0);
    \draw[densely dotted] (0,1)--(1,0);
    \node at (.5,.5){$1$};
    \draw[line width=1] (1,0)--(1,1);
    \node at (.5,0)[below, scale=.6]{$3$};
    \node at (.5,1)[above, scale=.6]{$2$};
    \end{tikzpicture}};

\node at (5.3,.87){\begin{tikzpicture}[scale=.5]
    \draw (0,0)--(1,0)--(1,1)--(0,1)--(0,0);
    \draw[densely dotted] (0,1)--(1,0);
    \node at (.5,.5){$2$};
    \draw[line width=1] (0,0)--(0,1);
    \node at (1,.5)[right, scale=.6]{$3$};
    \node at (.5,0)[below, scale=.6]{$1$};
    \end{tikzpicture}};

\draw[-open triangle 45] (6,1) -- node[rotate=0,above, scale=.7] {glue} (7,1);

\node at (8,1){\begin{tikzpicture}[scale=.5]
    \draw (0,0)--(1,0)--(1,1)--(0,1)--(0,0) (1,0)--(2,0)--(2,1)--(1,1);
    \draw[densely dotted] (0,1)--(1,0) (1,1)--(2,0);
    \node at (.5,.5){$1$};
    \node at (1.5,.5){$2$};
    \draw[line width=1] (1,0)--(1,1);
    \node at (.5,0)[below, scale=.6]{$3$};
    \node at (.5,1)[above, scale=.6]{$2$};
    \node at (2,.5)[right, scale=.6]{$3$};
    \node at (1.5,0)[below, scale=.6]{$1$};
    \end{tikzpicture}};

\node at (4,-.5){\begin{tikzpicture}[scale=.5]
    \draw (0,0)--(1,0)--(1,1)--(0,1)--(0,0) (1,0)--(2,0)--(2,1)--(1,1);
    \draw[densely dotted] (0,1)--(1,0) (1,1)--(2,0);
    \node at (.5,.5){$1$};
    \node at (1.5,.5){$2$};
    \node at (.5,0)[below, scale=.6]{$3$};
    \node at (.5,1)[above, scale=.6]{$2$};
    \node at (2,.5)[right, scale=.6]{$3$};
    \node at (1.5,0)[below, scale=.6]{$1$};
    \draw[line width=1.2, dashed] (1,1)--(2,1);
    \end{tikzpicture}};

\node at (5.3,-.5){\begin{tikzpicture}[scale=.5]
    \draw (0,0)--(1,0)--(1,1)--(0,1)--(0,0);
    \draw[densely dotted] (0,1)--(1,0);
    \node at (.5,.5){$3$};
    \draw[line width=1.2, dashed] (0,0)--(1,0);
    \node at (1,.5)[right, scale=.6]{$1$};
    \node at (0,.5)[left, scale=.6]{$2$};
    \end{tikzpicture}};

\draw[-open triangle 45] (6,-.5) -- node[rotate=0,above, scale=.7] {glue} (7,-.5);

\node at (8,-.5){\begin{tikzpicture}[scale=.5]
    \draw (0,0)--(1,0)--(1,1)--(0,1)--(0,0) (1,0)--(2,0)--(2,1)--(1,1) (1,1)--(1,2)--(2,2)--(2,1);
    \draw[densely dotted] (0,1)--(1,0) (1,1)--(2,0) (1,2)--(2,1);
    \draw[line width=1.2, dashed] (1,1)--(2,1);
    \node at (.5,.5){$1$};
    \node at (1.5,.5){$2$};
    \node at (1.5,1.5){$3$};
    \node at (.5,0)[below, scale=.6]{$3$};
    \node at (.5,1)[above, scale=.6]{$2$};
    \node at (2,.5)[right, scale=.6]{$3$};
    \node at (1.5,0)[below, scale=.6]{$1$};
    \node at (1,1.5)[left, scale=.6]{$2$};
    \node at (2,1.5)[right, scale=.6]{$1$};
    \end{tikzpicture}};

\node at (4,-2.5){\begin{tikzpicture}[scale=.5]
    \draw (0,0)--(1,0)--(1,1)--(0,1)--(0,0) (1,0)--(2,0)--(2,1)--(1,1) (1,1)--(1,2)--(2,2)--(2,1);
    \draw[densely dotted] (0,1)--(1,0) (1,1)--(2,0) (1,2)--(2,1);
    \draw[line width=1] (1,2)--(2,2);
    \node at (.5,.5){$1$};
    \node at (1.5,.5){$2$};
    \node at (1.5,1.5){$3$};
    \node at (.5,0)[below, scale=.6]{$3$};
    \node at (.5,1)[above, scale=.6]{$2$};
    \node at (2,.5)[right, scale=.6]{$3$};
    \node at (1.5,0)[below, scale=.6]{$1$};
    \node at (1,1.5)[left, scale=.6]{$2$};
    \node at (2,1.5)[right, scale=.6]{$1$};
    \end{tikzpicture}};

\node at (5.3, -2.5){\begin{tikzpicture}[scale=.5]
    \draw (0,0)--(1,0)--(1,1)--(0,1)--(0,0);
    \draw[densely dotted] (1,0)--(0,1);
    \node at (.5,.5){$1$};
    \draw[line width=1] (0,0)--(1,0);
    \node at (0,.5)[left, scale=.6]{$3$};
    \node at (1,.5)[right, scale=.6]{$2$};
    \end{tikzpicture}};

\draw[-open triangle 45] (6,-2.5) -- node[rotate=0,above, scale=.7] {glue} (7,-2.5);

\node at (8,-2.5){\begin{tikzpicture}[scale=.5]
    \draw (0,0)--(1,0)--(1,1)--(0,1)--(0,0) (1,0)--(2,0)--(2,1)--(1,1) (1,1)--(1,2)--(2,2)--(2,1) (1,2)--(1,3)--(2,3)--(2,2);
    \draw[densely dotted] (0,1)--(1,0) (1,1)--(2,0) (1,2)--(2,1) (1,3)--(2,2);
    \node at (.5,.5){$1$};
    \node at (1.5,.5){$2$};
    \node at (1.5,1.5){$3$};
    \node at (1.5,2.5){$1$};
    \node at (.5,0)[below, scale=.6]{$3$};
    \node at (.5,1)[above, scale=.6]{$2$};
    \node at (2,.5)[right, scale=.6]{$3$};
    \node at (1.5,0)[below, scale=.6]{$1$};
    \node at (1,1.5)[left, scale=.6]{$2$};
    \node at (2,1.5)[right, scale=.6]{$1$};
    \node at (1,2.5)[left, scale=.6]{$3$};
    \node at (2,2.5)[right, scale=.6]{$2$};
    \end{tikzpicture}};

\node at (-1,-5){$\begin{tikzpicture}[scale=.5]
    \draw (0,0)--(1,0)--(1,1)--(0,1)--(0,0) (1,0)--(2,0)--(2,1)--(1,1) (1,1)--(1,2)--(2,2)--(2,1) (1,2)--(1,3)--(2,3)--(2,2);
    \draw[densely dotted] (0,1)--(1,0) (1,1)--(2,0) (1,2)--(2,1) (1,3)--(2,2);
    \draw[line width=1.2, dashed] (1,3)--(2,3);
    \node at (.5,.5){$1$};
    \node at (1.5,.5){$2$};
    \node at (1.5,1.5){$3$};
    \node at (1.5,2.5){$1$};
    \node at (.5,0)[below, scale=.6]{$3$};
    \node at (.5,1)[above, scale=.6]{$2$};
    \node at (2,.5)[right, scale=.6]{$3$};
    \node at (1.5,0)[below, scale=.6]{$1$};
    \node at (1,1.5)[left, scale=.6]{$2$};
    \node at (2,1.5)[right, scale=.6]{$1$};
    \node at (1,2.5)[left, scale=.6]{$3$};
    \node at (2,2.5)[right, scale=.6]{$2$};
    \end{tikzpicture}$};

\node at (0.3,-5){$\begin{tikzpicture}[scale=.5]
    \draw (0,0)--(1,0)--(1,1)--(0,1)--(0,0);
    \draw[densely dotted] (0,1)--(1,0);
    \node at (.5,.5){$2$};
    \draw[line width=1.2, dashed] (0,0)--(1,0);
    \node at (.5,1)[above, scale=.6]{$3$};
    \node at (0,.5)[left, scale=.6]{$1$};
    \end{tikzpicture}$};

\draw[-open triangle 45] (1,-5) -- node[rotate=0,above, scale=.7] {glue} (2,-5);

\node at (3,-5){$\calg^{\Delta}_{\gamma}$};

\node at (4.2,-5){$\begin{tikzpicture}[scale=.5]
    \draw (0,0)--(1,0)--(1,1)--(0,1)--(0,0) (1,0)--(2,0)--(2,1)--(1,1) (1,1)--(1,2)--(2,2)--(2,1) (1,2)--(1,3)--(2,3)--(2,2) (1,3)--(1,4)--(2,4)--(2,3);
    \draw[densely dotted] (0,1)--(1,0) (1,1)--(2,0) (1,2)--(2,1) (1,3)--(2,2) (1,4)--(2,3);
    \node at (.5,.5){$1$};
    \node at (1.5,.5){$2$};
    \node at (1.5,1.5){$3$};
    \node at (1.5,2.5){$1$};
    \node at (1.5,3.5){$2$};
    \node at (.5,0)[below, scale=.6]{$3$};
    \node at (.5,1)[above, scale=.6]{$2$};
    \node at (2,.5)[right, scale=.6]{$3$};
    \node at (1.5,0)[below, scale=.6]{$1$};
    \node at (1,1.5)[left, scale=.6]{$2$};
    \node at (2,1.5)[right, scale=.6]{$1$};
    \node at (1,2.5)[left, scale=.6]{$3$};
    \node at (2,2.5)[right, scale=.6]{$2$};
    \node at (1,3.5)[left, scale=.6]{$1$};
    \node at (1.5,4)[above, scale=.6]{$3$};
    \end{tikzpicture}$};

\draw[-open triangle 45] (5,-5) -- node[rotate=0,above, scale=.7] {} (6,-5);

\node at (7,-5){$\calg_{\gamma}$};

\node at (8,-5){$\begin{tikzpicture}[scale=.5]
    \draw (0,0)--(1,0)--(1,1)--(0,1)--(0,0) (1,0)--(2,0)--(2,1)--(1,1) (1,1)--(1,2)--(2,2)--(2,1) (1,2)--(1,3)--(2,3)--(2,2) (1,3)--(1,4)--(2,4)--(2,3);
    \node at (.5,.5){$1$};
    \node at (1.5,.5){$2$};
    \node at (1.5,1.5){$3$};
    \node at (1.5,2.5){$1$};
    \node at (1.5,3.5){$2$};
    \node at (.5,0)[below, scale=.6]{$3$};
    \node at (.5,1)[above, scale=.6]{$2$};
    \node at (2,.5)[right, scale=.6]{$3$};
    \node at (1.5,0)[below, scale=.6]{$1$};
    \node at (1,1.5)[left, scale=.6]{$2$};
    \node at (2,1.5)[right, scale=.6]{$1$};
    \node at (1,2.5)[left, scale=.6]{$3$};
    \node at (2,2.5)[right, scale=.6]{$2$};
    \node at (1,3.5)[left, scale=.6]{$1$};
    \node at (1.5,4)[above, scale=.6]{$3$};
    \end{tikzpicture}$};

\end{tikzpicture}
\caption{Snake graph associated to the arc $\gamma$.} \label{ExampleSnakeGraph}
\end{figure}
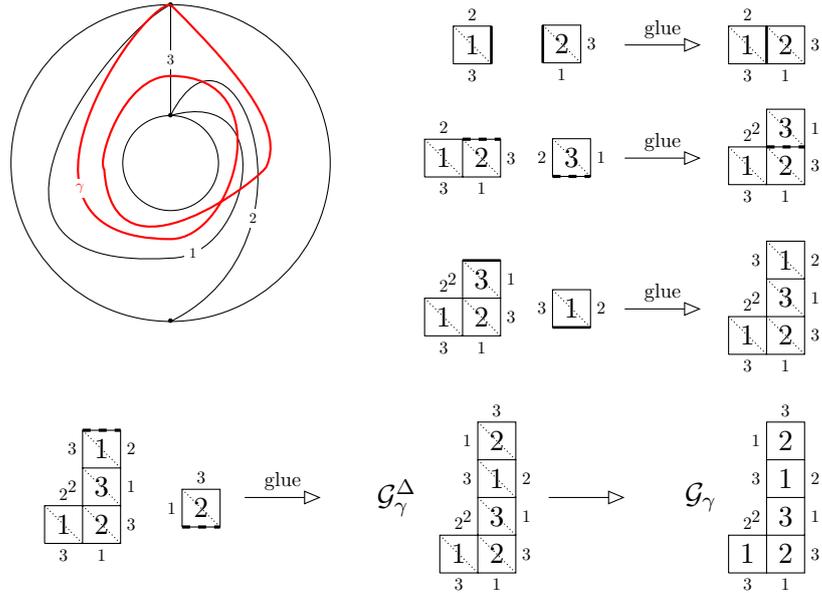

The labels on the edges of a snake graph, given by the corresponding arcs in the triangulation, 
are called \emph{weights}. Sometimes  snake graphs with weights are referred to as \emph{labelled snake graphs}. 
See Figure~\ref{ExampleSnakeGraph} for an example of a labelled snake graph associated to an arc. 

The $d-1$ edges corresponding to the arcs $\sigma_1, \ldots, \sigma_{d-1}$ which are contained in two tiles are called \emph{interior edges of } $\mathcal{G}_\gamma$. Denote this set 
by ${\rm{Int}}({\calg}_\gamma)$. The edges of  $\mathcal{G}_\gamma$ not in  ${\rm{Int}}(\calg_\gamma)$ are called \emph{boundary edges}. We define a subgraph $\mathcal{G}_\gamma[ i, i+t]$, for $1 \leq i \leq d$ and for  $0 \leq t \leq d-i$, to be 
the subgraph of $\mathcal{G}_\gamma$ consisting of the tiles $(G_i, \ldots, G_{i+t})$.

Let $_{SW} \mathcal{G}_\gamma$ (resp.  $ \mathcal{G}^{NE}_\gamma$) be the set containing the 2 elements corresponding to the south and west edge of $G_1$ (resp. the north and east edge of 
$G_d$). Define $\mathcal{G}_\gamma \setminus \pred (\sigma_i) = \mathcal{G}_\gamma [i+1, d]$. If $e$ is an edge in $ \mathcal{G}_\gamma^{NE}$ then 
$\mathcal{G}_\gamma \setminus \pred (e) = \{e\}$. 
Analogously, let  $\mathcal{G}_\gamma \setminus \suc (\sigma_i) = \mathcal{G}_\gamma [1, i]$. If $e$ is an edge in $_{SW} \mathcal{G}_\gamma$ then 
$\mathcal{G}_\gamma \setminus \suc (e) = \{e\}$.  

If all tiles of a snake graph $\mathcal{G}_\gamma$ are in a row or a column, we call $\mathcal{G}_\gamma$ \emph{straight} and we call it \emph{zigzag} if no three consecutive tiles are 
straight. 

Note that there is a notion of abstract snake graphs  as  combinatorial objects introduced in \cite{CS}. However, all snake graphs we consider here are snake graphs associated to arcs 
in triangulated surfaces as introduced above. In general, we will use the notation $\mathcal{G}$ if we do not need to refer to the associated arc or if the arc is clear from the context.

\begin{Definition}
A \emph{sign function} on a snake graph $\mathcal{G}$ is a function $f$ from the set of interior edges $\{\sigma_1, \ldots, \sigma_{d-1}\}$ to the set $\{+,-\}$ such that 

(1) if three consecutive tiles $G_{i-1}, G_i, G_{i+1}$ form a straight piece then $f(\sigma_{i-1}) = - f(\sigma_i)$, 

(2) if three   consecutive tiles $G_{i-1}, G_i, G_{i+1}$ form a zigzag piece then $f(\sigma_{i-1}) =  f(\sigma_i)$.

Extend the sign function to all edges of $\calg$ by the following rule: opposite edges have opposite signs and the south side and east side of each tile have the same sign as do the north and west side of each tile. 
\end{Definition}

Note that for every snake graph there are two sign functions, $f$ and $f'$ such that $f(\sigma_i) = - f'( \sigma_i)$, for each $ 1 \leq i \leq d-1$. 

\begin{figure}[ht]
\begin{tikzpicture}[trans/.style={thick,->,shorten >=2pt,shorten <=2pt,>=stealth}]
\node at (0,0){$\begin{tikzpicture}[scale=.7]
    
    \node[scale=.4](a) at (0,7) [label={[label distance=-.05cm, scale=.7]90:}, scale=.7] {$\bullet$};
    \node(a') at (-1.5,7.5) { };
    \node(a'') at (1.5,3){ };
    \node(bb) at (1.5,7.5) { };
    \node(b'') at (-1.5,3) { };
    \node[scale=.4](f) at (2,6) [label={[label distance=-.15cm, scale=.7]10:}, scale=.7] {$\bullet$};
    \node[scale=.4](ff) at (-2,6) [label={[label distance=-.15cm, scale=.7]170:}, scale=.7] {$\bullet$};
    \node[scale=.4](g) at (2,4) [label={[label distance=-.15cm, scale=.7]-10:}, scale=.7] {$\bullet$};
    \node[scale=.4](gg) at (-2,4) [label={[label distance=-.15cm, scale=.7]-170:}, scale=.7] {$\bullet$};
    \node[scale=.4](h) at (0,3) [label={[label distance=-.05cm, scale=.7]270:}, scale=.7] {$\bullet$};
    
    
    \draw (a.center) -- (f.center);
    \draw (a.center) -- (ff.center);
    \path (ff.center) edge (f.center);
    \path (gg.center) edge (g.center);
    \draw [dotted, line width=.8] (ff.center)--(gg.center) (f.center)--(g.center);
    \draw (gg.center)--(h.center) (g.center)--(h.center);
    
    \path[color=red, line width=.7] (a'.center) edge node [pos=.7, rotate=115, scale=.7]{$<$} ( $(a')!0.4!(a'') $ );
    \draw[color=red, line width=.7] (a'.center) -- ( $(a')!0.4!(a'') $ );
    \path[color=red, line width=.7] ( $(a')!0.4!(a'') $ )[dotted] edge node [pos=.6, fill=white,outer sep=1mm,scale=.6, bend right]{$\gamma_M$} ( $(a')!0.7!(a'') $ ) ;
    \draw[color=red, line width=.7] ( $(a')!0.7!(a'') $ ) -- (a''.center);
    
    \path[color=tealblue, line width=.7] (bb.center) edge node [pos=.65, rotate=65, scale=.7]{$<$} ( $(bb)!0.4!(b'') $ );
    \draw[color=tealblue, line width=.7] (bb.center) -- ( $(bb)!0.4!(b'') $ );
    \path[color=tealblue, line width=.7] ( $(bb)!0.4!(b'') $ )[dotted] edge node [pos=.65, fill=white,outer sep=1mm,scale=.6, bend right]{$\gamma_N$} ( $(bb)!0.7!(b'') $ ) ;
    \draw[color=tealblue, line width=.7] ( $(bb)!0.7!(b'') $ ) -- (b''.center);
    %
    \end{tikzpicture}$};

\node at (4,0){$\begin{tikzpicture}[trans/.style={thick,->,shorten >=2pt,shorten <=2pt,>=stealth}]
    \node at (3.7,0)[rotate=90]{$\begin{tikzpicture}[scale=.7]
        \node[scale=.4] (a) at (1,1.7) [label={[label distance=0cm, scale=.7,rotate=-90]170:}, scale=.7] {$\bullet$};
        \node[scale=.4] (b) at (0,-.8) [label={[label distance=-.1cm, scale=.7,rotate=-90]-10:}, scale=.7] {$\bullet$}; 
        \node[scale=.4] (c) at (2,-.8) [label={[label distance=-.1cm, scale=.7,rotate=-90]10:}, scale=.7] {$\bullet$};
        \draw (a.center)--(b.center)--(c.center)--(a.center);
        
        \path[color=red] (b.center) edge node[pos=.5, scale=.7, rotate=35]{$<$} (2.5,1.5);
        \path[color=red] (b.center) edge node[pos=.74, fill=white,outer sep=1mm, scale=.6, rotate=-90]{$\gamma_1$} (2.5,1.5);
        \draw[color=red, dotted] (2.5,1.5) -- ( 3,2 );
        \path[color=tealblue] (c.center) edge node[pos=.5, scale=.7, rotate=-40]{$<$} (-.5,1.5);
        \path[color=tealblue] (c.center) edge node[pos=.74, fill=white,outer sep=1mm, scale=.6, rotate=-90]{$\gamma_2$} (-.5,1.5);
        \draw[color=tealblue, dotted] (-.5,1.5) -- ( -1,2 );
        %
        \end{tikzpicture}$};
    \end{tikzpicture}$};

\node at (8,0){\begin{tikzpicture}[trans/.style={thick,->,shorten >=2pt,shorten <=2pt,>=stealth}]
    \node at (7,-.5){\begin{tikzpicture}[scale=.7]
        \node[scale=.4] (a) at (1,2)  [label={[label distance=-.1cm, scale=.7]10:}, scale=.7]{$\bullet$};
        \node[scale=.4] (b) at (0,0)  [label={[label distance=-.1cm, scale=.7]190:}, scale=.7]{$\bullet$}; 
        \node[scale=.4] (c) at (2,0)  [label={[label distance=-.1cm, scale=.7]-30:}, scale=.7]{$\bullet$};
        \node (d) at (.4,1.7)  { };
        \node (e) at (2.5, 1)  { };
        \node (f) at (1, -.5)  { };
        \node (g) at (1.2,-.6) { };
        \node (h) at (-.4,1.7) {};
        \draw (a.center)--(b.center)--(c.center)--(a.center);
        
        \path[line width=.7, color=red] (d.center) edge node[pos=.4, scale=.7, rotate=-65]{$>$} (g.center);
        \draw[color=red, dotted] (.15,2.5) -- (d.center);
        \draw[color=red, dotted] (g.center) -- (1.4,-1.3);
        \path[line width=.7, color=red] (d.center) edge node[pos=0, fill=white,outer sep=1mm, scale=.6]{$\gamma_1$} (g.center);
        
        \path[line width=.7, color=tealblue] (b.center) edge node[pos=.55, scale=.7, rotate=15]{$>$} (e.center);
        \path[line width=.7, color=tealblue] (b.center) edge node[pos=.88, fill=white,outer sep=1mm, scale=.6]{$\gamma_2$} (e.center);
        \draw[color=tealblue, dotted] (e.center) -- (3.4,1.3);
        \end{tikzpicture}};
    \end{tikzpicture}};
\end{tikzpicture}
\caption{The leftmost figure corresponds to an overlap crossing in terms of snake graphs and a module crossing in terms of string modules, the middle figure corresponds to grafting with $s=d$ in terms of snake graphs and an arrow crossing in terms of string modules, and the rightmost figure corresponds to grafting with $s \neq d$ in terms of snake graphs and a 3-cycle crossing in terms of string modules.} \label{ThreeTypes}
\end{figure}
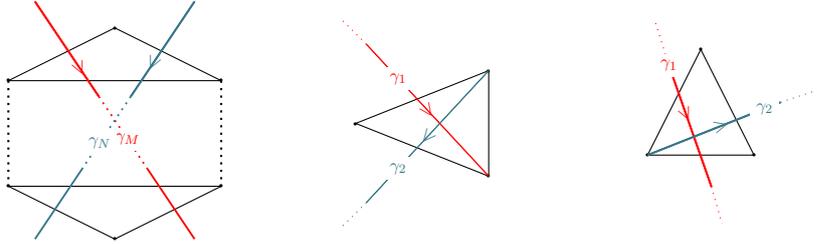

A crossing of two arcs $\gamma_1, \gamma_2$ has an interpretation in terms of the associated snake graphs as given in \cite{CS} and as
further explored in \cite{CS2}.
Depending on the triangulation and the arcs, there are three different ways in which the two arcs can cross, see Figure~\ref{ThreeTypes}. 
In the first case,  the arcs cross in what we refer to as an  \emph{overlap} since both arcs $\gamma_1$ and $ \gamma_2$ cross at least one  common arc in the triangulation and the crossing can be moved up to homotopy to any triangle adjoining an arc of the triangulation crossed by both $\gamma_1$ and $\gamma_2$.  In the last two cases the crossing occurs in a single triangle and cannot be moved outside of this triangle by homotopy. We say that the arcs {\it cross in a triangle}.

In terms of snake graphs and using the terminology of \cite{CS2} these crossings correspond 
to  a crossing 
overlap, grafting with $s=d$ and grafting with $s \neq d$. Note that the integer
$d$ corresponds to the number of tiles of the snake graph $\calg_1$ corresponding to $\gamma_1$ and $s$ is some integer $1 \leq s \leq d$ as defined in \cite[Section 3.3]{CS2} denoting the position at which the `grafting' takes place. We will later see in Section~{\ref{CrossingModulesSection}}  that these correspond to three types of module crossings, 
namely crossing in a module, arrow crossing,
and 3-cycle crossing, respectively.

We start by defining an overlap of two snake graphs and a self-overlap of a snake graph. 

\begin{Definition}\cite[Section 2.5]{CS2} 
Let $\calg_1=(\gi 12d)$ and $\calg_2=(\gii 12{d'})$ be two snake graphs such that there exists two embeddings of graphs  $i'_1: \calg \to \calg_1$  and $i'_2: \calg \to \calg_2$ such that,  for $j =1,2$, the image $i_j(\calg)$ is either identical to $\calg$,  a 180$^\circ$ rotation of $\calg$, or a reflection of $\calg$ at one of the lines $y = x$ or $y =x$. In particular, the south west vertex of the first tile of $\calg$ has to be mapped to the south west vertex of the first tile in $i_j(\calg)$ or to the north west vertex of the last tile in $i_j(\calg)$.  Moreover, we require that  $i_1$ and $i_2$ are maximal in the following sense:

(1) If $\calg$ has at least two tiles and if there exists a snake graph $\calg'$ with two embeddings $i'_1: \calg' \to \calg_1,$ $i'_2: \calg' \to \calg_2$ such that $i_1(\calg) \subseteq i'_1(\calg')$ and $i_2(\calg) \subseteq i'_2(\calg')$ then $i_1(\calg) = i'_1(\calg')$ and $i_2(\calg) = i'_2(\calg').$

(2) $\calg$ is a snake graph consisting of at least one tile.

If the above the hold then we say  that $\calg_1$ and $\calg_2$ have an \emph{overlap} at $\calg$.  
\end{Definition}

In the case of a self-overlap we might have $i_1(\calg) \cap i_2(\calg) \neq \emptyset$. 

Two snake graphs  might have several overlaps with respect to the same and different snake sub-graphs.  

\begin{Definition}\cite[{Definition 2.4}]{CS2}\label{DefinitionCrossingOverlap}
Let $\calg_1 = (G_1, \ldots, G_d)$ and $\calg_2 = (G'_1, \ldots, G'_{d'})$ be two snake graphs with overlap $\calg$ and
embeddings $i_1(\calg) = \calg_1[s,t]$ and 
$i_2(\calg) = \calg_2[s',t']$. Let $(\sigma_1, \ldots, \sigma_{d-1})$ (respectively $(\sigma'_1, \ldots, \sigma'_{d'-1})$ ) 
be the interior edges of $\calg_1$ (respectively $ \calg_2$) and 
let $f$ be a sign function on $\calg$. Then $f$ induces sign functions $f_1$ on $\calg_1$ and $f_2$ on $\calg_2$. We say that $\calg_1$ \emph{and} 
$\calg_2$ \emph{cross in $\calg$} if one of the following conditions hold. 

(1) $f_1(\sigma_{s-1}) = - f_1 (\sigma_t)$ if $s>1$ and $t <d$ or  $f_2(\sigma'_{s'-1}) = - f_2 (\sigma'_{t'})$ if $s'>1$ and $t' <d',$

(2) $f_1 (\sigma_t) = f_2(\sigma'_{s'-1})$ if $s=1$, $t<d$, $s' >1$, and $t' = d'$ or $f_1(\sigma_{s-1}) = f_2 (\sigma'_{t'})$ if 
$s>1$, $t=d$, $s'=1$, and $t' <d'$.

We call such an overlap a \emph{crossing }.
\end{Definition}

We have a similar definition of a self-crossing overlap. 

\begin{Definition}\cite[{Definition 2.6}]{CS2}\label{DefinitionSelfCrossingOverlap}
Let $\calg_1 = (G_1, \ldots, G_d)$ be a snake graph with self-overlap $\calg$ and
embeddings $i_1(\calg) = \calg_1[s,t]$ and 
$i_2(\calg) = \calg_1[s',t']$. Let $(\sigma_1, \ldots, \sigma_{d-1})$ 
be the interior edges of $\calg_1$  and 
let $f$ be a sign function on $\calg_1$.  We say that $\calg_1$ has a {\em self-crossing in $\calg$} if one of the following conditions hold. 

(1) $f(\sigma_{s-1}) = - f (\sigma_t)$   or $f(\sigma_{s'-1}) = - f (\sigma_{t'})$ if  $t' <d,$

(2) $f (\sigma_t) = f(\sigma_{s'-1})$.

We call such an overlap a \emph{self-crossing overlap}.
\end{Definition}

We will now see that a crossing or self-crossing overlap corresponds to a crossing or a self-crossing of arcs. 

\begin{Theorem}\cite[Theorem 6.1]{CS2}\label{Overlap=Crossing} Let $\gamma_1, \gamma_2$ be (generalized) arcs and $\calg_1, \calg_2$ their corresponding snake graphs.
\begin{itemize}
\item[(1)] $\gamma_1, \gamma_2$ cross with a nonempty local overlap $(\tau_{i_s}, \cdots, \tau_{i_t})=(\tau_{i'_{s'}}, \cdots, \tau_{i'_{t'}})$ if and only if $\calg_1, \calg_2$ cross in $\calg_1[s,t] \cong \calg_2[s',t'].$
\item[(2)] $\gamma_1$ crosses itself  with a non-empty local overlap $(\tau_{i_s}, \cdots, \tau_{i_t})=(\tau_{i'_{s'}}, \cdots, \tau_{i'_{t'}})$  if and only if  
$\calg_1$ has a self-crossing overlap $\calg_1[s,t] \cong \calg_1[s',t'].$
\end{itemize}
\end{Theorem}

In the proof of Theorem~\ref{Triangles} we  consider self-crossings of an arc. However, the way we  treat a self-crossing of an arc is by replacing the 
one arc by two copies of the same arc. We then smooth crossings of the `two' arcs as opposed to \cite{CS2} where self-crossings of the single arc are smoothed.  This can be 
done because of the following more general statement. 

\sloppy

\begin{Proposition}\label{self-crossing equals crossing overlap}
Let $\calg_\gamma$ be a snake graph with a self-overlap $\calg$. Then $\calg$ is a self-crossing overlap of  $\calg_\gamma$  if and only if $\calg$ is a crossing overlap
of two copies of $\calg_\gamma$. 
\end{Proposition}

We note that here we consider the two copies of $\calg_\gamma$ as `distinct' snake graphs.  

{\it Proof:}
Suppose $f$ is a sign function on $\calg_\gamma = (G_1, \ldots, G_s, \ldots, G_{s'}, \ldots, G_d)$ and that $\calg_\gamma$ has a self-crossing overlap $\calg$. Denote by 
$(\sigma_1, \ldots, \sigma_{d-1})$ the interior edges of 
$\calg_\gamma$. 
Then by Definition~\ref{DefinitionSelfCrossingOverlap} 
there exist two embeddings $i_1(\calg) = (G_s, \ldots, G_t)$ and $i_2(\calg) = (G_{s'}, \ldots, G_{t'})$ where $s < s'$ and $f(\sigma_t) = f (\sigma_{s' -1})$.  Now consider
two copies of $\calg_\gamma$, denote them by $\calg_\gamma$, $\calg'_\gamma$. Then we can canonically embed $(G_s, \ldots, G_t)$ into $\calg_\gamma$ and $(G_{s'}, \ldots, G_{t'})$
into $\calg'_\gamma$. Then by \cite{CS} this is an overlap for $\calg_\gamma$ and $\calg'_\gamma$ and the sign conditions for a crossing in the overlap $\calg$ is satisfied by the sign function
$f$ on 
$\calg_\gamma$ and on $\calg'_\gamma$.  

The converse immediately follows from a similar argument. 
\hfill $\Box$ 

It is clear that reversing the roles of $\calg_\gamma$ and $\calg'_\gamma$ gives a similar result for the second crossing of $\calg_\gamma$ and $\calg'_\gamma$ and it 
immediately implies that the number of self-crossings of an arc $\gamma$ with itself times two is equal to the number of crossings of two copies of $\gamma$. 

\sloppy

The smoothing of a crossing of two arcs $\gamma_1$ and $\gamma_2$ such that the associated snake graphs $\calg_1$ and $\calg_2$ cross in an overlap is 
called , in terms of snake graphs,  the \emph{resolution of the overlap} \cite{CS2}. 
It immediately follows from the definitions above and below that given a self-crossing of an arc $\gamma$ with associated snake graph $\calg_\gamma$, the two corresponding crossings of two copies of $\calg_\gamma$ coincide and give rise to the same
resolution. 

\sloppy

We now define four snake graphs $\calg_3$, $\calg_4$, $\calg_5$ and $\calg_6$.  For $\calg[i,j] = (G_i, \ldots, G_j)$, define 
$\overline{\calg} [j, i] = (G_j, \ldots, G_i)$. 
In order to define $\calg_5$ and $\calg_6$, we introduce the notation 
$\calg_5'=\calg_1[1,s-1] \cup \ocalg_2[s'-1,1]$ where we glue the two subgraphs
along the edge with weight $\sigma_s$ and  $\calg_6'=\ocalg_2[d',t'+1] \cup \calg_1[t+1,d]$
where we glue the two subgraphs along the edge with weight $\sigma_t$.

Let $f_5$ be a sign function on $\calg_5'$ and $f_6 $ a sign function on $\calg_6'$.
\sloppy

We then define four snake graphs as follows. 

$ \calg_3 =\calg_1[1,t] \cup \calg_2[t'+1, d']$ where the gluing of the two subgraphs is induced by  by the embedding $i_2$ of $\calg$ in $\calg_2; $

$ \calg_4 =\calg_2[1,t'] \cup \calg_1[t+1,d]$ where the gluing of the two subgraphs is induced by by the embedding $i_1$ of $\calg$ in $\calg_1;$
\sloppy \begin{align*}
\calg_5 &=
\begin{cases} \calg_5'&  \parbox[t]{.55\textwidth}{ if $s>1$, $s'>1$ ;}\\
\calg_5'\setminus \suc (\sigma) &     \parbox[t]{.65\textwidth}{ if $s'=1$ where $\sigma$ is the last edge in $\Int(\calg_5') \cup {}_{SW}\calg_5'$ such  that $f_5(\sigma)=f_5(\sigma_{s-1})$; }  \\
\calg_5'\setminus \pred (\sigma) &     \parbox[t]{.65\textwidth}{ if $s=1$ where $\sigma$ is the first edge in $\Int(\calg_5') \cup \calg_5'^{N\!E}$ such  that $f_5(\sigma)=f_5(\sigma'_{s-1})$; }  \\ 
\end{cases}\\ 
\calg_6&=
\begin{cases}
\calg_6'&  \parbox[t]{.55\textwidth}{ if $t<d$, $t' <d'$;}\\
\calg_6'\setminus \suc (\sigma) &     \parbox[t]{.65\textwidth}{ if $t=d$, where $\sigma$  is the last edge in $\Int(\calg_6') \cup {}_{S\!W}\calg_6'$ such  that $f_6(\sigma)=f_6(\sigma'_{t'})$; }   \\
\calg_6'\setminus \pred (\sigma) &     \parbox[t]{.65\textwidth}{ if $t'=d'$, where $\sigma$  is the first edge in $\Int(\calg_6') \cup \calg_6'^{N\!E}$ such  that $f_6(\sigma)=f_6(\sigma_{t})$. } 
\end{cases}\\
\end{align*}

In the above definition we write ${\rm{Res}}_\calg(\calg_1, \calg_2) = (\calg_3 \sqcup \calg_4, \calg_5 \sqcup \calg_6)$ for the resolution of the crossing
of $\calg_1$ and $\calg_2$.

\begin{Theorem}\cite[Theorem 6.2]{CS2}\label{Overlap=Crossing} 
Let $\gamma_1, \gamma_2$ be (generalized) arcs and $\calg_1, \calg_2$ their corresponding snake graphs.
The snake graphs of the four arcs obtained by smoothing the crossing of $\gamma_1$ and $\gamma_2$ in the overlap are given by the resolution $\re 12$ of the crossing of the snake graphs $\calg_1$ and $\calg_2$ at the overlap $\calg.$
\end{Theorem}

We now consider a particular crossing of   two arcs $\gamma_1$ and $\gamma_2$   such that this crossing does not correspond to a crossing overlap 
in the associated snake graphs $\calg_1$ and $\calg_2$. This situation occurs exactly if the crossing of $\gamma_1$ and $\gamma_2$ occurs in what we call in section 3 \textit{a crossing in an arrow or in a 3-cycle} in 
$(S, \M, T)$. In particular, one (or both) of the arcs will have at least one endpoint coinciding with a vertex in that triangle. 
Following \cite{CS2} there are two cases to consider. 

Let  $\calg_1=(\gi 12d)$ and $\calg_2=(\gii 12{d'})$ be two snake graphs such that $G_s \neq G'_1$ for some $1 \leq s \leq d$ and let $f_1$ be a sign function on $\calg_1$.
Let $\e$ be the unique common edge in $\calg^{N\!E}_s$ and  ${}_{S\!W}\calg'_1$ if it exists. 
Let $f_2$ be  a sign function on $ \calg_2$ such that $f_2(\e)=f_1(\e).$ 
Then define four snake graphs as follows.

\emph{Case 1.} Suppose $s=d.$  

\begin{align*}
\calg_3=&\calg_1 \cup \calg_2 \mbox{ where the  two subgraphs are glued along the edge $\e$;} &\\
\calg_4=& \{ \e \};
\\
\calg_5 =&\,\calg_1 \setminus \suc(\sigma) ,\  \parbox[t]{.75\textwidth}{where $\sigma \in\calgSW_1\cup\,\Int(\calg_1)$  is the last edge such that $f_1(\sigma)=f_1(\delta)$;}
\\
\calg_6=&\,\calg_2\setminus \pred(\sigma),\ \mbox{where $\sigma \in \Int(\calg_2)\cup \calgNE_2$ is the first edge such that $f_2(\sigma)=f_2(\delta).$}&
\end{align*}

\emph{Case 2.}
Suppose that $s \neq d.$

\begin{align*}
\calg_3=&\,\calg_1[1,s] \cup \calg_2 , \mbox{ where the  two subgraphs are glued along the edge $\e$}; \\
\calg_4 =&\,\calg_1 \setminus \pred(\sigma),\   \parbox[t]{.65\textwidth}{where $\sigma \in\Int(\calg_1[s\!+\!1,d])\cup\calgNE_1$  is the first edge such that $f_1(\sigma)=f_1(\delta)$;}
\\
\calg_5 =&\,\calg_1 \setminus \suc(\sigma) ,\  \parbox[t]{.75\textwidth}{where $\sigma \in\calgSW_1\cup\,\Int(\calg_1[1,s])$  is the last edge such that $f_1(\sigma)=f_1(\delta)$;}
\\
\calg_6=&\,\ocalg_2[d',1] \cup \calg_1 [s+1,d] ,\ \mbox{where the two subgraphs are glued along the edge $\sigma_s$.} &
\\
\end{align*}

In the above definitions we write $\gt s{\delta}12  = (\calg_3 \sqcup \calg_4, \calg_5 \sqcup \calg_6)$ and we call it grafting of $\calg_2$ on $\calg_1$
at $s$. 

Similarly, we consider a particular self-crossing of an arc $\gamma$ such that this crossing does not correspond to a crossing in an overlap in the associated snake graph $\calg_\gamma$. 
As in the case of two distinct arcs this situation occurs exactly if the self-crossing of $\gamma$ does not have any overlap.   It is immediate that in terms of 
snake graphs detecting a self-crossing corresponding to a self-crossing of $\gamma$ in a triangle is equivalent to considering the arc as two distinct arcs and two
distinct snake graphs and detecting a crossing of two arcs in a triangle in terms of the two snake graphs. 

\begin{Theorem}\cite[Theorem 6.4]{CS2}\label{Grafting=Crossing}
Let $\gamma_1$ and $\gamma_2$ be two arcs which cross in a triangle $\Delta$ with an empty  overlap, and let $\calg_1$ and $\calg_2$ be the corresponding snake graphs. Assume the orientation of $\gamma_2$ is such that $\Delta$ is the first triangle $\gamma_2$ meets. Then the snake graphs of the four arcs obtained by smoothing the crossing of $\gamma_1$ and $\gamma_2$ in $\Delta$ are given by the resolution $\gt s{\delta}12$ of the grafting of $\calg_2$ on $\calg_1$ at $s,$ where $0 \leq s \leq d$ is such that $\Delta_s = \Delta$ and if $s=0$ or $s=d$ then $\delta$ is the unique side of $\Delta$ that is not crossed by either $\gamma_1$ or $\gamma_2.$
\end{Theorem}

Given a triangulation  of $(S,\M)$, we now establish a one-to-one correspondence between the set of snake graphs with a sign function (excluding
snake graphs consisting of a single edge) and the set of 
string modules of the associated Jacobian algebra $J(Q,W)$. 
Let $M(w_\gamma)$ be the string module
corresponding to  an arc $\gamma$ and let $\calg_\gamma$ be the associated snake graph. Then the arrows and their formal inverses uniquely 
define a sign function $f_\gamma$ on 
$\calg_\gamma$. Namely, let $w_\gamma  =  \varepsilon_1 \ldots \varepsilon_{d-1}$  and let $(\sigma_1, \ldots, \sigma_{d-1})$
be the interior edges of $\calg_\gamma$. Define a sign function $f_\gamma$ on $\calg_\gamma$ by setting
$f_\gamma (\sigma_i) = \phi(\varepsilon_i)$ where $\phi(\varepsilon_i) = +$ if $\varepsilon_i \in Q_1$ and $\phi(\varepsilon_i) = -$ if 
$\varepsilon_i^{-1} \in Q_1$, for $1 \leq i \leq d-1$. 

Set ${\mathcal{R}} = \{ \; (\calg_\gamma, f_\gamma) \;| \; \gamma \mbox{ an arc in } (S,\M,T) \mbox{ and such that $\gamma$ is not in $T$} \, \}$. The following result is immediate. 

\begin{Proposition}\label{SnakeStrings}
There is a bijection between the set of strings $\S$ over $J(Q,W)$ and the set ${\mathcal{R}}$ given by the map that associates $(\calg_\gamma, f_\gamma)$ to the string $w_\gamma$ for every 
arc $\gamma$ in $(S,\M,T)$, $\gamma \notin T$.
\end{Proposition}

Such a correspondence in the setting of a triangulation of the once-punctured
torus has also been noted  in \cite{Ricke}.




\section{Extensions for the Jacobian algebra}\label{CrossingModulesSection}

Let $J(Q,W)$
be the Jacobian algebra  associated to a triangulation of a marked surface $(S, \mathcal{M})$ with all marked points in the boundary and such that every boundary component contains at least one marked point.

In this section we interpret  the crossing of arcs in terms of the corresponding string modules over $J(Q,W)$.   We do this by using the 
characterization of crossing
arcs in terms of snake graphs introduced in Section~\ref{SnakeSection} and the snake graph and string module correspondence given in Proposition~\ref{SnakeStrings}.

\subsection{Crossing String Modules}

Given two (not necessarily distinct) arcs in a surface $(S,\M)$, recall that there are three types of configurations in which these arcs can cross, see Figure~\ref{ThreeTypes}. 

Each of these crossings gives rise to a different structure of the corresponding string modules which leads to Definition~\ref{DefinitionMcrossesN}. 

We use the notation $\pred(\alpha)$ for the substring preceding an arrow or an inverse arrow $\alpha$ in a string $w$ and similarly we use the notation $\suc(\alpha)$ for the substring succeeding an arrow or an inverse arrow  $\alpha$ in a string $w.$

\begin{Definition}\label{DefinitionMcrossesN}
{\rm We say that a string module $M$ crosses a string module $N$ if there exist strings $w_M$ and $w_N$ such that $M \simeq M(w_M)$ and $N \simeq M(w_N)$ and if one of the following three 
    conditions hold

    (1)  there exists a string $w \in \S$, possibly consisting of a  single vertex  only, and  such that $w_M$ and $w_N$ do not both start at $s(w)$ or do not
    both end at $t(w)$ and if $$w_M = \Pred(\alpha) \stackrel{\alpha}{\longrightarrow} w  \stackrel{\beta}{\longleftarrow} \Succ(\beta)
    \;\;\;\; \mbox{               and              } \;\;\;\;
    w_N= \Pred(\varepsilon)  \stackrel{\varepsilon}{\longleftarrow} w \stackrel{\delta}{\longrightarrow} \Succ(\delta)$$
    where $\alpha, \beta, \varepsilon, \delta$ are arrows in $Q_1$.
    If $\alpha$ (resp. $\varepsilon$) doesn't exist then $w_M$ (resp. $w_N$) starts with $w$ and if $\beta$ (resp. $\delta$) doesn't exist then $w_M$ (resp. $w_M$) ends with $w$;
    
    (2) there exists an arrow $\alpha$ in $Q_1$ such that $w_M \stackrel{\alpha}{\longrightarrow} w_N \in \S;$ 
    
    (3) if there exists   a 3-cycle $ \xymatrix{
        a \ar[r]|\alpha           & b \ar[d]|\delta \\
        & \ar[ul] | \beta c }$
    in $Q$ with $\alpha \delta, \delta \beta,  \beta \alpha \in I$ and 
    such that $\alpha$ is in $w_M$ and $s(w_N) = c$ and $w_N$ does not start or end with $\delta$ 
    or $\beta$ nor their inverses. 
    
    Moreover, if (1) holds we say that $M$ crosses $N$ in \emph{a module}, if (2) holds we say that 
    $M$ crosses $N$ in \emph{an arrow} and if (3) holds we say that 
    $M$ crosses $N$  in \emph{a 3-cycle}.  Moreover, if $M=N$, we say $M$ has \emph{a self-crossing} in \emph{a module, an arrow, and a 3-cycle, respectively.}

    Furthermore,  if in (1) the arc $\gamma_{M(w)}$ is a self-crossing arc, we say $M_1$ crosses $M_2$ in a {\it self-crossing module $M(w)$},
    and if $M=N,$ we say that $M$ has a {\it self-crossing in a self-crossing module $M(w)$}.   In accordance with the terminology of snake graphs, we say that the corresponding arcs $\gamma_1$ and $\gamma_2$ have a {\it crossing in a self-crossing overlap}. If  the arc $\gamma_{M(w)}$ has no self-crossing then we say that the respective crossing is in a {\it non self-crossing module $M(w)$}. 
}
\end{Definition}

\begin{Remark} (1) We note that there is a direction in the crossings of modules. Namely, for a fixed crossing, we say that $M$ crosses $N$ or $N$ crosses $M$. We will see in Theorem 3.7 that for a particular crossing, $M$ crossing $N$ will give rise to an element in $\Ext^1_{\mathcal J}(M,N)$ if the crossing is in a module or in an arrow whereas there will be no element in $ \Ext^1_{\mathcal J}(N,M)$ corresponding to this crossing.

(2)It is possible for a string module $M$ to cross a string module $N$ simultaneously in 
a module, in an arrow and in a 3-cycle. 
In Section~\ref{ExampleSection} we give an example of all three types of crossings of modules as defined above. 

(3) Crossings in modules as described in~\ref{DefinitionMcrossesN}(1) have also been considered in \cite{BZ, GLS2, ZZZ}. It also immediately follows from \cite{CB} that
there is a non-zero homomorphism from $N$ to $M$ in that case.

(4) It is possible for $M$ to cross $N$ several times in modules, arrows and 3-cycles. For an example, see Section~\ref{ExampleSection}. 

(5) It is possible for $M$ to cross $N$ and for $N$ to cross $M$, see Section~\ref{ExampleSection} for an example. 
\end{Remark}

\begin{Proposition}\label{CrossingModules}
Let $M$ and $N$ be two string modules over $J(Q,W)$ with corresponding arcs $\gamma_M$ and $\gamma_N$ in $(S, \M)$. 
Then  $\gamma_M$ and $\gamma_N$ cross if and only if  $M$ crosses $N$ or $N$ crosses $M$.
\end{Proposition}

\begin{figure}[ht]
\begin{tikzpicture}[trans/.style={thick,->,shorten >=2pt,shorten <=2pt,>=stealth}]
\node at (0,0){$\begin{tikzpicture}[scale=.7]
    
    \node[scale=.4](a) at (0,7) [label={[label distance=-.05cm, scale=.7]90:$C$}, scale=.7] {$\bullet$};
    \node(a') at (-1.5,7.5) { };
    \node(a'') at (1.5,3){ };
    \node(bb) at (1.5,7.5) { };
    \node(b'') at (-1.5,3) { };
    \node[scale=.4](f) at (2,6) [label={[label distance=-.15cm, scale=.7]10:$B$}, scale=.7] {$\bullet$};
    \node[scale=.4](ff) at (-2,6) [label={[label distance=-.15cm, scale=.7]170:$A$}, scale=.7] {$\bullet$};
    \node[scale=.4](g) at (2,4) [label={[label distance=-.15cm, scale=.7]-10:$B'$}, scale=.7] {$\bullet$};
    \node[scale=.4](gg) at (-2,4) [label={[label distance=-.15cm, scale=.7]-170:$A'$}, scale=.7] {$\bullet$};
    \node[scale=.4](h) at (0,3) [label={[label distance=-.05cm, scale=.7]270:$C'$}, scale=.7] {$\bullet$};
    
    \draw (a.center) -- (f.center);
    \draw (a.center) -- (ff.center);
    \path (ff.center) edge (f.center);
    \path (gg.center) edge (g.center);
    \draw [dotted, line width=.8] (ff.center)--(gg.center) (f.center)--(g.center);
    \draw (gg.center)--(h.center) (g.center)--(h.center);
    
    \path[color=red, line width=.7] (a'.center) edge node [pos=.7, rotate=115, scale=.7]{$<$} ( $(a')!0.4!(a'') $ );
    \draw[color=red, line width=.7] (a'.center) -- ( $(a')!0.4!(a'') $ );
    \path[color=red, line width=.7] ( $(a')!0.4!(a'') $ )[dotted] edge node [pos=.6, fill=white,outer sep=1mm,scale=.6, bend right]{$\gamma_M$} ( $(a')!0.7!(a'') $ ) ;
    \draw[color=red, line width=.7] ( $(a')!0.7!(a'') $ ) -- (a''.center);
    
    \path[color=tealblue, line width=.7] (bb.center) edge node [pos=.65, rotate=65, scale=.7]{$<$} ( $(bb)!0.4!(b'') $ );
    \draw[color=tealblue, line width=.7] (bb.center) -- ( $(bb)!0.4!(b'') $ );
    \path[color=tealblue, line width=.7] ( $(bb)!0.4!(b'') $ )[dotted] edge node [pos=.65, fill=white,outer sep=1mm,scale=.6, bend right]{$\gamma_N$} ( $(bb)!0.7!(b'') $ ) ;
    \draw[color=tealblue, line width=.7] ( $(bb)!0.7!(b'') $ ) -- (b''.center);
        
    \draw[trans, color=brown] (-1.5,6.4)--(-1.3,5.9);
    \node[color=brown, scale=.65] at (-1.1,6.15){$\alpha$};
    \draw[trans, color=brown] (-1.4,4.1)--(-1.5,3.6);
    \node[color=brown, scale=.65] at (-1.2,3.8){$\delta$};
    \draw[trans, color=brown] (1.4,5.9)--(1.6,6.4);
    \node[color=brown, scale=.65] at (1.2,6.2){$\varepsilon$};
    \draw[trans, color=brown] (1.5,3.65)--(1.4,4.1);
    \node[color=brown, scale=.65] at (1.2,3.8){$\beta$};
    \end{tikzpicture}$};
\end{tikzpicture}
\caption{ Overlap crossing where both the arcs $\gamma_M$ and $\gamma_N$ might self-cross multiple times.} \label{FigMcrossesN}
\end{figure}
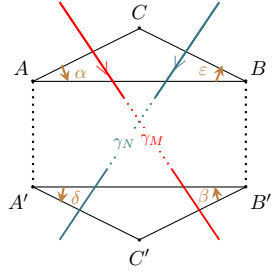

{\it Proof:} First assume, without loss of generality, that $M$ crosses $N${, that is, Definition}~\ref{DefinitionMcrossesN} (1), (2) or (3) holds. First suppose that $M$ crosses $N$ in a module. Let $w_M$ and $ w_N$ be the associated strings as defined in~\ref{DefinitionMcrossesN} (1). Then if the arrows $\alpha, \beta, \delta, \varepsilon$ all exist, we have  a local configuration as in Figure~\ref{FigMcrossesN}.

The endpoints $A,B,C,A',B'$ and $C'$ in Figure~\ref{FigMcrossesN} might (all) coincide. If $\alpha$ does not exists then in the orientation given in Figure~\ref{FigMcrossesN}, the arc $\gamma_M$ starts at $C$ and by definition the arrow $\varepsilon$ must exist since, by definition, the strings {$w_M$ and $w_N$} do not start at the same vertex. Similarly for the end points hence either $\delta$ or $\varepsilon$ exists and thus {$\gamma_M$ and $\gamma_N$} cross. 

\begin{figure}[ht]
\begin{tikzpicture}
\node at (0,0){$\begin{tikzpicture}[trans/.style={thick,->,shorten >=2pt,shorten <=2pt,>=stealth}]
\node at (3.7,0)[rotate=90]{$\begin{tikzpicture}[scale=.7]
\node[scale=.4] (a) at (1,1.7) [label={[label distance=0cm, scale=.7,rotate=-90]170:$A$}, scale=.7] {$\bullet$};
\node[scale=.4] (b) at (0,-.8) [label={[label distance=-.1cm, scale=.7,rotate=-90]-10:$C$}, scale=.7] {$\bullet$}; 
\node[scale=.4] (c) at (2,-.8) [label={[label distance=-.1cm, scale=.7,rotate=-90]10:$B$}, scale=.7] {$\bullet$};
\draw (a.center)--(b.center)--(c.center)--(a.center);

\path[color=red] (b.center) edge node[pos=.5, scale=.7, rotate=35]{$<$} (2.5,1.5);
\path[color=red] (b.center) edge node[pos=.74, fill=white,outer sep=1mm, scale=.6, rotate=-90]{$\gamma_M$} (2.5,1.5);
\draw[color=red, dotted] (2.5,1.5) -- ( 3,2 );
\path[color=tealblue] (c.center) edge node[pos=.5, scale=.7, rotate=-40]{$<$} (-.5,1.5);
\path[color=tealblue] (c.center) edge node[pos=.74, fill=white,outer sep=1mm, scale=.6, rotate=-90]{$\gamma_N$} (-.5,1.5);
\draw[color=tealblue, dotted] (-.5,1.5) -- ( -1,2 );

\draw[trans, color=brown] (1.3,1.2) -- (.7, 1.2);
\node[color=brown, scale=.65, bend left,rotate=-90] at (1,1){$\alpha$};

\end{tikzpicture}$};
\end{tikzpicture}
$};
\node at (0,-1.5){(a)};

\node at (5,0){$\begin{tikzpicture}[trans/.style={thick,->,shorten >=2pt,shorten <=2pt,>=stealth}]
\node at (7,-.5){\begin{tikzpicture}[scale=.7]
\node[scale=.4] (a) at (1,2)  [label={[label distance=-.1cm, scale=.7]10:$C$}, scale=.7]{$\bullet$};
\node[scale=.4] (b) at (0,0)  [label={[label distance=-.1cm, scale=.7]190:$A$}, scale=.7]{$\bullet$}; 
\node[scale=.4] (c) at (2,0)  [label={[label distance=-.1cm, scale=.7]-30:$B$}, scale=.7]{$\bullet$};
\node (d) at (.4,1.7)  { };
\node (e) at (2.5, 1)  { };
\node (f) at (1, -.5)  { };
\node (g) at (1.2,-.6) { };
\node (h) at (-.4,1.7) {};
\draw (a.center)--(b.center)--(c.center)--(a.center);

\path[line width=.7, color=red] (d.center) edge node[pos=.4, scale=.7, rotate=-65]{$>$} (g.center);
\draw[color=red, dotted] (.15,2.5) -- (d.center);
\draw[color=red, dotted] (g.center) -- (1.4,-1.3);
\path[line width=.7, color=red] (d.center) edge node[pos=0, fill=white,outer sep=1mm, scale=.6]{$\gamma_M$} (g.center);

\path[line width=.7, color=tealblue] (b.center) edge node[pos=.55, scale=.7, rotate=15]{$>$} (e.center);
\path[line width=.7, color=tealblue] (b.center) edge node[pos=.88, fill=white,outer sep=1mm, scale=.6]{$\gamma_N$} (e.center);
\draw[color=tealblue, dotted] (e.center) -- (3.4,1.3);

\draw[trans, color=brown] (.24,.57) -- (.45, -.1);
\node[color=brown, scale=.65] at (.38,-.2){$\alpha$};

\draw[trans, color=brown] (1.55, -.1) -- (1.85,.5);
\node[color=brown, scale=.65] at (2.1,.3){$\beta$};

\draw[trans, color=brown] (1.3,1.6) -- (.7,1.6);
\node[color=brown, scale=.65] at (1.4,1.6){$\gamma$};
\end{tikzpicture}};
\end{tikzpicture}$};

\node at (5,-1.5){(b)};
\end{tikzpicture}

\caption{Local configurations of crossings of $\gamma_M$ and $\gamma_N$ induced by the module $M$ crossing the module $N$, corresponding in (a) to a crossing as in  Def~\ref{DefinitionMcrossesN} (2) and in (b) to  Def~\ref{DefinitionMcrossesN} (3).}\label{FigArrowCrossing}
\end{figure}
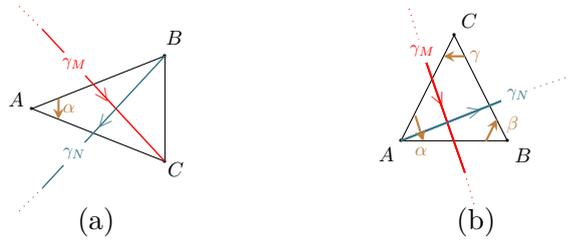

If $M$ and $N$ cross as in Def~\ref{DefinitionMcrossesN} (2), then with the induced orientation on $\gamma_M$ and $ \gamma_N,$ we have  a local configuration as in Figure~\ref{FigArrowCrossing} (a) and thus $\gamma_M$ and $ \gamma_N$ cross.

If $M$ and $N$ cross as in Def~\ref{DefinitionMcrossesN} (3), then with the induced orientation on $\gamma_M, \gamma_N,$ we have  a local configuration as in Figure~\ref{FigArrowCrossing}  (b) and $\gamma_M$ and $\gamma_N$ cross.

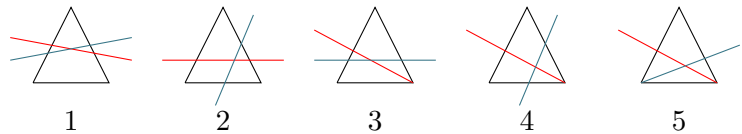
\begin{figure}[ht]
\begin{tikzpicture}
\node (a1) at (0,0){};
\node (b1) at (1,0){};
\node (c1) at (.5,1){};

\node (a2) at (2,0){};
\node (b2) at (3,0){};
\node (c2) at (2.5,1){};

\node (a3) at (4,0){};
\node (b3) at (5,0){};
\node (c3) at (4.5,1){};

\node (a4) at (6,0){};
\node (b4) at (7,0){};
\node (c4) at (6.5,1){};

\node (a5) at (8,0){};
\node (b5) at (9,0){};
\node (c5) at (8.5,1){};

\node at (.5,-.5){1};
\node at (2.5,-.5){2};
\node at (4.5,-.5){3};
\node at (6.5,-.5){4};
\node at (8.5,-.5){5};

\draw (a1.center)--(b1.center)--(c1.center)--(a1.center);
\draw[color=red] (-.3,.6)--(1.3,.3);
\draw[color=tealblue] (-.3,.3)--(1.3,.6);

\draw (a2.center)--(b2.center)--(c2.center)--(a2.center);
\draw[color=red] (1.7,.3)--(3.3,.3);
\draw[color=tealblue] (2.4,-.3)--(2.9,.9);

\draw (a3.center)--(b3.center)--(c3.center)--(a3.center);
\draw[color=red] (3.7,.7)--(b3.center);
\draw[color=tealblue] (3.7,.3)--(5.3,.3);

\draw (a4.center)--(b4.center)--(c4.center)--(a4.center);
\draw[color=red] (5.7,.7)--(b4.center);
\draw[color=tealblue] (6.4,-.3)--(6.9,.9);

\draw (a5.center)--(b5.center)--(c5.center)--(a5.center);
\draw[color=red] (7.7,.7)--(b5.center);
\draw[color=tealblue] (a5.center)--(9.3,.5);

\end{tikzpicture}
\caption{ The five possible ways two arcs can cross in any given triangle.}\label{5crossingconfigurations}
\end{figure}

Conversely, suppose that {$\gamma_M$ and $ \gamma_N$} cross. Then there are 5 possible local configurations of a crossing of $\gamma_M$ and $\gamma_N$ as in  Figure~\ref{5crossingconfigurations}.

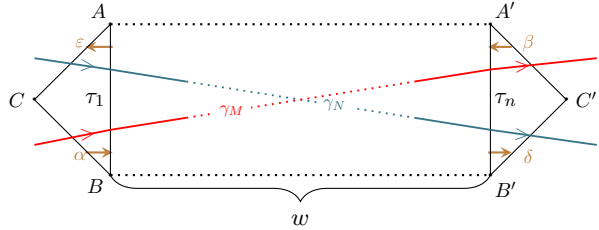
\begin{figure}[ht]
\begin{tikzpicture}[trans/.style={thick,->,shorten >=2pt,shorten <=2pt,>=stealth}]

\node[scale=.4](a) at (0,2) [label={[label distance=-.15cm, scale=.7]170:$A$}, scale=.7] {$\bullet$};
\node[scale=.4](b) at (0,0) [label={[label distance=-.15cm, scale=.7]-100:$B$}, scale=.7] {$\bullet$};

\node[scale=.4](a') at (5,2) [label={[label distance=-.15cm, scale=.7]10:$A'$}, scale=.7] {$\bullet$};
\node[scale=.4](b') at (5,0) [label={[label distance=-.15cm, scale=.7]-10:$B'$}, scale=.7] {$\bullet$};

\node[scale=.4](c) at (-1,1) [label={[label distance=-.05cm, scale=.7]180:$C$}, scale=.7] {$\bullet$};
\node[scale=.4](c') at (6,1) [label={[label distance=-.05cm, scale=.7]0:$C'$}, scale=.7] {$\bullet$};

\node (ab) at ( $(a)!0.3!(b) $){};
\node (a'b') at ( $(a')!0.3!(b') $){};
\node (ba) at ( $(a)!0.7!(b) $){};
\node (b'a') at ($(a')!0.7!(b') $){};

\node at (-.2,1)[scale=.8]{$\tau_1$};
\node at (5.2,1)[scale=.8]{$\tau_{n}$};

\draw[trans, color=brown] (.1,1.7)--(-.4,1.7);
\node[color=brown, scale=.65] at (-.4,1.75){$\varepsilon$};
\draw[trans, color=brown] (-.4,.3)--(.1,.3);
\node[color=brown, scale=.65] at (-.4,.25){$\alpha$};
\draw[trans, color=brown] (5.35,1.7)--(4.9,1.7);
\node[color=brown, scale=.65] at (5.5,1.75){$\beta$};
\draw[trans, color=brown] (4.9,.3)--(5.35,.3);
\node[color=brown, scale=.65] at (5.5,.25){$\delta$};

\draw (a.center) -- (b.center) (a'.center) -- (b'.center) (a.center) -- (c.center) (c.center) -- (b.center) (a'.center) -- (c'.center) (c'.center) -- (b'.center);
\draw [dotted, line width=.8] (a.center)--(a'.center) (b.center)--(b'.center);


\path[color=tealblue, line width=.7] (-1,1.55) edge node [pos=.7, rotate=-5, scale=.7]{$>$} (ab.center);
\draw[color=tealblue, line width=.7] (ab.center) -- ( $(ab)!0.2!(b'a')$);
\path[color=tealblue, line width=.7] ( $(ab)!0.2!(b'a') $ )[dotted] edge node [pos=.65, fill=white,outer sep=1mm,scale=.6, bend right]{$\gamma_N$} 
($(ab)!0.8!(b'a')$);
\draw [color=tealblue, line width=.7] ( $(ab)!0.8!(b'a') $ ) --(b'a'.center);
\path[color=tealblue, line width=.7] (b'a'.center) edge node [pos=.3, rotate=-5, scale=.7]{$>$} (6.4,.4);

\path[color=red, line width=.7] (-1,.4) edge node [pos=.7, rotate=5, scale=.7]{$>$} (ba.center);
\draw[color=red, line width=.7] (ba.center) -- ( $(ba)!0.2!(a'b')$);
\path[color=red, line width=.7] ( $(ba)!0.2!(a'b') $ )[dotted] edge node [pos=.2, fill=white,outer sep=1mm,scale=.6, bend right]{$\gamma_M$} 
($(ba)!0.8!(a'b')$);
\draw [color=red, line width=.7] ( $(ba)!0.8!(a'b') $ ) --(a'b'.center);
\path[color=red, line width=.7] (a'b'.center) edge node [pos=.3, rotate=5, scale=.7]{$>$} (6.4,1.55);

\draw [decorate,decoration={brace,amplitude=10pt},xshift=-4pt,yshift=-9pt]
(b'.center) -- (b.center); 
\node at (2.5,-.6) {\footnotesize $w$};
\end{tikzpicture}
\caption{ Crossing of $\gamma_M$ and $\gamma_N$ where $\gamma_M$ and $\gamma_N$ both successively cross $\tau_1, \ldots, \tau_n$.} \label{FigOverlapCrossinginProofProp}
\end{figure}

By homotopy, configuration (1) can always be reduced to either configuration (2) or (3). So suppose that $\gamma_M$ and $\gamma_N$ are as in configuration (2) or (3) of Figure~\ref{5crossingconfigurations}. Then  locally the crossing takes place in a configuration as in Figure~\ref{FigOverlapCrossinginProofProp} where $\gamma_M$ might start or end at $C$ or $C'$ and $\gamma_N$ might  start or end at $C$ or $C'$ but not both $\gamma_M$ and $\gamma_N$ start or end simultaneously at either point. Let $w_M$ and $w_N$ be the strings associated to the orientation of $\gamma_M$ or $\gamma_N$ in Figure~\ref{FigOverlapCrossinginProofProp}.  
 We could also have chosen the inverse orientation for both $\gamma_M$ and $\gamma_N$ and then have worked with $w^{-1}_M$ and $w^{-1}_N$.  
Given the chosen orientation, let $\gamma_{AB}$ be the first arc of $T$ both arcs $\gamma_M$ and $\gamma_N$ cross with respect to the crossing under consideration and suppose that $\gamma_M$ and $\gamma_N$ then both successively cross  the arcs $\gamma_{AB}=\tau_1,\tau_2, \cdots, \tau_n$   of the triangulation $T$ (where possibly $\tau_i=\tau_j$ for some $i,j$).

In order for $\gamma_M$ and $\gamma_N$ to cross {as in Figure~\ref{FigOverlapCrossinginProofProp} we must have either $\gamma_M$ crossing $\tau_{CB}$ and $\gamma_N$ crossing $\tau_{AC}$ or $\gamma_N$ starting at $C$ or if $\gamma_M$ starts at $C$ then $\gamma_N$ must cross $\tau_{AC}.$

If $\gamma_M$ crosses $\tau_{CB}$ then $w_M=u_M \alpha w  \beta^{-1} v_M$ and $w_N=u_N \varepsilon^{-1} w  \delta v_N$ or $w_N=w \delta v_N$ for some strings $w, u_M, u_N, v_M $ and $v_ N$. If we set $\Pred (\alpha) = u_M$,  $ \Pred (\varepsilon) = u_N$, $\Succ (\beta) = v_M$ and $\Succ (\delta) = v_n$ then   $w_M$ and $w_N$ are exactly as in  Definition~\ref{DefinitionMcrossesN} (1) and thus $M$ crosses $N$ in a module.

If $\gamma_M$ starts at $C$ then $w_M= w \beta^{-1} v_M$ and $w_N=u_N \varepsilon^{-1}w \delta v_N$ for some strings $w, u_M, u_N, v_M $ and $v_ N$. Again the result follows if we set $\Pred (\alpha) = u_M$,  $ \Pred (\varepsilon) = u_N$, $\Succ (\beta) = v_M$ and $\Succ (\delta) = v_n$. 

Similarly, if $\gamma_M$ and $ \gamma_N$ end at $D$ then the strings $w_M$ and $w_N$ are again exactly as in  Definition~\ref{DefinitionMcrossesN} (1) and $M$ crosses $N$ in a module.

Suppose $\gamma_M$ and $\gamma_N$ cross locally as in configuration $(4)$ in Figure~\ref{5crossingconfigurations}.  Then there is an orientation of $\gamma_M$ and $\gamma_N$ such that  we have a local configuration as in  Figure~\ref{FigArrowCrossing} (a).

Hence the crossing cannot be moved outside  the triangle $ABC$ by homotopy
 without increasing the number of intersections of $\gamma_M$ or $\gamma_N$ with
  the arcs of the triangulation.  Therefore $e(w_M)=v_{\tau_{AB}}$ and $s(w_N)=v_{\tau_{AC}}$ and thus we can form the string $w_M \alpha w_N$ 
  and $M$
   crosses $N$ in an arrow.

Suppose $\gamma_M$ and $\gamma_N$ cross locally as in configuration $(5)$ in Figure~\ref{5crossingconfigurations} which corresponds to a crossing of oriented arcs as in Figure~\ref{FigArrowCrossing} (b). Again the crossing cannot be moved  outside the triangle $ABC$ by homotopy without increasing the number of intersection of $\gamma_M$ and $\gamma_N$ with the arcs of  the triangulation.

Then $ \xymatrix{
         		 a \ar[r]|\alpha	      & b \ar[d]|\gamma \\
                           &  \ar[ul]|\sigma c }$ is a $3$-cycle such that  $\alpha\gamma, \gamma\sigma, \sigma\alpha$ are in $I$. Furthermore, $\alpha$ is in $w_M,$ $s(w_N)=c$ and $w_N$ does not start with $\alpha$ or $\gamma$ nor their inverses. Hence $M$ crosses $N$ in a $3$-cycle.

\hfill $\Box$

 We now show the existence of  some arrows that occur if $M$ crosses $N$ in a module. These arrows will be needed in Definition~\ref{DefinitionOfModulesInResolution} (1) below. 

\begin{Lemma}\label{ExistenceOfArrows} Let $Q$ be a quiver associated to a triangulation of  $(S,\M)$.
Suppose two strings $w$ and $v$ are of the form $w = \Pred(\alpha) \stackrel{\alpha}{\longrightarrow} u  \stackrel{\beta}{\longleftarrow} \Succ(\beta)$ and $v= \Pred(\gamma)
 \stackrel{\gamma}{\longleftarrow} u \stackrel{\delta}{\longrightarrow} \Succ(\delta)$ where $u \in \S$ and $\alpha, \beta, \gamma, \delta$ are arrows in $Q_1$. Then
  if the arrows $\alpha$ and $\gamma$ exist with $\alpha \gamma \in I$ then 
 there exists an arrow $\sigma$ in $Q$ such that    $ \xymatrix{
         		 a \ar[r]|\alpha	      & b \ar[d]|\gamma \\
                           &  \ar[ul]|\sigma c }$ is a 3-cycle in $Q$ and  $\gamma \sigma \in I$
 and $\sigma \alpha \in I$ and $b = s(u)$. 


Similarly,  if the arrows $\beta$ and $\delta $ exist with $\beta \delta \in I$ then
 there exists an arrow $\rho$ in $Q$ such that    $ \xymatrix{
         		 d    \ar[dr]|\delta  & e \ar[l]|\beta	 \\
                           &    \ar[u]|\rho  f }$ is a 3-cycle in $Q$  and $\beta \rho \in I$
 and $\rho \delta \in I$ and $d = t(u)$.

\end{Lemma}

{\it Proof:} Let $\xymatrix{u = u_1 \ar@{-}[r]^{ \mu_1} & u_2 \ar@{-}[r]^{ \mu_2}  & u_3 \; \cdots \; u_r \ar@{-}[r]^{ \mu_r} & u_{r+1}}$. 
Then $\alpha \mu_1 \in \S$ and
$\gamma^{-1} \mu_1 \in \S$ and $\alpha \gamma $ is a non-zero path in  $Q$. Since either $\mu_1 \in Q_1$ or $\mu_1^{-1} \in Q_1$ we have either $\alpha \mu_1 \in kQ$ and $\alpha \mu_1 \notin I$
 or $ \mu_1^{-1} \gamma  \in kQ$ and $\mu_1^{-1} \gamma \notin I$. Since $J(Q,W)$ is gentle, by (S2) we have $\alpha \gamma \in I$. Since $J(Q,W)$ is a gentle algebra
 coming from a surface triangulation this implies 
 that there exists an arrow $\sigma$ such that $\alpha \gamma \sigma$ is a 3-cycle in $Q$ and that $\gamma \sigma \in I$
 and $\sigma \alpha \in I$, see Figure \ref{ArrowSigma}. 
 
\begin{figure}[ht]
\begin{tikzpicture}[trans/.style={thick,->,shorten >=2pt,shorten <=2pt,>=stealth}]
\draw (0,0)--(2,0)--(1,1.2)--(0,0);
\draw[color=red, thick] (.3,1.7)  .. controls (.5, .5).. node [pos=.55, fill=white,outer sep=1mm,scale=.7]{$\gamma_1$}  (2, -.7);
\draw[color=tealblue, thick] (1.7,1.7) .. controls (1.5, .5).. node [pos=.55, fill=white,outer sep=1mm,scale=.7]{$\gamma_2$} (0, -.7);
\draw[trans, color=brown, densely dotted, thick] (1.25,1)--(.75,1);
\node[color=brown] at (1.4,1){$\sigma$};
\draw[trans, color=brown, thick] (.3,.45)--(.3,-.05);
\node[color=brown] at (0,.2){$\alpha$};
\draw[trans, color=brown, thick] (1.7,-.05)--(1.7,.45);
\node[color=brown] at (2,.2){$\gamma$};
\end{tikzpicture}
\caption{Existence of the arrow $\sigma$ in the case of a crossing in a module.}
\label{ArrowSigma}
\end{figure}
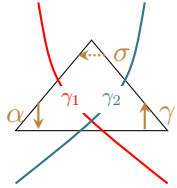

A similar argument proves the existence of the 3-cycle containing $\rho$. \hfill $\Box$

\subsection{Smoothing of crossings for string modules}

Given two arcs in $(S,\M)$ that cross, the smoothing of a crossing gives rise to four new arcs as in Figure~\ref{Figuresmoothingcrossing}. 
We interpret these arcs in terms of string modules.

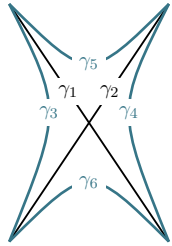
\begin{figure}[ht]
\begin{tikzpicture}[scale=.7]
\node(a') at (-1.5,7.5) { };
\node(a'') at (1.5,3){ };
\node(b') at (1.5,7.5) { };
\node(b'') at (-1.5,3) { };

\path[line width=.7] (a'.center) edge node [pos=.37, fill=white,outer sep=1mm,scale=.7, bend right]{$\gamma_1$} (a''.center);
\path[line width=.7] (b'.center) edge node [pos=.37, fill=white,outer sep=1mm,scale=.7, bend right]{$\gamma_2$} (b''.center);

\draw[color=tealblue, line width=1] (a'.center) .. controls +(1.5,-1.5) ..   node [pos=.5, fill=white,outer sep=1mm,scale=.7]{$\gamma_5$}(b'.center);
\draw[color=tealblue, line width=1] (a''.center) .. controls +(-1.5,1.5) ..  node [pos=.5, fill=white,outer sep=1mm,scale=.7]{$\gamma_6$}(b''.center);
\draw[color=tealblue, line width=1] (a'.center) .. controls +(1,-2) ..  node [pos=.5, fill=white,outer sep=1mm,scale=.7]{$\gamma_3$}(b''.center);
\draw[color=tealblue, line width=1] (b'.center) .. controls +(-1,-2) ..  node [pos=.5, fill=white,outer sep=1mm,scale=.7]{$\gamma_4$}(a''.center);
\end{tikzpicture}

\caption{Smoothing a crossing of $\gamma_1$ and $\gamma_2$.}\label{Figuresmoothingcrossing}
\end{figure}

Given two string modules $M_1 = M(w_1)$ and $M_2=M(w_2)$ in $\mod  J(Q,W)$ such that the associated arcs $\gamma_1$ and $\gamma_2$ cross
in $(S,\M)$, we define four new modules $M_3= M(w_3)$, $M_4= M(w_4)$, $M_5= M(w_5)$
and $M_6= M(w_6)$.  By possibly relabeling  $M_1$ and $M_2$ we can always  assume that $M_1$ crosses $M_2$.

\begin{Definition}\label{DefinitionOfModulesInResolution}
\rm{  Let $M_1$ and $M_2$ be two  string modules over $J(Q,W)$ with strings $w_1$ and $w_2$, respectively.

(1) Suppose $M_1$ crosses $M_2$ in a module $M(w)$. That is
$$w_1 = \Pred(\alpha) \stackrel{\alpha}{\longrightarrow} w  \stackrel{\beta}{\longleftarrow} \Succ(\beta)
\;\;\;\; \mbox{               and              } \;\;\;\;
w_2= \Pred(\gamma)  \stackrel{\gamma}{\longleftarrow} w \stackrel{\delta}{\longrightarrow} \Succ(\delta)$$
where $\alpha, \beta, \gamma, \delta$ are arrows in $Q_1$.  Define
\begin{itemize}
\item[] $w_3 = \Pred(\alpha) \stackrel{\alpha}{\longrightarrow}  w \stackrel{\delta}{\longrightarrow} \Succ(\delta)$
\item[] $w_4 =  \Pred(\gamma)  \stackrel{\gamma}{\longleftarrow} w  \stackrel{\beta}{\longleftarrow} \Succ(\beta)$
\item[] $w_5 = \left\{ \begin{array}{ll}
				\Pred(\alpha) 	\stackrel{\sigma}{\longleftarrow} \Pred^{-1}(\gamma) & \mbox{if $s(w_1) \neq s(w)$ and $s(w_2) \neq s(w)$} \\
				\Pred(\alpha)_c & \mbox{if $s(w_2) = s(w)$} \\
				\Pred(\gamma)_h  & \mbox{if $s(w_1) = s(w)$}
				\end{array}\right.$
\item[] $w_6 = \left\{ \begin{array}{ll}
				\Succ^{-1}(\beta) 	\stackrel{\rho}{\longleftarrow} \Succ(\delta) & \mbox{if $t(w_1) \neq t(w)$ and $t(w_2) \neq t(w)$} \\
			{_c}\Succ(\beta) & \mbox{if $t(w_2) = t(w)$} \\
				 {_h}\Succ(\delta)  & \mbox{if $t(w_1) = t(w)$}
				\end{array}\right.$				
\end{itemize}
where $\sigma, \rho \in Q_1$ are as in Lemma~\ref{ExistenceOfArrows}. 

(2) Suppose $M_1$ {crosses} $M_2$ in an arrow $\alpha$. Suppose without loss of generality that $s(\alpha) = t(w_1)$ and $t(\alpha) = s(w_2)$.  Define
\begin{itemize}
\item[] $w_3 = w_1 \stackrel{\alpha}{\longrightarrow}  w_2 $
\item[] $w_4 =0$
\item[] $w_5 = (w_1)_c$ 
\item[] $w_6 = {_h}w_2.$
\end{itemize}

(3) Suppose $M_1$ crosses $M_2$ in the 3-cycle  $ \xymatrix{
         		 s(w_2) \ar[r]|\gamma	       & s(\alpha) \ar[d]|\alpha \\
                           & \ar[ul]|\beta t(\alpha) }$
where (modulo inverting $w_1$) $w_1 = \Pred(\alpha) \stackrel{\alpha}{\longrightarrow}  \Succ(\alpha)$. Define
\begin{itemize}
\item[] $w_3 = \Pred(\alpha) \stackrel{\gamma}{\longleftarrow}  w_2 $
\item[] $w_4 = {_c}(s(\alpha) \stackrel{\alpha}{\longrightarrow}  \Succ(\alpha))$
\item[] $w_5 = (\Pred(\alpha) \stackrel{\alpha}{\longrightarrow} t(\alpha))_h$
\item[] $w_6 = w_2^{-1} \stackrel{\beta}{\longleftarrow} \Succ(\alpha)   .$
\end{itemize}
}\end{Definition}

 Just as is the case for  snake graphs, if $M_1 = M_2$  then one self-crossing gives rise to two crossings of modules when we consider two copies of the same module.
  However, the smoothing of either of the two crossings results in the same four modules. Therefore it is enough to consider only one of the crossings in each case.

\begin{Proposition}  Let $\gamma_1 $ and $\gamma_2$ be two crossing arcs in $(S,M,T)$ with strings $w_1$ and $ w_2$ respectively. Consider a given crossing of $\gamma_1 $ and $\gamma_2$. Then the arcs $\gamma_3, \gamma_4, \gamma_5$ and $ \gamma_6$ defined by the strings in Def~\ref{DefinitionOfModulesInResolution} correspond to the arcs obtained by smoothing the given crossing of $\gamma_1$ and $\gamma_2$ as in  Definition~\ref{SmoothingCrossingArcs}.
\end{Proposition}

{\it Proof:} This follows directly from the resolution and grafting of snake graphs in \cite{CS}  which for the convenience of the reader we have recalled in Section 2.4 in Theorems~\ref{Overlap=Crossing} and~\ref{Grafting=Crossing}. \hfill $\Box$

 We show that this enables us for any crossing to characterise whether it gives rise to a short exact sequence in $J(Q,W)$ or not and if there is a short exact sequence we describe the middle terms in terms of the arcs obtained from smoothing the crossing.

\begin{Theorem}\label{ShortExactSequences}
 Let $M_1$ and $M_2$ be  string modules in $\mod  J(Q,W)$. 

(1) If $M_1$ crosses $M_2$ in a module then the modules $M_3$ and $M_4$ defined in {Definition~\ref{DefinitionOfModulesInResolution}}(1) above give a non-split short exact sequence in $\mod  J(Q,W)$ of the form
$$ 0 \longrightarrow M_2 \longrightarrow M_3 \oplus M_4 \longrightarrow M_1 \longrightarrow 0.$$

(2) If $M_1$ crosses $M_2$ in an arrow then the module  $M_3$  defined in {Definition~\ref{DefinitionOfModulesInResolution}}(2) above gives a non-split short exact sequence in $\mod  J(Q,W)$ of the form
$$ 0 \longrightarrow M_2 \longrightarrow M_3  \longrightarrow M_1 \longrightarrow 0.$$

(3) If $M_1$ crosses $M_2$ in a 3-cycle then the modules $M_3$ and $M_4$ defined in {Definition~\ref{DefinitionOfModulesInResolution}}(3) above do not give rise to an element in $\Ext_{J(Q,W)}^1(M_1, M_2)$.
\end{Theorem}

{\it Proof:} We use the notation of Definition~\ref{DefinitionOfModulesInResolution}. 

(1) Set $w_1 = P_1 w S_1$ where  $P_1 = \Pred(\alpha) \stackrel{\alpha}{\longrightarrow}$ and $S_1 = \stackrel{\beta}{\longleftarrow} \Succ(\beta)$ and set
$w_2 = P_2 w S_2$ where $P_2 = \Pred(\gamma)  \stackrel{\gamma}{\longleftarrow} $ and  $S_2= \stackrel{\delta}{\longrightarrow} \Succ(\delta)$. By definition
we do not  simultaneously have $P_1 = 0$ and $P_2 = 0 $ or $S_1 = 0$ and $ S_2 = 0$. Thus by \cite{S} the sequence 
$\mbox{$ 0 \longrightarrow M_2 \longrightarrow M_3 \oplus M_4 \longrightarrow M_1 \longrightarrow 0$}$ is a non-split short exact sequence.

(2) This follows directly from the canonical embedding $M(w_2) \longrightarrow M(w_1 \stackrel{\alpha}{\longrightarrow}  w_2)$ and the canonical projection
$M(w_1 \stackrel{\alpha}{\longrightarrow}  w_2)\longrightarrow  M(w_1),$ see \cite{BR}.

(3) Since $\dim M({_c}(\stackrel{\alpha}{\longrightarrow}  \Succ(\alpha)) < \dim M(\stackrel{\alpha}{\longrightarrow}  \Succ(\alpha))$
it is immediate by comparing the dimensions of $M_1 \oplus M_2$ and $M_3 \oplus M_4$ that these four modules cannot give rise to a short exact  sequence in $\mod  J(Q,W)$. \hfill $\Box$

\begin{Remark} \label{Geometric interpretation of crossings}{\rm (1) The geometric interpretation of the module crossings in Theorem~\ref{ShortExactSequences} is as follows: 

For Theorem~\ref{ShortExactSequences} (1), the module $M_1$ crossing the module $M_2$  corresponds to a crossing of the corresponding arcs $\gamma_1$ and $\gamma_2$ as in  Figure~\ref{ModuleCrossingPictures} (1) and the modules $M_3$ and $M_4$ correspond to the arcs $\gamma_3$ and $\gamma_4$. 

If $M_1$ crosses $M_2$ as in  Theorem~\ref{ShortExactSequences} (2), in the geometric set-up this corresponds to a crossing of arcs as in Figure~\ref{ModuleCrossingPictures} (2). In this case $\gamma_4$ is either a boundary arc or an arc in the triangulation and the corresponding module $M_4$ is the zero module. 

For modules crossings as in Theorem~\ref{ShortExactSequences} (3), the geometric picture is as in Figure~\ref{ModuleCrossingPictures} (3).

\begin{figure}[ht]
\begin{tikzpicture}[trans/.style={thick,->,shorten >=2pt,shorten <=2pt,>=stealth}]
\node at (0,0){$\begin{tikzpicture}[scale=.7]

\node[scale=.4](a) at (0,7) [label={[label distance=-.05cm, scale=.7]90:}, scale=.7] {$\bullet$};
\node(a') at (-1.5,7.5) { };
\node(a'') at (1.5,3){ };
\node(bb) at (1.5,7.5) { };
\node(b'') at (-1.5,3) { };
\node[scale=.4](f) at (2,6) [label={[label distance=-.15cm, scale=.7]10:}, scale=.7] {$\bullet$};
\node[scale=.4](ff) at (-2,6) [label={[label distance=-.15cm, scale=.7]170:}, scale=.7] {$\bullet$};
\node[scale=.4](g) at (2,4) [label={[label distance=-.15cm, scale=.7]-10:}, scale=.7] {$\bullet$};
\node[scale=.4](gg) at (-2,4) [label={[label distance=-.15cm, scale=.7]-170:}, scale=.7] {$\bullet$};
\node[scale=.4](h) at (0,3) [label={[label distance=-.05cm, scale=.7]270:}, scale=.7] {$\bullet$};


\draw (a.center) -- (f.center);
\draw (a.center) -- (ff.center);
\path (ff.center) edge (f.center);
\path (gg.center) edge (g.center);
\draw [dotted, line width=.8] (ff.center)--(gg.center) (f.center)--(g.center);
\draw (gg.center)--(h.center) (g.center)--(h.center);

\path[color=red, line width=.7] (a'.center) edge node [pos=.7, rotate=115, scale=.7]{$<$} ( $(a')!0.4!(a'') $ );
\draw[color=red, line width=.7] (a'.center) -- ( $(a')!0.4!(a'') $ );
\path[color=red, line width=.7] ( $(a')!0.4!(a'') $ )[dotted] edge node [pos=.6, fill=white,outer sep=1mm,scale=.6, bend right]{$\gamma_1$} ( $(a')!0.7!(a'') $ ) ;
\draw[color=red, line width=.7] ( $(a')!0.7!(a'') $ ) -- (a''.center);

\draw[color=green, line width=1] (a'.center) .. controls +(1,-2) .. node[pos=.4, scale=.7, rotate=110]{$<$} node [pos=.5, fill=white,outer sep=1mm,scale=.6]{$\gamma_3$}(b''.center);
\draw[color=green, line width=1] (b'.center) .. controls +(-1,-2) .. node[pos=.4, scale=.7, rotate=90]{$<$} node [pos=.5, fill=white,outer sep=1mm,scale=.6]{$\gamma_4$}(a''.center);

\path[color=tealblue, line width=.7] (bb.center) edge node [pos=.65, rotate=65, scale=.7]{$<$} ( $(bb)!0.4!(b'') $ );
\draw[color=tealblue, line width=.7] (bb.center) -- ( $(bb)!0.4!(b'') $ );
\path[color=tealblue, line width=.7] ( $(bb)!0.4!(b'') $ )[dotted] edge node [pos=.65, fill=white,outer sep=1mm,scale=.6, bend right]{$\gamma_2$} ( $(bb)!0.7!(b'') $ ) ;
\draw[color=tealblue, line width=.7] ( $(bb)!0.7!(b'') $ ) -- (b''.center);
%
\end{tikzpicture}$};

\node at (4,0){$\begin{tikzpicture}[trans/.style={thick,->,shorten >=2pt,shorten <=2pt,>=stealth}]
\node at (3.7,0)[rotate=90]{$\begin{tikzpicture}[scale=.7]
\node[scale=.4] (a) at (1,1.7) [label={[label distance=0cm, scale=.7,rotate=-90]170:}, scale=.7] {$\bullet$};
\node[scale=.4] (b) at (0,-.8) [label={[label distance=-.1cm, scale=.7,rotate=-90]-10:}, scale=.7] {$\bullet$}; 
\node[scale=.4] (c) at (2,-.8) [label={[label distance=-.1cm, scale=.7,rotate=-90]10:}, scale=.7] {$\bullet$};
\draw (a.center)--(b.center)--(c.center)--(a.center);

\path[color=red] (b.center) edge node[pos=.5, scale=.7, rotate=35]{$<$} (2.5,1.5);
\path[color=red] (b.center) edge node[pos=.74, fill=white,outer sep=1mm, scale=.6, rotate=-90]{$\gamma_1$} (2.5,1.5);
\draw[color=red, dotted] (2.5,1.5) -- ( 3,2 );
\path[color=tealblue] (c.center) edge node[pos=.5, scale=.7, rotate=-40]{$<$} (-.5,1.5);
\path[color=tealblue] (c.center) edge node[pos=.74, fill=white,outer sep=1mm, scale=.6, rotate=-90]{$\gamma_2$} (-.5,1.5);
\draw[color=tealblue, dotted] (-.5,1.5) -- ( -1,2 );

\draw[color=green, line width=1] (2.3,1.6) .. controls +(-1.5,-1) .. node[pos=.2, scale=.7, rotate=35]{$<$} node [pos=.42, fill=white,outer sep=1mm,scale=.6, rotate=-90]{$\gamma_3$}(-.4,1.6);
\draw[color=green, dotted] ( 3,2 ) --(2.3,1.6)  (-1,2)--(-.4,1.6);

\draw[color=green, line width=1] (c.center) .. controls +(0,0) .. node[pos=.65, scale=.7, rotate=5]{$<$}  node [pos=.8, fill=white,outer sep=1mm,scale=.6, rotate=-90]{$\gamma_4$}(b.center);
%
\end{tikzpicture}$};
\end{tikzpicture}$};

\node at (8,0){\begin{tikzpicture}[trans/.style={thick,->,shorten >=2pt,shorten <=2pt,>=stealth}]
\node at (7,-.5){\begin{tikzpicture}[scale=.7]
\node[scale=.4] (a) at (1,2)  [label={[label distance=-.1cm, scale=.7]10:}, scale=.7]{$\bullet$};
\node[scale=.4] (b) at (0,0)  [label={[label distance=-.1cm, scale=.7]190:}, scale=.7]{$\bullet$}; 
\node[scale=.4] (c) at (2,0)  [label={[label distance=-.1cm, scale=.7]-30:}, scale=.7]{$\bullet$};
\node (d) at (.4,1.7)  { };
\node (e) at (2.5, 1)  { };
\node (f) at (1, -.5)  { };
\node (g) at (1.2,-.6) { };
\node (h) at (-.4,1.7) {};
\draw (a.center)--(b.center)--(c.center)--(a.center);

\path[line width=.7, color=red] (d.center) edge node[pos=.4, scale=.7, rotate=-65]{$>$} (g.center);
\draw[color=red, dotted] (.15,2.5) -- (d.center);
\draw[color=red, dotted] (g.center) -- (1.4,-1.3);
\path[line width=.7, color=red] (d.center) edge node[pos=0, fill=white,outer sep=1mm, scale=.6]{$\gamma_1$} (g.center);

\path[line width=.7, color=tealblue] (b.center) edge node[pos=.55, scale=.7, rotate=15]{$>$} (e.center);
\path[line width=.7, color=tealblue] (b.center) edge node[pos=.88, fill=white,outer sep=1mm, scale=.6]{$\gamma_2$} (e.center);
\draw[color=tealblue, dotted] (e.center) -- (3.4,1.3);

\draw[color=green, dotted] (.15,2.5) -- (.6,1.7);
\draw[color=green, line width=1] (.6,1.7) .. controls +(.55,-1) .. node[pos=.18, scale=.7, rotate=115]{$<$} node [pos=.38, 
fill=white,outer sep=1mm,scale=.6, rotate=0]{$\gamma_3$}(2,1);
\draw[color=green, dotted] (2,1) -- (3.4,1.3);

\draw[color=green, line width=1] (.9,-.65) .. controls +(-.3,.3) .. node[pos=.1, scale=.7, rotate=130]{$<$} node [pos=.6, fill=white,outer sep=1mm,scale=.6, rotate=0]{$\gamma_4$}(b.center);
\draw[color=green, dotted] (.9,-.65) -- (1.4,-1.3);

\end{tikzpicture}};
\end{tikzpicture}};

\node at (0,-2){(1)};
\node at (4,-2){(2)};
\node at (8,-2){(3)};
\end{tikzpicture}
\caption{Geometric interpretation of Theorem~\ref{ShortExactSequences}.  } \label{ModuleCrossingPictures} 
\end{figure}
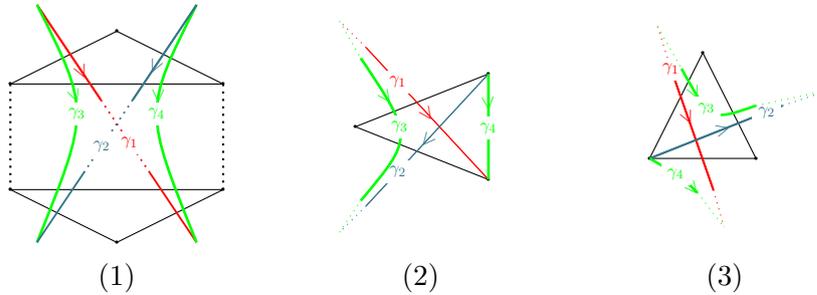

(2) In general, $M_5$ and $M_6$ as defined above never give rise to an element in $\Ext^1_{J(Q,W)}(M_1,M_2)$ nor  in $\Ext^1_{J(Q,W)}(M_2,M_1)$ since $\dim (M_5 \oplus M_6) \leq \dim (M_1 \oplus M_2) - 1$.

(3) If $M_2 = M(\gamma)$ and $M_1 = M(_s\gamma _e)$ for an arc $\gamma$ in $(S,\M)$ where  $_s\gamma _e$ is the arc rotated by the elementary pivot moves as defined in \cite{BZ} then in
Theorem~\ref{ShortExactSequences}(i) and (ii) above we recover the AR-sequences described in \cite{BZ}.}
\end{Remark}


\section{Extensions in the cluster category and dimension formula for Jacobian algebras}

 In this section let $\C(S,\M)$ be the cluster category of a marked surface $(S,\M)$ where all marked points lie in the boundary of $S$ and each boundary component has at least one 
marked point.

In Theorem~\ref{Triangles}  we explicitly describe the middle terms of triangles 
in  $\C(S,\M)$. 
Namely, we show that any crossing of two arcs gives rise to at least one triangle in the cluster category where the middle terms of the triangle are given by one of the pairs of arcs obtained by smoothing the crossing. In almost all cases, except
if the crossing has a self-overlap as defined in Definition~\ref{DefinitionMcrossesN}, the middle terms of the triangle in the opposite direction are given by the other pair of arcs obtained by smoothing the crossing.   

 Let  $M_1 = M(\gamma_1)$ and  $M_2= M(\gamma_2)$ be two  string $J(Q,W)$-modules corresponding to the string
objects $\gamma_1$ and $\gamma_2$ in $\C(S,\M)$. 
 It follows immediately from \cite{Plamondon} (see also \cite[Lemma 4.4]{KY}) that every short exact sequence in $\mod  J(Q,W)$  lifts to a triangle in $\C(S,\M)$  such that the image under
the canonical projection functor from $\C(S,\M)$ to $\mod  J (Q,W)$ is isomorphic to the short
exact sequence. Since both the indecomposable
objects in the cluster category and as well as the indecomposable modules over the
Jacobian algebra correspond to arcs in the surface (where the canonicity of this bijection follows from the Appendix to this paper), we can use the associated string to
explicitly check whether we have a short exact sequence in the module category. It then follows that 
this short exact sequence comes from the triangle corresponding to the same set of arcs
in the associated cluster category.  Therefore by  Theorem~\ref{ShortExactSequences}, we have that

(*) if $M_1$ crosses $M_2$ in a module and if $M_3 = M(\gamma_3)$ and $M_4 = M(\gamma_4)$ are  defined as in Definition~\ref{DefinitionOfModulesInResolution}(1) 
then by  there is a non-split triangle in $\C(S,\M)$ given by
$$\gamma_2 \longrightarrow \gamma_3 \oplus \gamma_4 \longrightarrow \gamma_1 
\longrightarrow \gamma_2[1].$$ Note that in this situation, even if the module crossing 
is self-crossing, we obtain this triangle in $\C(S,\M)$.

(**)  if  $M_1$ crosses $M_2$ in an arrow $\alpha$ and  if $M_3 = M(\gamma_3)$ is defined as 
in Definition~\ref{DefinitionOfModulesInResolution}(2) 
then there is a non-split triangle in $\C(S,\M)$ given by
 $$\gamma_2 \longrightarrow \gamma_3 \oplus \gamma_4 \longrightarrow \gamma_1 \longrightarrow 
 \gamma_2[1]$$ where $\gamma_4 \neq 0$ if and only if $\gamma_4$ is not a boundary arc.

\begin{Theorem}\label{Triangles}
Let $\gamma_1$ and $\gamma_2$ be two  string objects (not necessarily distinct) in $\C(S,\M)$
 such that their corresponding arcs cross in $(S,\M)$.
 Let $\gamma_3, \gamma_4, \gamma_5, \gamma_6$ be the string objects corresponding to 
 the smoothing of a crossing of a suitable orientation of the corresponding arcs 
 $\gamma_1$ and $\gamma_2$. Then there is a non-split triangle in $\C(S,\M)$ given by
 \begin{equation}
\gamma_2 \longrightarrow \gamma_3 \oplus \gamma_4 \longrightarrow \gamma_1 \longrightarrow \gamma_2[1]
\end{equation} 
and  if the crossing of $\gamma_1$ and $\gamma_2$ is not in a self-crossing  overlap in some triangulation of $(S,\M)$  then we obtain a non-split triangle given by 
\begin{equation} \gamma_1 \longrightarrow \gamma_5 \oplus \gamma_6 \longrightarrow \gamma_2 \longrightarrow \gamma_1[1] \end{equation}
where $\gamma_3, \gamma_4, \gamma_5, \gamma_6$ are zero objects in $\C(S,\M)$ if they correspond to boundary arcs.
\end{Theorem}

Before we give a proof of Theorem~\ref{Triangles}	in Section 4.1, we explore some consequences.

\begin{Remark} \label{Facts}{\rm Keeping the notations of Theorem~\ref{Triangles}, suppose that we are given a particular
triangulation of $(S, \M)$. Then we have the following two facts: 

(i) For every crossing of $\gamma_1$ and $\gamma_2$ that corresponds to either an arrow
 crossing  or a 3-cycle crossing
of the associated string modules, we always obtain exactly two triangles in the cluster category 
with middle terms 
given by $\gamma_3 \oplus \gamma_4$ and by $\gamma_5 \oplus \gamma_6$, respectively. 
 
 (ii) For every crossing of  $\gamma_1$ and $\gamma_2$ that corresponds to a non-selfcrossing
  module crossing
of the associated string modules   we obtain exactly two triangles in the cluster category 
with middle terms 
given by $\gamma_3 \oplus \gamma_4$ and by $\gamma_5 \oplus \gamma_6$. 
 If the module crossing is self-crossing, 
  we only 
 obtain one of the two triangles in the cluster category and  the middle terms of  that triangle are 
 given by one of the two pairs of arcs obtained from smoothing the crossing. In Theorem 4.1 we have denoted  this pair of arcs by  $\gamma_3 \oplus \gamma_4$.
 By the 2-Calabi Yau property of the cluster category and by \cite{ZZZ} we know that
  this crossing must also give rise to a second triangle.}
\end{Remark}

\begin{Question}\label{Question}
 It is an open question whether, in the self-crossing module crossing  case, the
  middle terms of the second triangle in Remark~\ref{Facts} (ii) are also induced by the `other' pair of arcs,  which in the notation of Theorem 4.1 are $\gamma_5$ and $\gamma_6$. 
\end{Question}

Combining Theorem~\ref{ShortExactSequences} with Theorem~\ref{Triangles} and \cite[Theorem 3.4]{ZZZ} we obtain a formula for the dimensions of the 
first extension space over the Jacobian algebra in the case where the arcs do not create a self-crossing overlap.

\begin{Corollary}\label{DimensionFormulaModules}
Let $M, N$ be two   string modules over $J(Q,W)$ and let $\gamma_M$ and $\gamma_N$ be the corresponding arcs in $(S,\M)$ such that $\gamma_M$ and $\gamma_N$ have no crossing with self-crossing overlap.  

(1) A basis of $\Ext^1_{J(Q,W)}(M,N)$ is given 
by all short exact sequences arising from $M$ crossing $N$ 
in a module or an arrow and where the middle terms are as described in Theorem~\ref{ShortExactSequences};

(2) We have 
$$\dim \Ext^1_{J(Q,W)}(M,N) + \dim \Ext^1_{J(Q,W)}(N, M) = \Int(\gamma_M, \gamma_N) - k-k'$$
 where $k$ (resp. $k'$) is the number of times that $M$ crosses $N$ (resp.  $N$ crosses $M$)
 in a 3-cycle. In particular, if $M=N$ we have 
 $$2 \dim \Ext^1_{J(Q,W)}(M, M) = \Int(\gamma_M, \gamma_M) - 2k.$$
\end{Corollary}





{\it Proof:}  (1) Fix a crossing $p$ of $\gamma_M$ and $\gamma_N$ and suppose without loss of generality that this crossing yields
 $M$ crosses $N$ in $\mod J(Q,W)$. Then $M$ crosses $N$ either in a  non-self-crossing module, in an arrow, or in a 3-cycle. Let $\gamma_3, \ldots, \gamma_6$ 
 be the arcs obtained by smoothing the crossing $p.$  By view of Theorem~\ref{Triangles}, it is sufficient to check whether 
 $M(\gamma_3), M(\gamma_4)$ or $M(\gamma_5), M(\gamma_6)$ appear as middle terms of a non-split short exact sequence from $N$ to $M.$ 
 If $M$ crosses $N$ in a  non self-crossing  module or in an arrow at $p$, then by Theorem~\ref{ShortExactSequences} there exists a non-split short exact sequence 
 with middle terms given by $M(\gamma_3)$ and $M(\gamma_4)$.  If $M$ crosses $N$ in a 3-cycle, then by Theorem~\ref{ShortExactSequences} 
 this does not result in a short exact sequence from $N$ to $M$
  with middle terms $M(\gamma_3)$ and $M(\gamma_4)$. The modules $M(\gamma_5)$ and $M(\gamma_6)$
  never induce an element in $\Ext^1_{J(Q,W)}(M,N)$.

(2) From the proof of (1) we see that a dimension formula
for $\Ext^1_{J(Q,W)}(M,N)$ accounts only for the number of times $M$ crosses $N$ in a (non self-crossing) module or in an arrow.  The  result then follows from 
the dimension formula in \cite{ZZZ} and
 Theorem~\ref{ShortExactSequences}. \hfill $\Box$

\subsection{Proof of Theorem~\ref{Triangles}}

Our general strategy for the proof of Theorem~\ref{Triangles} is as follows.

We consider each type of crossing separately. That is, given a fixed triangulation $T$ of $(S,\M)$ and two  string objects $\gamma_1$ and $\gamma_2$ in $\C(S,\M)$ corresponding to two crossing arcs in $(S,\M)$, we treat the different crossings of the corresponding string modules $M_1 = M(w_1)$ and $M_2 = M(w_2)$ one by one.

If the crossing under consideration is a crossing in a module then by (*) above  we obtain one triangle with two middle terms given by the string objects $\gamma_3$ and $\gamma_4$. The other triangle is obtained by possibly flipping the overlap to an \emph{orthogonal overlap} (see  proof of Theorem 4.1, Section 4.1.1, {\it Case 1 } below for the definition. However, sometimes this is not possible. In this case we adapt a strategy similar to the one in \cite{ZZZ}. That is, we increase the number of marked points in the surface by one or two points to obtain a surface $(S,\M')$ where $\M\subset \M'$.  We triangulate $(S,\M')$ by adding one or two arcs and flip the - now bigger - overlap to an orthogonal overlap.  This gives rise to a triangle from $\gamma_1$ to $\gamma_2$ with middle terms $\gamma_5$ and $\gamma_6$ in $\C(S,\M')$ where here $\gamma_1$ and $\gamma_2$ are considered as arcs in $(S,\M')$.  Then by flipping the new (orthogonal) arc and using the cutting procedure  described in \cite{MP} which is compatible with Iyama-Yoshino reduction \cite{IY},
we obtain the triangle involving the arcs $\gamma_5$ and $\gamma_6$ in $\C(S,\M)$. 

If the crossing under consideration is an  arrow crossing then by (**) above we obtain a triangle
with middle terms corresponding to the arcs $\gamma_3$ and $\gamma_4$. The other triangle is then obtained by either flipping an 
arc in the triangulation and thus creating 
an overlap ({\it{i.e.}} a module crossing) which we can flip to an orthogonal overlap or, if this is not possible, by adding a marked point  
to obtain a surface $(S,\M')$ with $\M \subset \M'$ and completing it to a triangulation of $(S,\M')$. In which case we obtain 
an overlap which we can flip to an orthogonal overlap. This gives rise to a triangle in $\C(S,\M')$. 
By cutting according to \cite{MP}, we obtain the 
corresponding triangle from $\gamma_1$ to $\gamma_2$ in $\C(S,\M)$ with middle terms $\gamma_5$ and $\gamma_6$. 

If the crossing is a crossing in a 3-cycle then none of the triangles in $\C(S,\M)$ are obtained from non-split short exact sequences in
the Jacobian algebra corresponding to the given triangulation $T$. Instead we change the triangulation to create a crossing in a module - by possibly adding a marked point. Once we are in the case of a module crossing we can adapt the described strategy for module crossings above.

\subsubsection{Module crossing}

Here we consider the case that $M_1 = M(w_1)$ crosses $M_2 = M(w_2)$  in a module. In terms of snake graphs a crossing in a module corresponds to a crossing in an overlap  $\calg$. Therefore as explained above,  by Theorem 3.5 there always is a non-split triangle in $\C(S,\M)$ given by
 $$\gamma_2 \longrightarrow \gamma_3 \oplus \gamma_4 \longrightarrow \gamma_1 \longrightarrow \gamma_2[1]. $$
  When we give an orientation to an  arc $\gamma$, we call $s(\gamma)$, the marked point at which $\gamma$ starts and $t(\gamma)$ the marked point at which $\gamma$ ends.
In order to prove the existence of the triangle involving $\gamma_5$ and $\gamma_6$, 
there are several cases to consider depending on where the arcs $\gamma_1$ and $\gamma_2$ start and end with respect to the overlap. 

Let $w_1 = P_1 w S_1$ and $w_2 = P_2 w S_2$ where $w$ corresponds to the overlap  $\calg$.   Let $\tau_1, \tau_2, \ldots, \tau_n$ be the arcs corresponding to the  overlap $\calg$,
 that is $s(w) = \tau_1$ and $t(w) = \tau_n$. 

\emph{ Case 1:} $P_1 \neq 0, S_1 \neq 0, P_2 \neq 0, S_2 \neq 0$.

In $(S,\M)$ this corresponds to the local configuration as in Figure~\ref{FigureCase1ModuleCrossing}.

\begin{figure}[ht]
\begin{tikzpicture}
\node(a) at (0,7) [label={[label distance=-.3cm]100:$A$}, scale=.7] {$\bullet$};
\node(a') at (-1.5,7.5) [label={[label distance=-.3cm, scale=.7]100:$s(\gamma_1)$}, scale=.7] { };
\node(a'') at (1.7,3) [label={[label distance=-.3cm, scale=.7]300:$t(\gamma_1)$}, scale=.7] { };
\node(bb) at (1.5,7.5) [label={[label distance=-.3cm, scale=.7]100:$s(\gamma_2)$}, scale=.7] { };
\node(b'') at (-1.7,3) [label={[label distance=-.3cm, scale=.7]100:$t(\gamma_2)$}, scale=.7] { };

\node (f) at (2,6) [scale=.5] {$\bullet$};
\node (ff) at (-2,6) [scale=.5] {$\bullet$};
\node (g) at (2,4) [scale=.5] {$\bullet$};
\node (gg) at (-2,4) [scale=.5] {$\bullet$};
\node (h) at (0,3)  [label={[label distance=-.15cm]270:$B$}, scale=.7] {$\bullet$};

\node at (-1.2,6.7)[scale=.8]{$i$};
\node at (1.5,6.7)[scale=.8]{$j$};
\node at (-1.25,3.3)[scale=.8]{$k$};
\node at (1.25,3.3)[scale=.8]{$l$};

\draw (a.center) -- (f.center);
\draw (a.center) -- (ff.center);
\path (ff.center) edge node [pos=.5, fill=white,outer sep=1mm,scale=.7]{$s(w)$} (f.center);
\path (gg.center) edge node [pos=.5, fill=white,outer sep=1mm,scale=.7]{$t(w)$} (g.center);
\draw [dotted] (ff.center)--(gg.center) (f.center)--(g.center);
\draw (gg.center)--(h.center) (g.center)--(h.center);

\path[color=red, line width=1] (a'.center) edge node [pos=.3, rotate=115]{$<$} (a''.center);
\path[color=red, line width=1] (a'.center) edge node [pos=.6, fill=white,outer sep=1mm,scale=.7, bend right]{$\gamma_1$} (a''.center);
\path[color=tealblue, line width=1] (bb.center) edge node [pos=.3, rotate=65]{$<$} (b''.center);
\path[color=tealblue, line width=1] (bb.center) edge node [pos=.6, fill=white,outer sep=1mm,scale=.7, bend right]{$\gamma_2$} (b''.center);

\end{tikzpicture}
\caption{\emph{Case 1}: Local configuration of $\gamma_1$ and $\gamma_2$ crossing in a module where $A$ and $B$ might coincide. }\label{FigureCase1ModuleCrossing}
\end{figure}
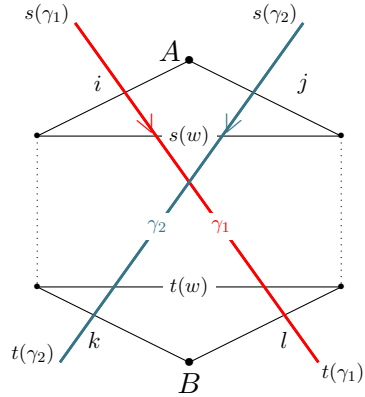

In  particular, $M_1$ crosses $M_2$ 
such that $t(P_1) =i$ and $t(P_2)
= j$ and $s(S_1) = l$ and $s(S_2) = k$.

 Since the crossing is not a self-crossing overlap, there exists an arc $\tau_{AB}$ from $A$ to $B$ crossing the arcs $\tau_1, \ldots, \tau_n$. Let $T'$ be the triangulation containing  $\tau_{AB}$ and such that the flips of $\tau_1, \ldots, \tau_n$ in 
$T'$ are connected to A and do not cross the arc $\tau_{AB}$. Locally  $\tau_{AB}$ lies in $T'$ as in Figure~\ref{AB}. 

\begin{figure}[ht]
\begin{tikzpicture}
\node (a) at (0,4) [label={[label distance=-.2cm]100:$A$}, scale=.7] {$\bullet$};
\node (b) at (0,0) [label={[label distance=-.1cm]270:$B$}, scale=.7] {$\bullet$};
\node (c) at (-1,1.5) [label={[label distance=-.2cm]100:$C$}, scale=.7] {$\bullet$};
\node (d) at (1.2,2.5) [label={[label distance=-.2cm]20:$D$}, scale=.7] {$\bullet$};
\path (a.center) edge[right] node[pos=.2, right]{$a$} (b.center);
\path (a.center) edge[left] node{$b$} (c.center);
\path (a.center) edge[right] node[pos=.4, right]{$c$} (d.center);
\path (c.center) edge[left] node[pos=.6, left]{$e$} (b.center);
\path (d.center) edge[right] node{$d$} (b.center);
\path[color=red, line width=1] (-1.2,4) edge node[pos=.43, rotate=115]{$<$} (1.2,0);
\path[color=red, line width=1] (-1.2,4) edge[left] node[pos=.45]{$\gamma_1$} (1.2,0);
\path[color=tealblue, line width=1] (1.2,4) edge node[pos=.43, rotate=65]{$<$} (-1.2,0);
\path[color=tealblue, line width=1] (1.2,4) edge[right] node[pos=.4]{$\gamma_2$} (-1.2,0);
\end{tikzpicture}
\caption{Triangulation $T'$ where $a$ denotes the arc $\tau_{AB}$. } \label{AB}
\end{figure}
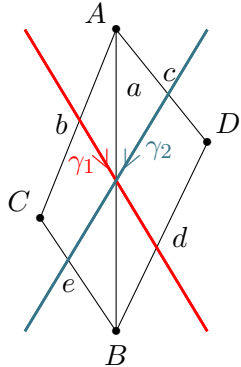

We denote by    $J(Q',W')$ the Jacobian 
algebra  with respect to  the new triangulation $T'$. Let $M'_1$ and $M'_2$ be the string modules over $J(Q',W')$  corresponding to the arcs $\gamma_1$ and $\gamma_2$. Now  $M'_2$ crosses $M'_1$ in a new overlap corresponding to $\tau_{AB}$.  We call this new overlap an \emph{orthogonal flip} of the overlap  $\calg$.

More explicitly, let  $w'_1$ (resp. $w'_2$) be the string of $\gamma_1$ (resp. $\gamma_2$) in $(S, \M, T')$. Then 
$w'_1$ contains the subword $b \longleftarrow a \longrightarrow d$ and 
$(w'_2)^{-1}$ contains the subword $e \longrightarrow a \longleftarrow c$ where, as in Figure~\ref{AB}, $a$ denotes the arc $\tau_{AB}$. This gives rise to  $M'_2 = M(w'_2)$  crossing $M'_1 = M(w'_1)$  in the simple module $M(a)$. 

By Theorem~\ref{ShortExactSequences}(i)  this gives rise to a non-split short exact sequence 
$$ 0 \longrightarrow M'_1 \longrightarrow   M'_3 \oplus M'_4 \longrightarrow M'_2 \longrightarrow 0$$
where $M'_3$ (resp. $M'_4$) is the string module over $J(Q',W')$ corresponding to the arc $\gamma_5$ (resp. $\gamma_6$) with respect to  $T'$.
Thus in $\C(S,\M)$ there is a triangle 
$$\gamma_1 \longrightarrow \gamma_5 \oplus \gamma_6 \longrightarrow  \gamma_2
\longrightarrow  \gamma_1[1].$$

\emph{ Case 2 (i):} $P_1 = 0, S_1 \neq 0, P_2 \neq 0, S_2 \neq 0$.

(a) Suppose $P_2$ is not a direct string. In this case we have the local  configuration as in Figure~\ref{FigureModuleCrossingNotDirectString}.

\begin{figure}[ht]
\begin{tikzpicture}
\node (aa) at (.3,8) [label={[label distance=-.2cm]100:$A$}, scale=.7] {$\bullet$};
\node(a) at (0,7) [label={[label distance=-.3cm]100:$s(\gamma_1)$}, scale=.7] {$\bullet$};
\node (b) at (2,8) [scale=.5] {$\bullet$};
\node (c) at (2.5,7.5) [scale=.5] {$\bullet$};
\node (d) at (2.8,7) [scale=.5] {$\bullet$};
\node (f) at (2,6) [scale=.5] {$\bullet$};
\node (ff) at (-2,6) [scale=.5] {$\bullet$};
\node (g) at (2,4) [scale=.5] {$\bullet$};
\node (gg) at (-2,4) [scale=.5] {$\bullet$};
\node (h) at (0,3)  [label={[label distance=-.15cm]270:$B$}, scale=.7] {$\bullet$};

\path (aa.center) edge node [pos=.5, fill=white,outer sep=1mm,scale=.7]{$k$} (b.center);
\path (a.center) edge node [pos=.5, fill=white,outer sep=1mm,scale=.6]{$f_1$} (b.center);
\path (a.center) edge node [pos=.5, fill=white,outer sep=1mm,scale=.6]{$f_2$} (c.center);
\path (a.center) edge node [pos=.5, fill=white,outer sep=1mm,scale=.6]{$f_3$} (d.center);
\draw (b.center) -- (c.center) -- (d.center); 
\draw (a.center)--(ff.center); 
\draw [loosely dotted] (d.center) -- (f.center);
\path (a.center) edge node [pos=.5, fill=white,outer sep=1mm,scale=.6]{$f_m$} (f.center);
\path (ff.center) edge node [pos=.5, fill=white,outer sep=1mm,scale=.7]{$s(w)$} (f.center);
\path (gg.center) edge node [pos=.5, fill=white,outer sep=1mm,scale=.7]{$t(w)$} (g.center);
\draw [loosely dotted] (ff.center)--(gg.center) (f.center)--(g.center) (1,8.6) -- (1,8);
\draw (gg.center)--(h.center) (g.center)--(h.center);

\path[color=red, thick] (a.center) edge[bend right] node[pos=.35, scale=.9, rotate=110]{$<$} (2,2.5);
\path[color=red, thick] (a.center) edge[bend right] node[pos=.5, fill=white,outer sep=1mm,scale=.6]{$\gamma_1$} (2,2.5);
\path[color=tealblue, thick] (1.3,8.8) edge[bend left] node[pos=.35, scale=.9, rotate=75]{$<$} (-2,2.5);
\path[color=tealblue, thick] (1.3,8.8) edge[bend left] node[pos=.48, fill=white,outer sep=1mm,scale=.6]{$\gamma_2$} (-2,2.5);

\end{tikzpicture}
\caption{\emph{ Case 2 (i) a):} Local configuration of $\gamma_1$ and $\gamma_2$ where $A$, $s(\gamma_1)$ and $s(\gamma_2)$ may coincide.}
\label{FigureModuleCrossingNotDirectString}
\end{figure}
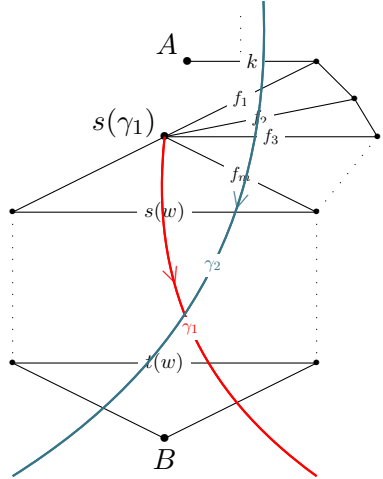

Since $P_2$ is not a direct string, there must exist a marked point $A$ and an arc $k$ and a maximal fan $f_1, \ldots, f_n$ such that $\gamma_2$ crosses $k$ and $f_1, \ldots, f_n$, as in Figure~\ref{FigureModuleCrossingNotDirectString}. Let $B$ be the marked point in the triangle $B s(\tau_n) t(\tau_n) $ where $\gamma_2$ crosses the arc from $s(\tau_n)$ to $B$, as in Figure~\ref{FigureModuleCrossingNotDirectString}. 
Then there exists a triangulation $T'$ of $(S,\M)$ containing the arc $\tau_{AB}$ starting at $A$, ending at $B$,  crossing both $\gamma_1$ and $\gamma_2$ and such that locally
no other arc in $T'$ crosses both $\gamma_1$ and $\gamma_2$.

Therefore by Theorem~\ref{ShortExactSequences}(i) we obtain a non-split  short exact sequence 
$$ 0 \longrightarrow M'_1 \longrightarrow   M'_3 \oplus M'_4 \longrightarrow M'_2 \longrightarrow 0$$
where $M'_3$ (resp. $M'_4$) is the string module over $J(Q',W')$ corresponding to the arc $\gamma_5$ (resp. $\gamma_6$) with respect to  $T'$.
Thus in $\C(S,\M)$ there is a triangle 
$$\gamma_1 \longrightarrow \gamma_5 \oplus \gamma_6 \longrightarrow  \gamma_2
\longrightarrow  \gamma_1[1].$$

(b) Suppose now that $P_2$ is a direct string, see Figure~\ref{overlap direct string}. 

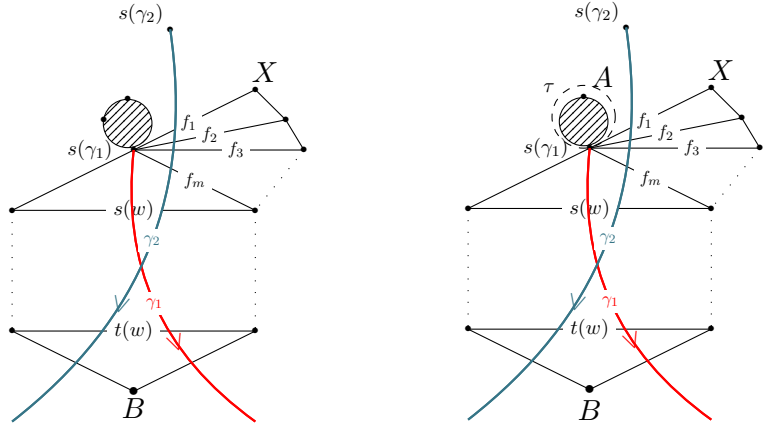
\begin{figure}[ht]
\begin{tikzpicture}
 \node at (0,0){$\begin{tikzpicture}[scale=.8]
\node (aaa) at (.6,9) [label={[label distance=-.2cm, scale=.7 ]100:$s(\gamma_2)$}, scale=.5] {$\bullet$};
\node (aa) at (.3,8) [label={[label distance=-.2cm]100: }, scale=.7] { };
\node(a) at (0,7) [scale=.7] [label={[label distance=.1cm, scale=.7]180:$s(\gamma_1)$}, scale=.7] {$\bullet$};
\node (b) at (2,8) [scale=.7] [label={[label distance=-.3cm]10:$X$}, scale=.7] {$\bullet$};
\node (c) at (2.5,7.5) [scale=.5] {$\bullet$};
\node (d) at (2.8,7) [scale=.5] {$\bullet$};
\node (f) at (2,6) [scale=.5] {$\bullet$};
\node (ff) at (-2,6) [scale=.5] {$\bullet$};
\node (g) at (2,4) [scale=.5] {$\bullet$};
\node (gg) at (-2,4) [scale=.5] {$\bullet$};
\node (h) at (0,3)  [label={[label distance=-.2cm]-90:$B$}, scale=.7] {$\bullet$};

\node (B) at (-.1,7.85) [scale=.7] [label={[label distance=0cm]40: }, scale=.7] {$\bullet$};
\node (C) at (-.49, 7.5) [scale=.5] {$\bullet$};
\draw[fill=white,pattern=north east lines] (-.1,7.44) circle (0.4);

\path (a.center) edge node [pos=.45, fill=white,outer sep=1mm,scale=.6]{$f_1$} (b.center);
\path (a.center) edge node [pos=.5, fill=white,outer sep=1mm,scale=.6]{$f_2$} (c.center);
\path (a.center) edge node [pos=.6, fill=white,outer sep=1mm,scale=.6]{$f_3$} (d.center);
\draw (b.center) -- (c.center) -- (d.center);
\draw (a.center)--(ff.center);
\draw [loosely dotted] (d.center) -- (f.center);
\path (a.center) edge node [pos=.5, fill=white,outer sep=1mm,scale=.6]{$f_m$} (f.center);
\path (ff.center) edge node [pos=.5, fill=white,outer sep=1mm,scale=.7]{$s(w)$} (f.center);
\path (gg.center) edge node [pos=.5, fill=white,outer sep=1mm,scale=.7]{$t(w)$} (g.center);
\draw [loosely dotted] (ff.center)--(gg.center) (f.center)--(g.center); 
\draw (gg.center)--(h.center) (g.center)--(h.center);

\path[color=red, thick] (a.center) edge[bend right] node[pos=.65, scale=.9, rotate=110]{$<$} (2,2.5);
\path[color=red, thick] (a.center) edge[bend right] node[pos=.5, fill=white,outer sep=1mm,scale=.6]{$\gamma_1$} (2,2.5);
\path[color=tealblue, thick] (aaa.center) edge[bend left] node[pos=.65, scale=.9, rotate=75]{$<$} (-2,2.5);
\path[color=tealblue, thick] (aaa.center) edge[bend left] node[pos=.48, fill=white,outer sep=1mm,scale=.6]{$\gamma_2$} (-2,2.5);
\end{tikzpicture}$};

\node at (6,0){$\begin{tikzpicture}[scale=.8]
\node (aaa) at (.6,9) [label={[label distance=-.2cm, scale=.7 ]100:$s(\gamma_2)$}, scale=.5] {$\bullet$};
\node (aa) at (.3,8) [label={[label distance=-.2cm]100: }, scale=.7] { };
\node(a) at (0,7) [scale=.7] [label={[label distance=.1cm, scale=.7]180:$s(\gamma_1)$}, scale=.7] {$\bullet$};
\node (b) at (2,8) [scale=.7] [label={[label distance=-.3cm]10:$X$}, scale=.7] {$\bullet$};
\node (c) at (2.5,7.5) [scale=.5] {$\bullet$};
\node (d) at (2.8,7) [scale=.5] {$\bullet$};
\node (f) at (2,6) [scale=.5] {$\bullet$};
\node (ff) at (-2,6) [scale=.5] {$\bullet$};
\node (g) at (2,4) [scale=.5] {$\bullet$};
\node (gg) at (-2,4) [scale=.5] {$\bullet$};
\node (h) at (0,3)  [label={[label distance=-.15cm]270:$B$}, scale=.7] {$\bullet$};

\node (B) at (-.1,7.85) [scale=.7] [label={[label distance=-.15 cm]10: $A$}, scale=.7] {$\bullet$};

\node at (-.67,8) [scale=.7] {$\tau$};

\draw[dashed] (-.1,7.52) circle (0.52);

\draw[fill=white,pattern=north east lines] (-.1,7.44) circle (0.4);

\path (a.center) edge node [pos=.45, fill=white,outer sep=1mm,scale=.6]{$f_1$} (b.center);
\path (a.center) edge node [pos=.5, fill=white,outer sep=1mm,scale=.6]{$f_2$} (c.center);
\path (a.center) edge node [pos=.6, fill=white,outer sep=1mm,scale=.6]{$f_3$} (d.center);
\draw (b.center) -- (c.center) -- (d.center);
\draw (a.center)--(ff.center); 
\draw [loosely dotted] (d.center) -- (f.center);
\path (a.center) edge node [pos=.5, fill=white,outer sep=1mm,scale=.6]{$f_m$} (f.center);
\path (ff.center) edge node [pos=.5, fill=white,outer sep=1mm,scale=.7]{$s(w)$} (f.center);
\path (gg.center) edge node [pos=.5, fill=white,outer sep=1mm,scale=.7]{$t(w)$} (g.center);
\draw [loosely dotted] (ff.center)--(gg.center) (f.center)--(g.center); 
\draw (gg.center)--(h.center) (g.center)--(h.center);

\path[color=red, thick] (a.center) edge[bend right] node[pos=.65, scale=.9, rotate=110]{$<$} (2,2.5);
\path[color=red, thick] (a.center) edge[bend right] node[pos=.5, fill=white,outer sep=1mm,scale=.6]{$\gamma_1$} (2,2.5);
\path[color=tealblue, thick] (aaa.center) edge[bend left] node[pos=.65, scale=.9, rotate=75]{$<$} (-2,2.5);
\path[color=tealblue, thick] (aaa.center) edge[bend left] node[pos=.48, fill=white,outer sep=1mm,scale=.6]{$\gamma_2$} (-2,2.5);

%

\end{tikzpicture}$};
\end{tikzpicture}
\caption{\emph{ Case 2 (i) b):} Left figure: local configuration in $(S, \M, T)$ where the boundary component 
 has several marked points. \newline
Right figure: 
local configuration in $(S,\M, T)$ where $s(\gamma_1)$ is the only marked point on the boundary component. 
}
 \label{overlap direct string}
\end{figure}


Our argument is based on the boundary component $\mathcal{B}$ containing $s(\gamma_1)$.

(I) Suppose that $\mathcal{B}$ contains another marked point $A$ which is not equal to $s(\gamma_1)$ and $s(\gamma_2)$ is not in $\mathcal{B}$. Then there is an arc $\tau_{AB} $ from $A$ to $B$ 
crossing $\gamma_1$ and $\gamma_2$, see left hand side of Figure~\ref{overlap direct string} and we conclude as in part (a). 

(II) Suppose $s(\gamma_1)$ is the only marked point on $\mathcal{B}$ and that $s(\gamma_1) \neq s(\gamma_2)$. Then consider instead the surface $(S,\M')$ where $\M' \supset \M$ has exactly one more marked point $A$ than $\M$ lying on 
the boundary component $\mathcal{B}$, see right hand side of Figure~\ref{overlap direct string}. Complete $T$ to a triangulation $T'$ on $(S,\M')$ by adding one new arc $\tau$ from 
$A$ to $X$, where $X$ is as in figure~\ref{overlap direct string}.
 Therefore the same method as in part (a) can be applied and we obtain a triangle in 
$\C{(S,\M')}$ 
$$\gamma_1 \longrightarrow \gamma_5 \oplus \gamma_6 \longrightarrow  \gamma_2
\longrightarrow  \gamma_1[1].$$

Note that since $\M \subset \M'$, whenever we have an arc in $(S,\M')$ between marked points $a, b \in \M'$ such that $a, b \in \M$, by a slight abuse of notation we use the same notation for this arc as an arc in $(S,\M)$ and as an arc in $(S,\M')$.

Now flip the triangulation  $T'$ to a triangulation $T''$ such that $T''$ contains the arc $\tau'$ around the boundary component $\mathcal{B}$ from $s(\gamma_1)$ to $s(\gamma_1)$,
 see right hand side of Figure~\ref{overlap direct string}. Then there is a triangle $s(\gamma_1) s(\gamma_1) A$ in $T''$. 
Cutting $\tau'$ as defined in \cite{MP} gives a surface isotopic to $(S,\M)$ since  we delete any component homeomorphic to a triangle after
the cut. 

Note that by \cite{MP} the arcs in $(S,\M')$ not crossing $\tau'$ are in bijection with the arcs in 
$(S,\M)$. Since $\tau'$ is a boundary segment in $(S,\M)$ the corresponding 
object in $\C(S,\M)$ is the zero object.  Thus by   Proposition 5 in \cite{MP} we obtain a triangle in $\C(S,\M)$
$$\gamma_1 \longrightarrow  \gamma_5 \oplus \gamma_6 \longrightarrow  \gamma_2
\longrightarrow  \gamma_1[1].$$

(III) Suppose that $s(\gamma_2)$ is in $\mathcal{B}$ and that there is a marked point $A$ in $\mathcal{B}$ between $s(\gamma_1)$ and $s(\gamma_2)$ such that the arcs $f_1$ and 
$s(\gamma_1) A$ are two sides of a triangle, where $f_1, \ldots, f_n$ is a fan as in left hand side of Figure~\ref{overlap direct string}. Then we conclude as in part (a). 
If no such point $A$ exists then $s(\gamma_1) s(\gamma_2)$ is a boundary segment.  We add a marked point $A$ on this boundary arc and argue as in case (II). 


(IV) Suppose that $s(\gamma_1) = s(\gamma_2)$.  Then there is a triangulation $T'$ containing an arc from  some point $A$ on $\mathcal{B}$ (we add the point $A$ if it does not already exist)
  to $B$ crossing $\gamma_1$ and $\gamma_2$ and such that locally no other 
arc in $T'$ crosses both $\gamma_1$ and $\gamma_2$ and we conclude as above.

\emph{Case 2 (ii):} $P_1 \neq 0$, $S_1 = 0$ and $P_2 \neq 0$, $S_2 \neq 0$ follows 
from Case 2(i) by
changing the orientation of $\gamma_1$ and $\gamma_2$. 

\emph{Case 2 (iii) and (iv):} The case $P_2 = 0$ and $P_1, S_1, S_2$ non-zero and the case
$S_2 = 0 $ and $P_1, S_1, P_2$ non-zero follow by similar arguments as above. 

\emph{Case 3 (i):} $P_1 = 0$ and $S_2 = 0$ and $S_1$ and $P_2$ non-zero. 

(a) Suppose that neither $P_2$ nor   $S_1$ is  a direct string.  By a similar argument as in Case 2(i)(a) there are marked points $A$ 
and $B$ such that there is a triangulation $T'$ of $(S,\M)$ containing an arc 
corresponding to the arc $\tau_{AB}$ and we obtain a triangle in $\C(S,\M)$
$$\gamma_1 \longrightarrow \gamma_5 \oplus \gamma_6 \longrightarrow  \gamma_2
\longrightarrow  \gamma_1[1].$$

(b) Suppose that $P_2$ is a direct string and that $S_1$ is not a direct string or that  $P_2$ is not a direct string and that $S_1$ is a direct string. Then 
we use a similar argument as in case 
2(i)(b) above to obtain a triangle in $\C(S,\M)$
$$\gamma_1 \longrightarrow \gamma_5 \oplus \gamma_6 \longrightarrow  \gamma_2
\longrightarrow  \gamma_1[1].$$

(c) Suppose both $P_2$ and $S_1$ are direct strings. Then if   $s(\gamma_1) s(\gamma_2)$ and 
$t(\gamma_1) t(\gamma_2)$  are not both boundary segments then the argument is a combination of the above arguments. 

Note that if  $s(\gamma_1) s(\gamma_2)$ and 
$t(\gamma_1) t(\gamma_2)$  are boundary segments then we obtain a trivial triangle. 
Namely, consider a surface $(S,\M')$ where $\M'$ contains 2 more marked points than $\M$, 
one in each of the boundary segments $s(\gamma_1) s(\gamma_2)$ and 
$t(\gamma_1) (\gamma_2)$. As before this gives 
rise to a triangle in $\C{(S,\M')}$
$$\gamma_1 \longrightarrow \gamma_5 \oplus \gamma_6 \longrightarrow  \gamma_2
\longrightarrow  \gamma_1[1].$$

Applying the  construction in \cite{MP} and cutting twice, we obtain the trivial triangle in  
$\C(S,\M)$
$$\gamma_1 \longrightarrow \gamma_1\longrightarrow 0
\longrightarrow  \gamma_2.$$


\emph{Case 3 (ii):} $P_1 = 0$ and $P_2 = 0$ and $S_1$ and $S_2$ non-zero, follows
from the above by changing the orientation of $\gamma_1$ and $\gamma_2$. 

\emph{Case 3 (iii) and (iv):} The case $P_1 = 0$ and $P_2 = 0$ and $S_2$ and $P_2$ non-zero, 
and the case $P_2 = 0$ and $S_2 = 0$ and $P_1$ and $S_1$ non-zero, follow
by similar arguments to the above.

\subsubsection{Arrow crossing}

Here we consider the case that $M_1$ crosses $M_2$  in an arrow.

\emph{Case 1:} Suppose the crossing occurs in an inner triangle of $T$. Let $\tau$ be the 
arc corresponding to the segment $s(\gamma_2) t(\gamma_1)=BC$, see Figure~\ref{FigureArrowCrossingCase1-2}  (a).

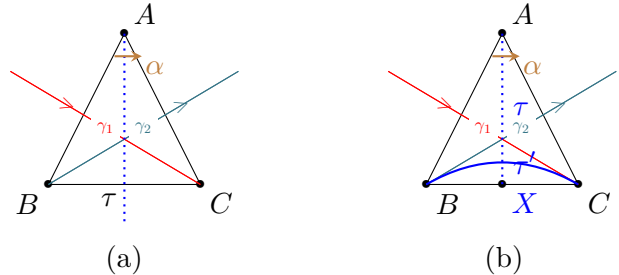
\begin{figure}[ht]
\begin{tikzpicture}
\node at (0,0){$\begin{tikzpicture}[trans/.style={thick,->,shorten >=2pt,shorten <=2pt,>=stealth}]
\node (a) at (1,2) [label={[label distance=-.2cm]30:$A$}, scale=.7] {$\bullet$};
\node (b) at (0,0) [label={[label distance=-.2cm]-150:$B$}, scale=.7] {$\bullet$}; 
\node (c) at (2,0) [label={[label distance=-.2cm]-60:$C$}, scale=.7] {$\bullet$};
\draw (a.center)--(b.center)--(c.center)--(a.center);
\path[color=tealblue] (b.center) edge node[pos=.7, scale=.7, rotate=20]{$>$} (2.5,1.5);
\path[color=tealblue] (b.center) edge node[pos=.5, fill=white,outer sep=1mm, scale=.6]{$\gamma_2$} (2.5,1.5);
\path[color=red] (c.center) edge node[pos=.7, scale=.7, rotate=-20]{$>$} (-.5,1.5);
\path[color=red] (c.center) edge node[pos=.5, fill=white,outer sep=1mm, scale=.6]{$\gamma_1$} (-.5,1.5);
\path[thick, dotted, color=blue] (a.center) edge[right] (1,-.5);
\draw[trans, color=brown] (.8,1.7)--(1.3,1.7);
\node[color=brown] at (1.4,1.55){$\alpha$};
\node at (.8,-.2) {$\tau$};
\end{tikzpicture}$};
\node at (0,-2){(a)};

\node at (5,0){$\begin{tikzpicture}[trans/.style={thick,->,shorten >=2pt,shorten <=2pt,>=stealth}]
\node (a) at (1,2) [label={[label distance=-.2cm]30:$A$}, scale=.7] {$\bullet$};
\node (b) at (0,0) [label={[label distance=-.2cm]-30:$B$}, scale=.7] {$\bullet$}; 
\node (c) at (2,0) [label={[label distance=-.2cm]-60:$C$}, scale=.7] {$\bullet$};
\node(d) at (1,0) [label={[color=blue, label distance=-.2cm]-60:$X$}, scale=.7] {$\bullet$};
\draw (a.center)--(b.center)--(c.center)--(a.center);
\path[color=tealblue] (b.center) edge node[pos=.7, scale=.7, rotate=20]{$>$} (2.5,1.5);
\path[color=tealblue] (b.center) edge node[pos=.5, fill=white,outer sep=1mm, scale=.6]{$\gamma_2$} (2.5,1.5);
\path[color=red] (c.center) edge node[pos=.7, scale=.7, rotate=-20]{$>$} (-.5,1.5);
\path[color=red] (c.center) edge node[pos=.5, fill=white,outer sep=1mm, scale=.6]{$\gamma_1$} (-.5,1.5);\path[thick, dotted, color=blue] (a.center) edge[right] node{$\tau$} (d.center);
\path[thick, color=blue] (b.center) edge[bend left, right] node{$\tau'$} (c.center);
\draw[trans, color=brown] (.8,1.7)--(1.3,1.7);
\node[color=brown] at (1.4,1.55){$\alpha$};
\end{tikzpicture}$};
\node at (5,-2){(b)};
\end{tikzpicture}
\caption{(a) Arrow crossing in an inner triangle, (b) Arrow crossing in a triangle where $BC$ is a boundary segment.} 
\label{FigureArrowCrossingCase1-2}
\end{figure}

Flipping $\tau$ to $\tau'$ in 
its quadrilateral gives rise to an overlap  given by $\tau'$ and we use the module crossing methods above to 
obtain a triangle in $\C(S,\M)$ 
$$\gamma_1 \longrightarrow \gamma_5 \oplus \gamma_6 \longrightarrow  \gamma_2
\longrightarrow  \gamma_1[1].$$

We remark that the points $A, B$ and  $C$ are not necessarily distinct.  
 
\emph{Case 2:} Suppose the crossing occurs in a triangle of $T$ where the 
segment $s(\gamma_2) t(\gamma_1)=BC$ is a boundary segment, see Figure~\ref{FigureArrowCrossingCase1-2}  (b).
%
%
%
%

Consider the surface $(S,\M')$ where $\M' = \M \cup \{X\}$  and $X$ lies on the boundary segment $BC$, see Figure~\ref{FigureArrowCrossingCase1-2}  (b).
 We complete $T$ to a triangulation of $(S,\M')$ by adding an arc $\tau$ corresponding to the segment $AX$. This
gives rise to an overlap given by  $\tau$. Again we use the module crossing methods 
above to obtain a triangle in $\C(S,\M')$ 
$$\gamma_1 \longrightarrow \gamma_5 \oplus \gamma_6 \longrightarrow  \gamma_2
\longrightarrow  \gamma_1[1].$$
We apply the construction in \cite{MP} and cut along $\tau'$. Then $\gamma_6$ corresponds to a boundary segment and as above
 we obtain a triangle in $\C(S,\M)$ 
$$\gamma_1 \longrightarrow \gamma_5  \longrightarrow  \gamma_2
\longrightarrow  \gamma_1[1].$$

Note that the points $B$ and $C$ may coincide.

\subsubsection{3-cycle crossing}

Here we consider the case that $M_1$ crosses $M_2$  in a 3-cycle.

\emph{Case 1:} Suppose that $s(\gamma_1) \neq s(\gamma_2)$, see Figure~\ref{3-cycle crossing}.

\begin{figure}[ht]
\begin{tikzpicture}[trans/.style={thick,->,shorten >=2pt,shorten <=2pt,>=stealth}]
\node (a) at (1,2) [label={[label distance=-.2cm]30:$A$}, scale=.7] {$\bullet$};
\node (b) at (0,0) [label={[label distance=-.2cm]-30:$B$}, scale=.7] {$\bullet$}; 
\node (c) at (2,0) [label={[label distance=-.2cm]-60:$C$}, scale=.7] {$\bullet$};
\node (d) at (-1,1) [label={[label distance=-.2cm, color=red]180:$\ldots$}, scale=.7] { };
\node (e) at (2.5, 1.2) [label={[label distance=-.2cm, color=red]0:$\ldots$}, scale=.7] { };
\node (f) at (.9, -1.2) [label={[label distance=-.25cm, color=tealblue] -3:$\vdots$}, scale=.7] { };
\node (g) at (-.3,.4) [label={[label distance=-.7cm, color=blue]0:flip}, scale=.7] { };
\node (h) at (-.4,1.7) [label={[label distance=-.2cm, color=blue]180:$X$}, scale=.7] {$\bullet$ };

\path[line width=.7, color=red] (d.center) edge node[pos=.75, scale=.7]{$>$} (e.center);
\path[line width=.7, color=red] (d.center) edge node[pos=.2, fill=white,outer sep=1mm, scale=.6]{$\gamma_1$} (e.center);

\draw (a.center)--(b.center)--(c.center)--(a.center);

\path[line width=.7, color=tealblue] (a.center) edge node[pos=.7, scale=.7, rotate=90]{$<$} (f.center);
\path[line width=.7, color=tealblue] (a.center) edge node[pos=.5, fill=white,outer sep=1mm, scale=.6]{$\gamma_2$} (f.center);
\draw[color=blue] (h.center) -- (c.center);
\draw[trans, color=blue] (g.center) -- (.3, .4);
\end{tikzpicture}
\caption{\emph{Case 1:} 3-cycle crossing.}
\label{3-cycle crossing}
\end{figure}
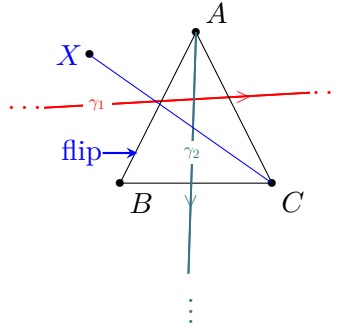

Remark that $s(\gamma_1)$ and $s(\gamma_2)$ may or may not lie in the same boundary component of $(S,\M)$. 
In either case we can add a marked point $X$ to obtain a surface $(S,\M')$ such that the segment $XC$ 
 lies between the segments $AC$ and $BC$. We flip the arc corresponding to the segment
$AB$. This gives rise to a crossing of $\gamma_1$ and $\gamma_2$ with overlap corresponding to $XC$. 
Applying the module crossing methods described above, we obtain two triangles in $\C(S,\M')$. By   \cite{MP} we 
obtain the corresponding  triangles in $\C(S,\M)$.

\emph{Case 2:} Suppose that $s(\gamma_1) = s(\gamma_2)$.

Denote by $P = \gamma_1(t_0) $ for some $t_0  \in [0, 1]$ the intersection of $\gamma_1$ and $\gamma_2$ corresponding to the crossing of $\gamma_1$ and $\gamma_2$ under consideration. Without loss of generality, assume that $P$ is also equal to $ \gamma_2(t_0)$.  Consider the closed curve $\sigma = \sigma[0,1]$ which is a union of the segments
 $\sigma[0,t_0]$  and $\sigma[t_0,1]$ such that $\sigma[0]=\sigma[1]$ and where 
$ \sigma[0,t_0] = \gamma_1[0, t_0]$ and $\sigma[t_0,1] = \gamma_2^{-1}[t_0, 0]$. Then $\sigma$ cannot be homotopic to a point, 
since otherwise there would not be a 3-cycle crossing.  Therefore the local configuration in this case corresponds to one of the two cases illustrated in Figure~\ref{3-cycle crossing type 2i}.

\begin{figure}[ht]
\begin{tikzpicture}
\node at (0,0){$\begin{tikzpicture}
\node (a) at (1,2) [label={[label distance=-.2cm]30:$A$}, scale=.7] {$\bullet$};
\node (b) at (0,0) [label={[label distance=-.2cm]-30:$B$}, scale=.7] {$\bullet$}; 
\node (c) at (2,0) [label={[label distance=-.2cm]-60:$C$}, scale=.7] {$\bullet$};
\node (f) at (1, -1.2) [label={[label distance=-.27cm, color=tealblue]-90:$\vdots$}, scale=.7] { };
\coordinate (A) at (-.6,1.2);

\coordinate (C) at (-.9,1.35);
\coordinate (C') at (-.9,1.2);
\coordinate (C'') at (-.9,1.05);

\draw[fill=white,pattern=north east lines] (0,1.2) circle (0.3);
\draw (a.center)--(b.center)--(c.center)--(a.center);

\path[line width=.7, color=tealblue] (a.center) edge node[pos=.7, scale=.7, rotate=90]{$<$} (f.center);
\path[line width=.7, color=tealblue] (a.center) edge node[pos=.34, fill=white,outer sep=1mm, scale=.6]{$\gamma_2$} (f.center);

\draw[color=red, line width=.7] (a.center) to [out=180, in=80] (-1.3,1.2) node[left, scale=.6]{$\gamma_1$}
to [out=-80, in=180] (3,1)node[color=red, right]{$\cdots$};

\draw[color=red, line width=.75] (a.center) to [out=180, in=80] (-1.3,1.2) node[scale=.6, rotate=-90]{$>$} to [out=-80, in=180] (3,1)node[color=red, right]{$\cdots$};

\draw[color=blue, dashed, line width=.7] (c)  to [out=120, in=50] (C);
\node[color=blue, scale=.9] at (C'.center){$\tau$};
\draw[color=blue, dashed, line width=.7] (C'') to [out=-50, in=150] (c);

\end{tikzpicture}
$};
\node at (7,0){$\begin{tikzpicture}
\node (a) at (1,2) [label={[label distance=-.2cm]30:$A$}, scale=.7] {$\bullet$};
\node (b) at (0,0) [label={[label distance=-.2cm]-30:$B$}, scale=.7] {$\bullet$}; 
\node (c) at (2,0) [label={[label distance=-.2cm]-60:$C$}, scale=.7] {$\bullet$};
\node (f) at (1, -1.2) [label={[label distance=-.27cm, color=tealblue]-90:$\vdots$}, scale=.7] { };
\coordinate (A) at (-.7,1.2);
\coordinate (A') at (-.1,1.2);
\coordinate (B) at (-.55,1.15);
\coordinate (B') at (-.25,1.15);

\coordinate (C) at (-.9,1.35);
\coordinate (C') at (-.9,1.2);
\coordinate (C'') at (-.9,1.05);

\draw (A) to [out=-30, in=210] (A');
\draw (B) to [out=30, in=150] (B');

\draw (a.center)--(b.center)--(c.center)--(a.center);

\path[line width=.7, color=tealblue] (a.center) edge node[pos=.7, scale=.7, rotate=90]{$<$} (f.center);
\path[line width=.7, color=tealblue] (a.center) edge node[pos=.35, fill=white,outer sep=1mm, scale=.6]{$\gamma_2$} (f.center);

\draw[color=red, line width=.7] (a.center) to [out=180, in=80] (-1.3,1.2) node[left, scale=.6]{$\gamma_1$}
to [out=-80, in=180] (3,1)node[color=red, right]{$\cdots$};

\draw[color=red, line width=.75] (a.center) to [out=180, in=80] (-1.3,1.2) node[scale=.6, rotate=-90]{$>$} to [out=-80, in=180] (3,1)node[color=red, right]{$\cdots$};

\draw[color=blue, dashed, line width=.7] (c)  to [out=120, in=50] (C);
\node[color=blue, scale=.9] at (C'.center){$\tau$};
\draw[color=blue, dashed, line width=.7] (C'') to [out=-50, in=150] (c);

\end{tikzpicture}
$};
\end{tikzpicture}
\caption{Possible 3-cycle crossing: $A= s(\gamma_1) = s(\gamma_2)$.}\label{3-cycle crossing type 2i}
\end{figure}
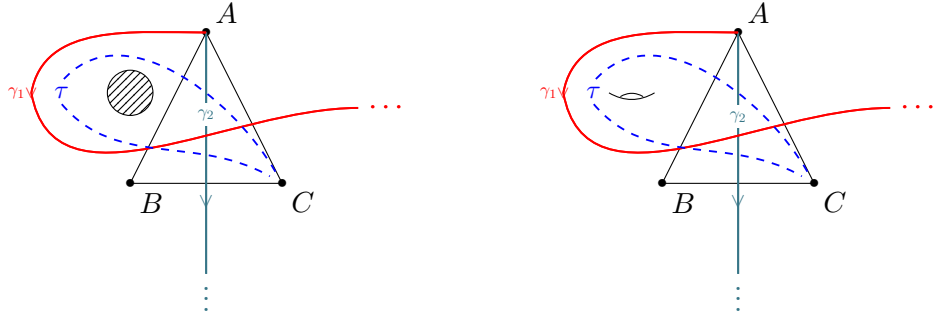

 Consider now a curve $\tau = \tau[0,1]$ which is a union of the segments
 $\tau[0,t_0]$, $\tau[t_0,t_1]$ and $\tau[t_1,1]$ where 
$ \tau[0,t_0] = CP$, $\tau[t_0,t_1]$ is a non-contractible non-selfcrossing curve such that $\tau[t_0] = P = \tau[t_1]$ and $\tau[t_1,1] = CP$. Note
that we choose $CP$ to be a curve with no self-intersection. 
Then $\tau[0,1]$  is homotopic to a closed curve on $S$ without any self-intersection in its interior, see Figure~\ref{3-cycle crossing type 2i}.

There exists a triangulation $T'$  containing the arc corresponding to $\sigma$. 
Thus $\gamma_1$ and $\gamma_2$ cross in a non self-crossing overlap containing 
at least $\sigma$. The rest follows as in \emph{Case 1}  for module crossings.   

The points $A, B$ and $C$ may coincide two by two. 

This completes the proof of Theorem~\ref{Triangles}.
\hfill $\Box$

\section{Example}\label{ExampleSection}

Let $J$ be the Jacobian algebra corresponding to the triangulation in Figure~\ref{Example}. 
We see that the example contains two crossing arcs $\gamma_1$ and $\gamma_2$ such that
the corresponding string modules $M_1$ and $M_2$ cross four times and such that $M_1$ has one self-crossing. Each type of crossing (in a module, see crossings \scalebox{1.5}{$\circled{1}$}, 
\scalebox{1.5}{$\circled{2}$} and \scalebox{1.5}{$\circled{5}$}, in an arrow, see crossing \scalebox{1.5}{$\circled{4}$}, in a 3-cycle, see crossing \scalebox{1.5}{$\circled{3}$}) 
occurs at least once. 
We remark that there is a module crossing in both directions, that is $M_1$ crosses $M_2$ in a module,   see 
crossing \scalebox{1.5}{$\circled{1}$}, 
and $M_2$ crosses $M_1$  in a module, see crossing \scalebox{1.5}{$\circled{2}$}. 

\begin{figure}[ht]
\includegraphics[scale=1.2]{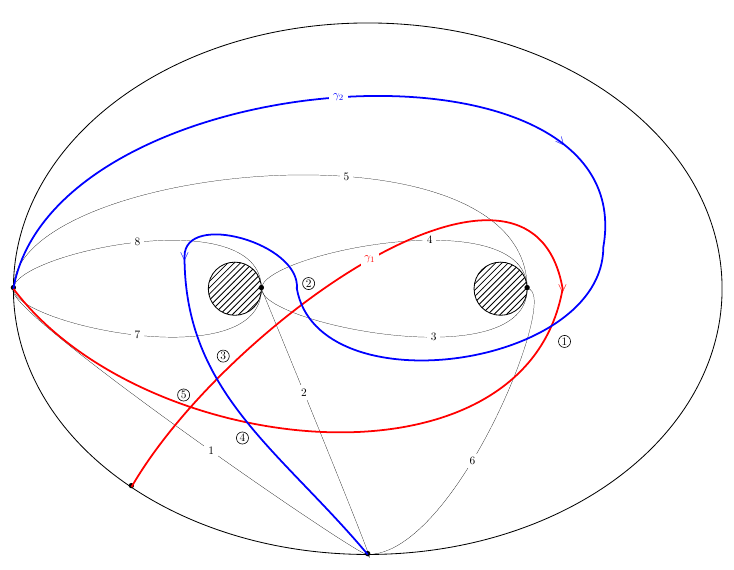}
\caption{Triangulation of a surface with two arcs crossing each other 4 times and with one arc with one self-crossing.}\label{Example}
\end{figure}

The arc $\gamma_1$ corresponds to the string module $M_1 =  M(w_1)$
and $\gamma_2$ corresponds to $M_2 = M(w_2)$ where  
$$ w_1=  \begin{tikzpicture}[baseline=-0.85ex]
\node at (0,0){1};
\node at (.25,-.25){2};
\node at (.5,0){3};
\node at (.75,.25){4};
\node at (1,0){5};
\node at (1.25,-.25){6};
\node at (1.5,0){2};
\end{tikzpicture}   \mbox{ \hspace{.5cm} and \hspace{.5cm} } w_2=\begin{tikzpicture}[baseline=-0.85ex]
\node at (0,0){6};
\node at (.25,-.25){3};
\node at (.5,0){4};
\node at (.75,.25){8};
\node at (1,0){7};
\end{tikzpicture}.$$

For each of the five crossings  \scalebox{1.5}{$\circled{1}$}-\scalebox{1.5}{$\circled{5}$}, we now explicitly give the modules $M_3 = M(w_3), M_4=M(w_4), M_5=M(w_5)$ and $M_6=M(w_6)$ corresponding to the smoothing of each of the given 
crossings as defined in Definition~\ref{DefinitionOfModulesInResolution}(1). Note that $M_3, \ldots, M_6$ depend on (and therefore change with) the given crossing whereas
$M_1$ and $M_2$ do not. Similarly, in terms of the corresponding arcs $\gamma_1 = \gamma(w_1), \ldots, \gamma_6 = \gamma(w_6)$, the arcs $\gamma_1$ and $\gamma_2$ are fixed whereas
$\gamma_3, \ldots, \gamma_6$ depend on the crossing under consideration.

Crossing \scalebox{1.5}{$\circled{1}$}: $M_1$ crosses $M_2$ in the simple module $\mathbf{6}$  and

\[
\begin{tikzpicture}
\node at (-1.75,0){$w_3=$};
\node at (0,0){\begin{tikzpicture}[baseline=-0.85ex]
	\node at (0,0){1};
	\node at (.25,-.25){2};
	\node at (.5,0){3};
	\node at (.75,.25){4};
	\node at (1,0){5};
	\node at (1.25,-.25){\bf 6};
	\node at (1.5,-.5){3};
	\node at (1.75,-.25){4};
	\node at (2,0){8};
	\node at (2.25,-.25){7};
	\end{tikzpicture}};
\node at (2.5,0){$w_4=$};
\node at (3.5,0){\begin{tikzpicture}[baseline=-0.85ex]
	\node at (.25,-.25){\bf 6};
	\node at (.5,0){2};
	\end{tikzpicture}};
\node at (5,0){$w_5=$};
\node at (6,0){\begin{tikzpicture}[baseline=-0.85ex]
	\node at (0,0){1};
	\node at (.25,-.25){2};
	\node at (.5,0){3};
	\node at (.75,.25){4};
	\end{tikzpicture}};
\node at (8,0){$w_6=$};
\node at (9,0){\begin{tikzpicture}[baseline=-0.85ex]
	\node at (.25,-.25){2};
	\node at (.5,0){3};
	\node at (.75,.25){4};
	\node at (1,.5){8};
	\node at (1.25,.25){7};
	\end{tikzpicture}};
\end{tikzpicture}
\]

Crossing \scalebox{1.5}{$\circled{2}$}: $M_2$ crosses $M_1$ in the module \begin{tikzpicture}[baseline=-0.85ex]
	\node at (.25,-.25){\bf{3}};
	\node at (.5,0){\bf{4}};
	\end{tikzpicture}  and 

\[
\begin{tikzpicture}
\node at (-1.35,0){$w_3=$};
\node at (0,0){\begin{tikzpicture}[baseline=-0.85ex]
	\node at (.25,-.25){6};
	\node at (.5,-.5){\bf 3};
	\node at (.75,-.25){\bf 4};
	\node at (1,-.5){5};
	\node at (1.25,-.75){6};
	\node at (1.5,-.5){2};
	\end{tikzpicture}};
\node at (2,0){$w_4=$};
\node at (3.4,0){\begin{tikzpicture}[baseline=-0.85ex]
	\node at (0,0){1};
	\node at (.25,-.25){2};
	\node at (.5,0){\bf 3};
	\node at (.75,.25){\bf 4};
	\node at (1,.5){8};
	\node at (1.25,.25){7};
	\end{tikzpicture}};
\node at (5.5,0){$w_5=$};
\node at (6.5,0){\begin{tikzpicture}[baseline=-0.85ex]
	\node at (0,0){6};
	\node at (.25,.25){2};
	\node at (.5,.5){1};
	\end{tikzpicture}};
\node at (8,0){$w_6=$};
\node at (9.2,0){\begin{tikzpicture}[baseline=-0.85ex]
	\node at (.25,-.25){7};
	\node at (.5,0){8};
	\node at (.75,.25){5};
	\node at (1,0){6};
	\node at (1.25,.25){2};
	\end{tikzpicture}};
\end{tikzpicture}
\]

Crossing \scalebox{1.5}{$\circled{3}$}: $M_1$ crosses $M(w_2^{-1}) \simeq M_2$ in the 3-cycle \begin{tikzpicture}[baseline=-0.85ex]
	\node at (0,0){1};
	\node[rotate=360, scale=1.2] at (.5,0){$\stackrel{\alpha}{\to}$};
	\node at (1,0){2};
	\node[rotate=135, scale=1.2] at (.35,-.3){$\to$};
	\node[rotate=-135, scale=1.2] at (.7,-.3){$\to$};
	\node at (.5,-.65){7};
	\end{tikzpicture}  and 

\[
\begin{tikzpicture}
\node at (-1.35,0){$w_3=$};
\node at (0,0){\begin{tikzpicture}[baseline=-0.85ex]
	\node at (0,0){1};
	\node at (.25,.25){7};
	\node at (.5,.5){8};
	\node at (.75,.25){4};
	\node at (1,0){3};
	\node at (1.25,.25){6};
	\end{tikzpicture}};
\node at (2,0){$w_4=$};
\node at (3.4,0){\begin{tikzpicture}[baseline=-0.85ex]
	\node at (.5,-.5){3};
	\node at (.75,-.25){4};
	\node at (1,-.5){5};
	\node at (1.25,-.75){6};
	\node at (1.5,-.5){2};
	\end{tikzpicture}};
\node at (5.5,0){$w_5=$};
\node at (6.3,0){0};
\node at (8,0){$w_6=$};
\node at (10,0){\begin{tikzpicture}[baseline=-0.85ex]
	\node at (0,0){6};
	\node at (.25,-.25){3};
	\node at (.5,0){4};
	\node at (.75,.25){8};
	\node at (1,0){7};
	\node at (1.25,.25){2};
	\node at (1.5,.5){3};
	\node at (1.75,.75){4};
	\node at (2,.5){5};
	\node at (2.25,.25){6};
	\node at (2.5,.5){2};
	\end{tikzpicture}};
\end{tikzpicture}
\]

Crossing \scalebox{1.5}{$\circled{4}$}: $M_1$ crosses  $M(w_2^{-1}) \simeq M_2$ in the arrow $2 \stackrel{\alpha}{\to}  7$ and 

\[
\begin{tikzpicture}
\node at (-1.35,0){$w_3=$};
\node at (.58,0){\begin{tikzpicture}[baseline=-0.85ex]
	\node at (0,0){1};
	\node at (.25,-.25){2};
	\node at (.5,0){3};
	\node at (.75,.25){4};
	\node at (1,0){5};
	\node at (1.25,-.25){6};
	\node at (1.5,0){2};
	\node at (1.75,-.25){7};
	\node at (2,0){8};
	\node at (2.25,-.25){4};
	\node at (2.5,-.5){3};
	\node at (2.75,-.25){6};
	\end{tikzpicture}};
\node at (3,0){$w_4=$};
\node at (3.8,0){0};
\node at (5.5,0){$w_5=$};
\node at (6.8,0){\begin{tikzpicture}[baseline=-0.85ex]
	\node at (0,0){1};
	\node at (.25,-.25){2};
	\node at (.5,0){3};
	\node at (.75,.25){4};
	\node at (1,0){5};
	\end{tikzpicture}};
\node at (8.8,0){$w_6=$};
\node at (10,0){\begin{tikzpicture}[baseline=-0.85ex]
	\node at (0,0){4};
	\node at (.25,-.25){3};
	\node at (.5,0){6};
	\end{tikzpicture}};
\end{tikzpicture}
\]

Crossing \scalebox{1.5}{$\circled{5}$}: $M_1$ crosses  itself, that is $M_1$ crosses $ M(w_1^{-1}) \simeq M_1$,  in the simple module $\mathbf{2}$ and 

\[
\begin{tikzpicture}
\node at (-1.35,0){$w_3=$};
\node at (.58,0){\begin{tikzpicture}[baseline=-0.85ex]
	\node at (0,0){1};
	\node at (.25,-.25){\bf 2};
	\node at (.5,-.5){6};
	\node at (.75,-.25){5};
	\node at (1,0){4};
	\node at (1.25,-.25){3};
	\node at (1.5,-.5){2};
	\node at (1.75,-.25){1};
		\end{tikzpicture}};
\node at (2.5,0){$w_4=$};
\node at (3.9,0){\begin{tikzpicture}[baseline=-0.85ex]
	\node at (.5,-.25){\bf 2};
	\node at (.75,0){3};
	\node at (1,.25){4};
	\node at (1.25,0){5};
	\node at (1.5,-.25){6};
	\node at (1.75,0){2};
	\end{tikzpicture}};
\node at (5.5,0){$w_5=$};
\node at (6.2,0){\begin{tikzpicture}[baseline=-0.85ex]
	\node at (0,0){0};
	\end{tikzpicture}};
\node at (7.3,0){$w_6=$};
\node at (9.3,0){\begin{tikzpicture}[baseline=-0.85ex]
	\node at (0,.25){2};
	\node at (.25,0){6};
	\node at (.5,.25){5};
	\node at (.75,.5){4};
	\node at (1,.25){3};
	\node at (1.25,.5){6};
	\node at (1.5,.75){5};
	\node at (1.75,1){4};
	\node at (2,.75){3};
	\node at (2.25,.5){2};
	\node at (2.5,.75){1};
	\end{tikzpicture}};
\end{tikzpicture}
\]

Smoothing crossings \scalebox{1.5}{$\circled{1}$}, \scalebox{1.5}{$\circled{3}$}   and \scalebox{1.5}{$\circled{4}$}
gives rise to two triangles in the cluster category

$$\gamma_2 \longrightarrow \gamma_3 \oplus \gamma_4 \longrightarrow \gamma_1 \longrightarrow \gamma_2[1],$$ 
$$  \gamma_1 \longrightarrow \gamma_5 \oplus \gamma_6 \longrightarrow \gamma_2 \longrightarrow \gamma_1[1],$$

where in the case of crossing \scalebox{1.5}{$\circled{3}$} the 
 indecomposable object in the cluster category corresponding to $\gamma_5$ is the zero object and 
where in the case of crossing \scalebox{1.5}{$\circled{4}$} the arc $\gamma_4$ corresponds to the arc labelled `1' in the triangulation.

Smoothing crossing \scalebox{1.5}{$\circled{2}$} gives rise to two triangles in the cluster category

$$\gamma_1 \longrightarrow \gamma_3 \oplus \gamma_4 \longrightarrow \gamma_2 \longrightarrow \gamma_1[1],$$ 
$$  \gamma_2 \longrightarrow  \gamma_5 \oplus \gamma_6 \longrightarrow \gamma_1 \longrightarrow \gamma_2[1].$$

Smoothing crossing \scalebox{1.5}{$\circled{5}$} gives rise to two triangles in the cluster category

$$\gamma_1 \longrightarrow \gamma_3 \oplus \gamma_4 \longrightarrow \gamma_1 \longrightarrow \gamma_1[1],$$ 
$$  \gamma_1 \longrightarrow  \gamma_6 \longrightarrow \gamma_1 \longrightarrow \gamma_1[1].$$

In the Jacobian algebra, smoothing crossings \scalebox{1.5}{$\circled{1}$}  and \scalebox{1.5}{$\circled{4}$}  gives rise to the following short exact sequences

$$0 \longrightarrow \begin{tikzpicture}[baseline=-0.85ex]
\node at (0,0){\bf 6};
\node at (.25,-.25){3};
\node at (.5,0){4};
\node at (.75,.25){8};
\node at (1,0){7};
\end{tikzpicture} \longrightarrow \begin{tikzpicture}[baseline=-0.85ex]
	\node at (0,0){1};
	\node at (.25,-.25){2};
	\node at (.5,0){3};
	\node at (.75,.25){4};
	\node at (1,0){5};
	\node at (1.25,-.25){\bf 6};
	\node at (1.5,-.5){3};
	\node at (1.75,-.25){4};
	\node at (2,0){8};
	\node at (2.25,-.25){7};
	\end{tikzpicture} \oplus \begin{tikzpicture}[baseline=-0.85ex]
	\node at (.25,-.25){\bf 6};
	\node at (.5,0){2};
	\end{tikzpicture} \longrightarrow \begin{tikzpicture}[baseline=-0.85ex]
\node at (0,0){1};
\node at (.25,-.25){2};
\node at (.5,0){3};
\node at (.75,.25){4};
\node at (1,0){5};
\node at (1.25,-.25){\bf 6};
\node at (1.5,0){2};
\end{tikzpicture} \longrightarrow 0,$$ 
 
$$0 \longrightarrow \begin{tikzpicture}[baseline=-0.85ex]
\node at (0,0){7};
\node at (.25,.25){8};
\node at (.5,0){4};
\node at (.75,-.25){3};
\node at (1,0){6};
\end{tikzpicture} \longrightarrow 
	\begin{tikzpicture}[baseline=-0.85ex]
	\node at (0,0){1};
	\node at (.25,-.25){2};
	\node at (.5,0){3};
	\node at (.75,.25){4};
	\node at (1,0){5};
	\node at (1.25,-.25){6};
	\node at (1.5,0){2};
	\node at (1.75,-.25){7};
	\node at (2,0){8};
	\node at (2.25,-.25){4};
	\node at (2.5,-.5){3};
	\node at (2.75,-.25){6};
	\end{tikzpicture}
	 \longrightarrow \begin{tikzpicture}[baseline=-0.85ex]
\node at (0,0){1};
\node at (.25,-.25){2};
\node at (.5,0){3};
\node at (.75,.25){4};
\node at (1,0){5};
\node at (1.25,-.25){ 6};
\node at (1.5,0){2};
\end{tikzpicture} \longrightarrow 0,$$ 
 
respectively. However, crossings \scalebox{1.5}{$\circled{2}$} and \scalebox{1.5}{$\circled{3}$} do not give short exact sequences from $M_2$ to $M_1$. 
 Therefore,  as stated in Corollary~\ref{DimensionFormulaModules}, $\dim \Ext^1_J (M_1, M_2) = 2$.
 
 Crossing \scalebox{1.5}{$\circled{2}$} is the only crossing that gives a short exact sequence from $M_1$ to $M_2$
 
$$0 \longrightarrow \begin{tikzpicture}[baseline=-0.85ex]
	\node at (0,0){1};
	\node at (.25,-.25){2};
	\node at (.5,0){\bf 3};
	\node at (.75,.25){\bf 4};
	\node at (1,0){5};
	\node at (1.25,-.25){6};
	\node at (1.5,0){2};
	\end{tikzpicture}
	 \longrightarrow 
	 \begin{tikzpicture}[baseline=-0.85ex]
	\node at (0,0){1};
	\node at (.25,-.25){2};
	\node at (.5,0){\bf 3};
	\node at (.75,.25){\bf 4};
	\node at (1,.5){8};
	\node at (1.25,.25){7};
	\end{tikzpicture} 
	\oplus 
	\begin{tikzpicture}[baseline=-3.85ex]
	\node at (.25,-.25){6};
	\node at (.5,-.5){\bf 3};
	\node at (.75,-.25){\bf 4};
	\node at (1,-.5){5};
	\node at (1.25,-.75){6};
	\node at (1.5,-.5){2};
	\end{tikzpicture}
	 \longrightarrow 
	\begin{tikzpicture}[baseline=-0.85ex]
\node at (0,0){ 6};
\node at (.25,-.25){\bf 3};
\node at (.5,0){\bf 4};
\node at (.75,.25){8};
\node at (1,0){7};
\end{tikzpicture}
	 \longrightarrow 0.$$ 
 
 Therefore $\dim \Ext^1_J (M_2, M_1) = 1$.
   And finally, we have $$ \dim \Ext^1_J (M_1, M_2) + \dim \Ext^1_J (M_2, M_1) = 3 = \Int(\gamma_1, \gamma_2) - k$$ where $\Int(\gamma_1, \gamma_2) =4$ and 
  $k=1$ corresponds to the only
  3-cycle crossing corresponding to crossing \scalebox{1.5}{$\circled{3}$}.

Crossing \scalebox{1.5}{$\circled{5}$}  gives a short exact sequence from $M_1$ to $M_1$

$$0 \longrightarrow \begin{tikzpicture}[baseline=-0.85ex]
	\node at (0,0){\bf 2};
	\node at (.25,-.25){6};
	\node at (.5,0){5};
	\node at (.75,.25){4};
	\node at (1,0){3};
	\node at (1.25,-.25){2};
	\node at (1.5,0){1};
	\end{tikzpicture}
	 \longrightarrow 
	 \begin{tikzpicture}[baseline=-1.85ex]
	\node at (0,0){1};
	\node at (.25,-.25){\bf 2};
	\node at (.5,-.5){6};
	\node at (.75,-.25){5};
	\node at (1,0){4};
	\node at (1.25,-.25){3};
	\node at (1.5,-.5){2};
	\node at (1.75,-.25){1};
	\end{tikzpicture} 
	\oplus 
	\begin{tikzpicture}[baseline=-0.85ex]
	\node at (.5,-.25){\bf 2};
	\node at (.75,0){3};
	\node at (1,.25){4};
	\node at (1.25,0){5};
	\node at (1.5,-.25){6};
	\node at (1.75,0){2};
	\end{tikzpicture}
	 \longrightarrow 
	\begin{tikzpicture}[baseline=-0.85ex]
	\node at (0,0){1};
	\node at (.25,-.25){\bf 2};
	\node at (.5,0){3};
	\node at (.75,.25){4};
	\node at (1,0){5};
	\node at (1.25,-.25){6};
	\node at (1.5,0){2};
	\end{tikzpicture}
	 \longrightarrow 0.$$ 
 
Since crossing  \scalebox{1.5}{$\circled{5}$}   is the only self-crossing of $\gamma_1$, we have $\dim \Ext^1_J (M_1, M_1) = 1$. 
 
Since $\gamma_2$ has no self-crossings, 
 we have $\dim \Ext^1_J (M_2, M_2) = 0$. 

\appendix

\section{On the canonicity of the generalized cluster category associated with a surface}

Claire Amiot

Let $k$ be an algebraically closed field and $(S,\mathcal M)$ be a marked surface (such that all marked points are in the boundary of $S$).  The cluster category $\mathcal{C}(S,\mathcal M)$ associated to the marked surface $(S,\mathcal M)$ is defined to be the generalized cluster category $\mathcal{C}_T:=\mathcal{C}(Q(T),W(T))$ (as defined in \cite{A}) where $T$ is a triangulation of $(S,\mathcal M)$ and $(Q(T),W(T))$ is the quiver with potential associated to $T$ by \cite{L1}.
 This category is well defined in the following sense: if $T'$ is another triangulation of $(S,\mathcal M)$, then combining the main results of \cite{L1} and \cite{KY} one gets an equivalence of triangulated categories $\mathcal{C}_T\simeq \mathcal{C}_{T'}$. This implies that $\mathcal{C}(S,\mathcal{M})$  is only well-defined up to equivalence of categories. Indeed, a priori given $T$ and $T'$ there is no canonical equivalence  $\mathcal{C}_T\simeq \mathcal{C}_{T'}$. The aim of this appendix is to exhibit some questions and problems this non-canonicity may pose. 
 
 More precisely let us recall the following result due to Br\"ustle and Zhang.
 
\begin{Theorem} \cite[Thm 1.1]{BZ} Let $(S,\mathcal M)$ be a marked surface such that all marked points are on the boundary of $S$. A parametrization of the isoclasses of indecomposable objects in $\mathcal{C}(S,\mathcal M)$ is given by string objects and band objects, where
 \begin{itemize}
 \item[(1)] the string objects are indexed by the homotopy classes of non contractible curves in $(S,\mathcal M)$ with end points in $\mathcal M$ which are not homotopic to a boundary segment of $(S,\mathcal M)$, subject to the equivalence relation $\gamma\sim\gamma^{-1}$;
 \item[(2)] the band objects are indexed by $k^*\times \pi_1^*(S)/\sim$ where $\pi_1^*(S)/\sim$ is given by the nonzero elements of the fundamental group of $S$ subject to the equivalence relation generated by cyclic permutation and $\gamma\sim\gamma^{-1}$.
\end{itemize}  
\end{Theorem}
 
Regarding this result one could first think that we get a description of the objects of $\mathcal{C}(S,\mathcal M)$ independent of the choice of a triangulation. However the parametrization depends on the choice of a triangulation. So given a triangulation $T$ let us denote by $s^T$ (resp. $b^T$) the above bijections that send a curve (resp. a curve with a scalar) to an indecomposable object in $\mathcal{C}_T$. The different facts shown in this appendix can be summarized as the following:
 
\begin{Proposition}\label{main prop} Let $(S,\mathcal M)$ be a marked surface such that all marked points are on the boundary of $S$. 
\begin{enumerate}
\item 
There exist triangulations $T$ and $T'$, and an equivalence $\Phi:\mathcal{C}_T\to \mathcal{C}_{T'}$ of triangulated categories that sends a string object to a band object.

\item For any triangulations $T$ and $T'$, there exists an equivalence $\Phi:\mathcal{C}_T\to \mathcal{C}_{T'}$ of triangulated categories such that $\Phi\circ s^T=s^{T'}$.

\item For the equivalences $\Phi$ of $(2)$, we may have $\Phi\circ b^T\neq b^{T'}$.
\end{enumerate}
\end{Proposition} 

As a consequence of $(2)$ we obtain a bijection between string objects in $\mathcal{C}(S,\mathcal{M})$  and homotopy classes of curves in $(S,\mathcal M)$ that does not depend on the choice of a triangulation.

\begin{Remark} \label{remark} A generalization of the bijection $s^T$ has been constructed in \cite{QZ} in the case where $(S,\mathcal M)$ is a marked surface with punctures and $T$ is an admissible triangulation (see \cite[Def 5]{QZ}). In the same paper, the authors state a more general analogue of (2) in the case where $T$ and $T'$ are admissible \cite[Thm 3.10]{QZ}. The proof given there is not fully detailed.  First, the fact that there exists a canonical equivalence between $\mathcal{C}_T$ and $\mathcal{C}_{T'}$ \cite[(3.3) in subsection 3.1]{QZ} is used without being proved. Indeed the equivalence constructed in \cite{KY} is not canonical since it depends of a choice of a right equivalence (see next subsection for more details). Secondly the proof in \cite[Appendix C]{QZ} does not stress the importance of signs in the computation of the mutation of decorated representations in the sense of \cite{DWZ}. Point (3)  of Proposition~\ref{main prop} above and subsection \ref{subsection exemple} below show that the manipulation of signs is actually a subtle issue in the computation. Though the results presented here are not strictly speaking original, I have thought it would be useful to the community to clarify the aforementionned issues.
\end{Remark}

\subsection{The problem on reduction}

\subsubsection{Construction of the equivalence $\mathcal{C}_T\to\mathcal{C}_{T'}$}

Let $T$ and $T'$ be triangulations of $(S,\mathcal M)$, and $\bf s$ be a sequence of flips such that $T'=\mathfrak{f}_{\bf s}(T)$. Then Labardini constructed in \cite{L1} a right equivalence between the associated quivers with potentials:
\[ \varphi_{\bf s}: (Q(T'),W(T'))\to \mu_{\bf s}(Q(T),W(T)),\] where $\mu_{\bf s}$ is the mutation of the quiver with potential defined by Derksen, Weyman and Zelevinsky in \cite{DWZ}.
Then by \cite{KY}, there exists an equivalence $\Phi_{\bf s}:\mathcal{C}_T\to \mathcal{C}_{T'}$. 
This equivalence depends  not only on the choice of the sequence of flips $\bf s$ but it also depends on the choice of a right equivalence at each flip. 
 
Let us concentrate on this second dependence and assume that $\bf s=\bf i$, that is $T$ and $T'$ differ only by one flip of an arc $\bf i$. Denote by $(Q,W):=(Q(T),W(T))$ the QP associated with $T$ and by $(\widetilde{Q'},\widetilde{W'}):= \widetilde{\mu_{\bf i}}(Q,W))$ the unreduced mutation of $(Q,W)$ at vertex $\bf i$. 
Now fix the right equivalence $\varphi:(\widetilde{Q'},\widetilde{W'})\to (Q',W')=(Q(T'), W(T'))$ corresponding to a reduction of the quiver with potential $(\widetilde{Q'},\widetilde{W'})$. The functor $\Phi_{\bf i}$ constructed by Keller and Yang is the composition of two equivalences:

\[\xymatrix{\mathcal{C}_T\ar[rr]^-{\widetilde{\Phi_{\bf i}}} && \mathcal{C}_{(\widetilde{Q'},\widetilde{W'})}\ar[rr]^-{R^\varphi} && \mathcal{C}_{T'}},\]
where $R^\varphi$ is induced by $\varphi_*$ on the corresponding Ginzburg DG algebras (see  \cite[Lemma 2.9]{KY}). 

\subsubsection{Link with mutations of (decorated) representations}

The authors in  \cite{DWZ} define a notion of mutation of decorated representations of a (non degenerate) quiver with potential: for a module $M\in \mod \Jac(Q,W)$ and a vertex $\textbf i$ of $Q$, they define a module $\widetilde{\mu_{\textbf i}}(M)\in \mod\Jac(\widetilde{Q'},\widetilde{W'})$. Then any reduction $\varphi:(\widetilde{Q'},\widetilde{W'})\to  (Q',W')$ induces an isomorphism of algebras $\varphi: \Jac(\widetilde{Q'},\widetilde{W'})\to \Jac (Q',W')$ by \cite[Prop 3.7]{DWZ}, and so an equivalence 

\[ \xymatrix{\mod \Jac(\widetilde{Q'},\widetilde{W'}) \ar[rr]^-{\varphi_*}&& \mod\Jac(Q',W')}.\]

Note that in \cite[Def 10.4]{DWZ} the equivalence is defined to be the restriction functor $(\varphi^{-1})^*$ but since $\varphi$ is an isomorphism of algebras we have $(\varphi^{-1})^*=\varphi_*$.

The mutation $\mu_{\textbf i}(M)$ of $M$ at $\textbf i$ is then defined to be $\varphi_*(\widetilde{\mu}_{\textbf i}(M))$. This implies that $\mu_{\textbf i}(M)$ is only defined up to right equivalence of representation (cf Remark 10.3 in \cite{DWZ}) and not up to isomorphism of module.

By construction the cluster category $\mathcal{C}_T$ comes with a canonical cluster-tilting object that we denote by $X$, with the property $\End_{\mathcal{C}_T}(X)\simeq \Jac(Q,W)$.  The functor $F:=\Hom_{\mathcal{C}_T}(X[-1],-):\mathcal{C}_T\to \mod \Jac(Q,W)$ is dense and sends any indecomposable object not isomorphic to a summand of $X$ to an indecomposable module. Similarly we denote by $\widetilde{X'}$ (resp. $X'$) the canonical cluster-tilting objects in $\mathcal{C}_{(\widetilde{Q'},\widetilde{W'})}$ (resp. in $ \mathcal{C}_{T'}$). Then we have the following:

\begin{equation}\label{diagram}\xymatrix{\mathcal{C}_T\ar[rr]^-{\widetilde{\Phi_{\textbf i}}}\ar[d]_{F} && \mathcal{C}_{(\widetilde{Q'},\widetilde{W'})}\ar[rr]^-{R^\varphi}\ar[d]^{\widetilde{F'}} && \mathcal{C}_{T'}\ar[d]^{F'}\\ 
\mod \Jac(Q,W)\ar@{-->}[rr]^-{\widetilde{\mu}_{\textbf i}} &&\mod \Jac(\widetilde{Q'},\widetilde{W'}) \ar[rr]^-{\varphi_*}&& \mod\Jac(Q',W')}.
\end{equation}

In this diagram, the right hand square is clearly commutative, since $R^\varphi(\widetilde{X'})$ is isomorphic to $X'$ in the category $\mathcal{C}_{T'}$.  
Moreover Plamondon showed in \cite[Prop 4.1]{Pla} that if $M$ is an object of $\mathcal{C}_T$, then $\widetilde{\mu}_{\textbf i}(F(M))$ and $\widetilde{F'}(\widetilde{\Phi_{\textbf i}}(M))$ are isomorphic in $\mod \Jac (\widetilde{Q'},\widetilde{W'})$. Hence the left square also commutes.

Since for any reduction $\varphi$, if we compose $\varphi$ with an automorphism of the algebra $\Jac(Q',W')$ that fixes each vertex we obtain another reduction, it is easy to construct an equivalence $\Phi_{\textbf i}:\mathcal{C}_T\to \mathcal{C}_{T'}$ sending certain string objects to band objects (see Remark 10.3 in \cite{DWZ}). For instance if $(S,\mathcal M)$ is the annulus with two marked points then for any triangulation $T$, $(Q_T,W_T)$ is the Kronecker quiver with the zero potential. Denote by $a$ and $b$ its arrows, the automorphism sending $a$ to $a$ and $b$ to $a+b$ sends the string $a$ to a band  of type $ab^{-1}$. The automorphism exchanging $a$ and $b$ sends the string $a$ to the string $b$. 

\[\scalebox{0.7}{\begin{tikzpicture}[scale=0.4, >=stealth]

\draw [black, thick ](-14,-14) circle (5);
\draw [black, thick ](-14,-14) circle (1);
\node at (-19,-14) {$\bullet$};
\node at (-14,-15){$\bullet$};

\draw[red,thick] (-14,-14) circle (3);
\node at (-14,-21) {band $ab^{-1}$};

\draw [black, thick ](14,-14) circle (5);
\draw [black, thick ](14,-14) circle (1);
\node at (9,-14) {$\bullet$};
\node at (14,-15){$\bullet$};

\draw [red,thick]  (14,-15).. controls (12,-16) and (10,-10).. (14,-10);
\draw[red, thick] (14,-10).. controls (18,-10) and (18,-17)..(14,-17);
\draw[red, thick] (14,-17)..controls (11,-17) and (11,-12)..(14,-12);
\draw[red, thick] (14,-12).. controls (16,-12) and (16,-16)..(14,-15);

\draw [black, thick ](0,-14) circle (5);
\draw [black, thick ](0,-14) circle (1);
\node at (-5,-14){$\bullet$};
\node at (0,-15){$\bullet$};

\draw[red, thick] (-5,-14)..controls (-3,-10) and (2,-10).. (2,-14);
\draw[red, thick] (2,-14).. controls (2,-17) and (-2,-17)..(-2,-14);
\draw[red, thick] (-2,-14).. controls (-2,-10) and (3,-10)..(3,-14);
\draw[red, thick] (3,-14).. controls (3,-18) and (-4,-18).. (-5,-14);

\node at (0,-21) { string $a$};
\node at (14,-21) {string $b$};
\end{tikzpicture}}\]

\subsection{Canonical bijection for strings}

In this subsection, we prove that the bijection $s^T$ constructed in \cite{BZ} is independent of $T$ as soon as we only allow certain reductions. 

\subsubsection{Canonical reduction}

One way to handle this problem is to allow only certain kind of reductions when constructing the triangle equivalences $\Phi_{\textbf i}$. This can be done easily in the case where $(S,\mathcal{M})$ is a marked surface with all marked points are on the boundary. Indeed in this case the quiver with potential is especially simple.  

Let $T$ be a triangulation, $\textbf i $ be an arc of $T$ and $T':=\mathfrak f_{\textbf i}(T)$ be the flip of $T$ at $\textbf i $. Denote by $\Delta_1$ and $\Delta_2$ the triangles in $T$ having $\textbf i $ as a side. Then $\Delta_1$ and $\Delta_2$ are distinct. If $\Delta_1$ and $\Delta_2$ are internal triangles then locally the quiver with potential $(Q_T,W_T)$ associated to $T$ looks as follows:

\[\scalebox{0.8}{\begin{tikzpicture}[scale=0.9,>=stealth]
\draw[loosely dotted] (0,0)--node (J1)[fill=white, inner sep=0pt]{${\textbf j}_1$}(8,0);
\draw[loosely dotted] (8,0)--node (K2)[fill=white, inner sep=0pt]{${\textbf k}_2$}(12,4);
\draw[loosely dotted] (12,4)--node (J2)[fill=white, inner sep=0pt]{${\textbf j}_2$}(4,4);
\draw[loosely dotted] (4,4)--node (K1)[fill=white, inner sep=0pt]{${\textbf k}_1$}(0,0);
\draw[loosely dotted] (4,4)--node (I)[fill=white, inner sep=0pt]{${\textbf i}$}(8,0);
\draw[<-] (I)--node [fill=white, inner sep=0pt]{$c_2$}(K2);
\draw[<-] (I)--node [fill=white, inner sep=0pt]{$c_1$}(K1);
\draw[<-] (K1)--node [fill=white, inner sep=0pt]{$a_1$}(J1);
\draw[<-] (K2)--node [fill=white, inner sep=0pt]{$a_2$}(J2);
\draw[<-] (J1)--node [fill=white, inner sep=0pt]{$b_1$}(I);
\draw[<-] (J2)--node [fill=white, inner sep=0pt]{$b_2$}(I);
\end{tikzpicture}}\]
with $W=c_1b_1a_1+c_2b_2a_2$. Note that we may have $\textbf j_1=\textbf j_2$ (or $\textbf k_1=\textbf k_2$), but in that case, there is a way to distinguish the arrow $b_1$ from the arrow $b_2$ since each arrow is canonically associated to a triangle of the triangulation. 
By \cite{DWZ} we obtain the following quiver with potential $\widetilde{\mu_i}(Q,W)=( \widetilde{Q'},\widetilde{W'})$ after `unreduced' mutation at $\textbf i $. 

\[\scalebox{0.8}{\begin{tikzpicture}[scale=0.9,>=stealth]

\node at (0,2) {$\widetilde{Q'}=$};
\draw[loosely dotted] (0,0)--node (J1)[fill=white, inner sep=0pt]{${\textbf j}_1$}(8,0);
\draw[loosely dotted] (8,0)--node (K2)[fill=white, inner sep=0pt]{${\textbf k}_2$}(12,4);
\draw[loosely dotted] (12,4)--node (J2)[fill=white, inner sep=0pt]{${\textbf j}_2$}(4,4);
\draw[loosely dotted] (4,4)--node (K1)[fill=white, inner sep=0pt]{${\textbf k}_1$}(0,0);
\draw[loosely dotted] (0,0)--node (I)[fill=white, inner sep=0pt]{${\textbf i}$}(12,4);

\draw[->] (I)--node [fill=white, inner sep=0pt]{$c_2^*$}(K2);
\draw[->] (I)--node [fill=white, inner sep=0pt]{$c_1^*$}(K1);
\draw[<-] (2.3,1.8)--node [fill=white, inner sep=0pt]{$a_1$}(4,0.4);
\draw[->] (2,1.6)--node [fill=white, inner sep=0pt,yshift=-3mm]{$[c_1b_1]$}(3.7,0.2);
\draw[<-] (10,2.4)--node [fill=white, inner sep=0pt]{$a_2$}(8.3,3.8);
\draw[->] (9.7,2.2)--node [fill=white, inner sep=0pt, yshift=-3mm]{$[c_2b_2]$}(8,3.6);

\draw[->] (J1)--node [fill=white, inner sep=0pt]{$b_1^*$}(I);
\draw[->] (J2)--node [fill=white, inner sep=0pt]{$b_2^*$}(I);

\draw[<-] (J2)--node [fill=white, inner sep=0pt]{$[c_1b_2]$}(K1);
\draw[<-] (J1)--node [fill=white, inner sep=0pt]{$[c_2b_1]$}(K2);
\end{tikzpicture}}\]

$$\widetilde{W'}=[c_1b_1]a_1+[c_2b_2]a_2+[c_1b_1]b_1^*c_1^*+[c_2b_2]b_2^*c_2^*+[c_2b_1]b_1^*c_2^*+[c_1b_2]b_2^*c_1^*.$$

Let $\varphi^T_{\textbf i }:k\widetilde{Q'}\to k\widetilde{Q'}$ be the map defined by 
\[\left \{ \begin{array}{rl}  \varphi_{\textbf i }^T (\alpha) = &\alpha \textrm{ if } \alpha\neq a_1,a_2\\
\varphi_{\textbf i }^T(a_\ell)= &a_\ell-b_\ell^*c_\ell^* \textrm{ for }\ell=1,2. \\
 \end{array}\right.\]
 
It is immediate to see that $\varphi^T_{\textbf i }$ is a right equivalence between the quiver with potential $(\widetilde{Q'},\widetilde{W'})$ and the direct sum of the quiver with potential $(Q',W')$ with a trivial quiver with potential, where $(Q',W')$ is the quiver with potential associated with the triangulation $T'=\mathfrak{f}_{\textbf i }(T)$. Note that if $\Delta_1$ or $\Delta_2$ have boundary sides, then the arcs $\textbf j_1$, $\textbf j_2$, $\textbf k_1$ and $\textbf k_2$ may not exist, so the quiver with potential $(Q_T,W_T)$ is simpler and the definition of the reduction $\varphi_{\textbf i }^T$ is similar. 
Moreover in the case where $(\widetilde{Q'},\widetilde{W'})$ does not have any $2$-cycle, then $\varphi^T_{\textbf i }$ is the identity map. This leads to introduce the following.
 
 \begin{Definition} For any arc $\textbf i $ in $T$,  the map $\varphi_{\textbf i }^T$ is called the \emph{canonical reduction} at $\textbf i $. 
 \end{Definition}
 
\subsubsection{The bijection $s^T$}
 Before proving the main result, let us recall the construction of the map $s^T$, and some properties on string modules that will be used in the proof.
 
By \cite{ABCP} the Jacobian algebra $\Jac(Q,W)$ is a string algebra. In such an algebra, a word $w=\alpha_1\ldots \alpha_n$ of arrows or formal inverse of arrows of $Q^{\rm op}$ is called a \emph{string} if $\alpha_{i+1}\neq \alpha_i^{-1}$ and no subword nor its inverse belongs to the Jacobian ideal. 
 
Let $w=\alpha_n\ldots \alpha_1$ be a string and $\lambda_1,\ldots,\lambda_n$ be in $k ^*$. Define a module $M^T(w:\lambda_1,\ldots,\lambda_n)$ in $\mod \Jac(Q,W)={\rm Rep}(Q^{\rm op},\partial W)$ as follows:  For any $\ell=0\ldots, n$ let $M_\ell$ be a $1$-dimensional $k$-vector space. Then for any vertex $\textbf i$ in $Q_0$ set 

\[ M^T_{\textbf i}  = \left\{ \begin{array}{lr} \bigoplus _{\ell, t(\alpha_\ell)=\textbf i}M_\ell \oplus M_0& \textrm{ if } s(\alpha_1)= \textbf i\\ 
  \bigoplus _{\ell, t(\alpha_\ell)=\textbf i} M_\ell & \textrm{ else }\end{array}\right.\] For any arrow $\alpha:\textbf i\to \textbf j$ in $Q^{\rm op}_1$, if there exists $\ell$ such that $\alpha=\alpha_\ell$  (resp. $\alpha^{-1}=\alpha_\ell$ ) then $M_{\ell-1}$ is a direct summand of  $M_{\textbf i}^T$ (resp. $M_{\textbf j}^T$) and $M_\ell$ a summand of $M^T_{\textbf j}$ (resp. $M^T_{\textbf i}$) and the restriction of $M^T_\alpha$ to $M_{\ell -1}$ (resp. $M_\ell$) is the multiplication by $\lambda_\ell$ from $M_{\ell -1}$ to $M_\ell$ (resp. the multiplication by $\lambda_{\ell}^{-1}$ from $M_\ell$ to $M_{\ell -1}$) .  

\begin{Definition}
The \emph{string module} associated to $w$ is defined to be $M^T(w):=M^T(w:1,\ldots,1)$.
\end{Definition}

The following isomorphisms are easy to check and classical \cite{BR}.

\begin{Lemma} We have the following isomorphisms in $\mod \Jac(Q,W)$:

\[ \begin{array}{rll}M^T(w:\lambda_1,\ldots,\lambda_n) &=M^T(w^{-1}:\lambda_n^{-1},\ldots,\lambda_1^{-1}) \\ &\simeq M^T(w:1,\lambda_2,\ldots,\lambda_n) &\\& \simeq M^T(w:\lambda_1,\ldots,\lambda_{\ell-1}\lambda_{\ell},1,\lambda_{\ell+1},\ldots,\lambda_n)& \forall \ell\end{array}\]
\end{Lemma}

\begin{Corollary}\label{lemmastrings}
The module $M^T(w:\lambda_1,\ldots,\lambda_n)$ is isomorphic to the string module $M^T(w)=M^T(w^{-1})$.
\end{Corollary}

\medskip

Now let $\gamma$ be an oriented curve on $S$ with endpoints in $\mathcal M$ which is not homotopic to a boundary segment or to an arc of $T$. Up to homotopy, we may assume that $\gamma$ intersects each arc of $T$ transversally and does not cut an arc of $T$ twice in succession. Then one can associate to $\gamma$ a sequence $w$ of arrows or inverse arrows in $Q_1^{\rm op}$ corresponding to the angles of $T$ intersected by $\gamma$. The map $\gamma\to w(\gamma)$ is shown to be a bijection between non trivial homotopy classes of such oriented curves in $(S,\mathcal M)$ and strings in $\Jac(Q,W)$ in \cite{ABCP}. Moreover, the string associated to $\gamma^{-1}$ is clearly $w(\gamma)^{-1}$.

The bijection $s^T$ is defined in \cite{BZ} as follows: if $\gamma=\textbf i$ is an arc of $T$, then $s^T(\gamma)$ is defined to be the object $X_{\textbf i}$ which is the indecomposable summand of the canonical cluster-tilting object $X$ in $\mathcal{C}_T$ corresponding to the vertex $\textbf i$ of $Q$; if $\gamma$ is not an arc of $T$, then $s^T(\gamma)$ is the indecomposable object $X(\gamma)$ such that $F(X(\gamma))\simeq M^T(w(\gamma))$ in $\mod \Jac(Q,W)$.  This indecomposable object is unique up to isomorphism.
 
\subsubsection{Compatibility for strings} 
 
The following result is the main result of this appendix.
 
\begin{Theorem}
Let $T$ and $T'$ be triangulations of a marked surface with marked points on the boundary, and $s^T$ and $s^{T'}$ be the bijections described above. For a sequence of flips $\textbf s$ such that $\mathfrak{f}_{\textbf s}(T)=T'$, denote by $\Phi_{\textbf s}:\mathcal{C}_T\to \mathcal{C}_{T'}$ the equivalence defined in \cite{KY} where at each mutation we apply the canonical reduction. Then for any such sequence $\textbf s$ we have $$\Phi_{\textbf s}\circ s^T=s^{T'}.$$
\end{Theorem}
 
\begin{proof}
 Let $T$ be a triangulation of $(S,\mathcal M)$ and $\textbf i$ be an arc of $T$. Denote by $T'$ the triangulation $\mathfrak f _{\textbf i}(T)$. Denote by $\varphi$ the canonical reduction $(\widetilde{Q'},\widetilde{W'})\to (Q',W')$ defined above. It is enough to show that $R^\varphi\circ \widetilde{\Phi_{\textbf i}}\circ s^T=s^{T'}$. 
 
Let $\gamma$ be a curve on $S$ with endpoints in $\mathcal{M}$. We consider the following cases.

\medskip
\noindent
\emph{Case 1: $\gamma=\textbf j$ is an arc of $T$ and $T'$.} 
Then we have $s^T(\textbf j)=X_{\textbf j}$ and  from the definition of $\widetilde{\Phi_{\textbf i}}$ we have $\widetilde{\Phi_{\textbf i}}(X_{\textbf j})\simeq \widetilde{X'}_{\textbf j}$ and $R^\psi\circ\widetilde{\Phi_{\textbf i}}(X_{\textbf j})\simeq X'_{\textbf j}=s^{T'}(\textbf j)$ for any choice of reduction $\psi$. 

\medskip
\noindent
\emph{Case 2: $\gamma=\textbf i$.}
Then $s^T(\textbf i)=X_{\textbf i}$. Denote by ${\textbf i}'$ the arc of $T'$ which is not in $T$. By definition, for any choice of reduction $\psi$, the object $R^\psi\circ \widetilde{\Phi_{\textbf i}}(X_{\textbf i})$ is the cone of the map 
$$X'_{\textbf i'}\to \bigoplus_{{\textbf i}'\to k \in Q'_1} X'_k$$ which is isomorphic to the cocone of the map
$$ \bigoplus_{ j\to {\textbf i'} \in Q'_1} X'_j\to X'_{\textbf i'}$$ by the properties of exchange triangles in a $2$-Calabi-Yau category with cluster-tilting objects. Applying the functor $F'=\Hom_{\mathcal{C}_{T'}}(X'[-1],-)$ we obtain that $F'(\Phi_{\textbf i}(X_{\textbf i}))$ is isomorphic to the cokernel of the map 
$$ \bigoplus_{ j\to {\textbf i'} \in Q'_1} P'_j\to P'_{\textbf i'}$$ where $P'_\ell$ is the projective associated to the vertex $\ell$ in $\Jac(Q',W')$. Therefore $F'(\Phi_{\textbf i}(X_{\textbf i}))$ is isomorphic to the simple module $S_{\textbf i'}$ of $\Jac(Q',W')$ associated to the vertex $\textbf i'$, which is the module $M^{T'}(\textbf i)$.

\medskip
\noindent
\emph{Case 3: $\gamma=\textbf i'$. } This case is similar to the previous case.

\medskip
\noindent
\emph{Case 4: $\gamma$ is not an arc of $T$ and not an arc of $T'$.} Denote by $w$ (resp. $w'$) the string in $\Jac(Q,W)$ (resp. $\Jac (Q',W)$) corresponding to $\gamma$. Then by the commutative diagram (\ref{diagram}) and the definition of the bijections $s^T$ and $s^{T'}$, it is enough to show that $\varphi_*(\widetilde{\mu_{\textbf i}}M^T(w))$ is isomorphic to $M^{T'}(w')$ in $\mod \Jac(Q',W')$.

As in the previous section, we denote by $\Delta_1$ and $\Delta_2$ the triangles of $T$ sharing $\textbf i$. There is exactly 14 ways for $\gamma$ to cross $\Delta_1\cup\Delta_2$, which are described in the following pitcures:

\[\scalebox{0.8}{\begin{tikzpicture}[scale=1,>=stealth]
\draw (0,0)-- (4,0)--(6,2)--(2,2)--(0,0);
\draw (2,2)--(4,0);

\draw[blue] (0,0.5)..controls (0.5,0.5) and (0.8,0)..(0.75,-0.5);

\draw[blue] (6,1.5)..controls (5.5,1.5) and (5.2,2)..(5.25,2.5);

\draw[blue] (1,1.5)..controls (1.5,1.5) and (3,2)..(3,2.5);

\draw[blue] (3,-0.5)..controls (3,0) and (4,0.5)..(5,0.5);

\draw[blue] (0,1)--(6,1);
\draw[blue] (2,-1)--(4,3);
\node[fill=white, inner sep=0pt] at (0.5,0.25) {$1$};
\node[fill=white, inner sep=0pt] at (4,0.3) {$2$};
\node[fill=white, inner sep=0pt] at (5.5,1.75) {$3$};
\node[fill=white, inner sep=0pt] at (2,1.7) {$4$};
\node[fill=white, inner sep=0pt] at (0.2,1) {$5$};
\node[fill=white, inner sep=0pt] at (2.25,-0.5) {$6$};
\end{tikzpicture}}\]

\[\scalebox{0.8}{\begin{tikzpicture}[scale=1,>=stealth]
\draw (0,0)-- (4,0)--(6,2)--(2,2)--(0,0);
\draw (2,2)--(4,0);

\draw[blue] (0,0)--(3.5,2.5);
\draw[blue] (0,0)--(6,1);
\draw[blue] (0,1)--(6,2);
\draw[blue] (2.5,-0.5)--(6,2);

\node[fill=white, inner sep=0pt] at (3.5,2.5) {$9$};
\node[fill=white, inner sep=0pt] at (6,1) {$8$};
\node[fill=white, inner sep=0pt] at (0,1) {$10$};
\node[fill=white, inner sep=0pt] at (2.5,-0.5) {$7$};

\draw (10,2)--(12,0)--(8,0)--(10,2)--(14,2)--(12,0);

\draw[blue] (10,-1)--(10,2)--(14.5,0.5);
\draw[blue] (12,3)--(12,0)--(7.5,1.5);

\node[fill=white, inner sep=0pt] at (10,-1) {$11$};
\node[fill=white, inner sep=0pt] at (14.5,0.5) {$12$};
\node[fill=white, inner sep=0pt] at (12,3) {$13$};
\node[fill=white, inner sep=0pt] at (7.5,1.5) {$14$};
\end{tikzpicture}}\]

Denote by $n_\ell$, $\ell=1,\ldots,14$ the number of times $\gamma$ intersects $\Delta_1\cup\Delta_2$ in the way $\ell$ (in both directions), and define $V_\ell=k^{n_\ell}$. Note that $\sum_{\ell=7}^{14}n_\ell\leq 2$ since these crossings correspond to endpoints of $\gamma$. 

Then the restriction of the representation $M^T(w)$ to the quiver $Q_{\Delta_1\cup\Delta_2}$ is the following:  

\[\scalebox{0.7}{\begin{tikzpicture}[scale=1.8,>=stealth]
\node at (0,2) {$M^T(w)=$};

\draw[loosely dotted] (0,0)--node (J1)[fill=white, inner sep=0pt]{$V_1\oplus V_2\oplus V_6$}(8,0);
\node at (4,-0.25) {$\oplus V_{7}\oplus V_{11}$};
\draw[loosely dotted] (8,0)--node (K2)[fill=white, inner sep=0pt]{$V_2\oplus V_3\oplus V_5$}(12,4);
\node at (10,1.75) {$\oplus V_8\oplus V_{12}$};

\draw[loosely dotted] (12,4)--node (J2)[fill=white, inner sep=0pt]{$\oplus V_{9}\oplus V_{13}$}(4,4);
\node at (8,4.25) {$V_3\oplus V_4\oplus V_6$};

\draw[loosely dotted] (4,4)--node (K1)[fill=white, inner sep=0pt]{$\oplus V_{10}\oplus V_{14}$}(0,0);
\node[fill=white, inner sep=0pt] at (2,2.25){$V_1\oplus V_4\oplus V_5$};

\draw[loosely dotted] (4,4)--(8,0);
\node (I2)[fill=white, inner sep=0pt] at (6,2.15) {$V_2\oplus V_4\oplus V_5\oplus V_6$};
\node (I1)[fill=white, inner sep=0pt] at (6,1.85) {$V_7\oplus V_8\oplus V_9\oplus V_{10}$};
 
\draw[->] (I1)--node [fill=white, inner sep=0pt]{$1_{V_2}\oplus 1_{V_5}\oplus 1_{V_8}$}(K2);
\draw[->] (I2)--node [fill=white, inner sep=0pt]{$1_{V_4}\oplus 1_{V_5}\oplus 1_{V_{10}}$}(K1);
\draw[->] (K1)--node [fill=white, inner sep=0pt]{$1_{V_1}$}(J1);
\draw[->] (K2)--node [fill=white, inner sep=0pt]{$1_{V_3}$}(J2);
\draw[->] (J1)--node [fill=white, inner sep=0pt]{$1_{V_2}\oplus 1_{V_6}\oplus 1_{V_7}$}(I1);
\draw[->] (J2)--node [fill=white, inner sep=0pt]{$1_{V_4}\oplus 1_{V_6}\oplus 1_{V_9}$}(I2);
\end{tikzpicture}}\]

Then a direct calculation gives the following representation for $\widetilde{\mu_{\textbf i}}(M^T(w))$:

\[\scalebox{0.7}{\begin{tikzpicture}[scale=1.8,>=stealth]
\draw[loosely dotted] (0,0)--node (J1)[fill=white, inner sep=0pt]{$V_1\oplus V_2\oplus V_6$}(8,0);
\node at (4,-0.25) {$\oplus V_{7}\oplus V_{14}$};
\draw[loosely dotted] (8,0)--node (K2)[fill=white, inner sep=0pt]{$V_2\oplus V_3\oplus V_5$}(12,4);
\node at (10,1.75) {$\oplus V_8\oplus V_{12}$};

\draw[loosely dotted] (12,4)--node (J2)[fill=white, inner sep=0pt]{$\oplus V_{7}\oplus V_{14}$}(4,4);
\node at (8,4.25) {$V_3\oplus V_4\oplus V_6$};

\draw[loosely dotted] (4,4)--node (K1)[fill=white, inner sep=0pt]{$\oplus V_{10}\oplus V_{14}$}(0,0);
\node[fill=white, inner sep=0pt] at (2,2.25){$V_1\oplus V_4\oplus V_5$};

\draw[loosely dotted] (0,0)--(12,4);
\node (I2)[fill=white, inner sep=0pt] at (6,2.15) {$V_1\oplus V_3\oplus V_5\oplus V_6$};
\node (I1)[fill=white, inner sep=0pt] at (6,1.85) {$V_{11}\oplus V_{12}\oplus V_{13}\oplus V_{14}$};
 
\draw[<-] (I1)--node [fill=white, inner sep=0pt]{$1_{V_3}\oplus (-1)_{V_5}\oplus 1_{V_{12}}$}(K2);
\draw[<-] (I2)--node [fill=white, inner sep=0pt]{$(-1)_{V_1}\oplus 1_{V_5}\oplus 1_{V_{14}}$}(K1);
\draw[->] (K1)--node [fill=white, inner sep=0pt]{$1_{V_1}$}(J1);
\draw[->] (K2)--node [fill=white, inner sep=0pt]{$1_{V_3}$}(J2);
\draw[<-] (J1)--node [fill=white, inner sep=0pt]{$1_{V_1}\oplus 1_{V_6}\oplus 1_{V_{14}}$}(I1);
\draw[<-] (J2)--node [fill=white, inner sep=0pt]{$(-1)_{V_3}\oplus (-1)_{V_6}\oplus 1_{V_{13}}$}(I2);

\draw[->] (7.5,4)--node[fill=white, inner sep=0pt]{$1_{V_4}$}(2.5,2.5);
\draw[->] (4.5,0)--node[fill=white, inner sep=0pt]{$1_{V_2}$}(9.5,1.5);

\draw[->] (3.5,0)--node[fill=white, inner sep=0pt]{$0$}(1.75,1.75);
\draw[->] (8.5,4)--node[fill=white, inner sep=0pt]{$0$}(10.25,2.25);
\end{tikzpicture}}\]

Note that in the computation of $\widetilde{\mu}_{\textbf i}(M^T(w))$, the splitting datas $(10.8)$ and $(10.9)$ of \cite{DWZ} are always $0$ or identity, hence the representation does not depend on the choice of these datas. 

Applying the canonical reduction we obtain the same representation except that the action of the arrows $a_1$ and $a_2$ are $0$ instead of $1_{V_3}$ and $1_{V_2}$. This representation is of the form $M^{T'}(w':\lambda_1,\ldots,\lambda_n)$ with $\lambda_\ell=\pm 1$ so using Corollary \ref{lemmastrings} we obtain the isomorphism

 \[ \varphi_*(\widetilde{\mu}_{\textbf i}(M^T(w))\simeq M^{T'}(w')\quad {\rm in }\ \mod \Jac(Q',W')\]
 which ends Case 4 and the proof.
\end{proof}

\begin{Remark}
The same kind of questions can be asked in the case where $(S,\mathcal M)$ is a surface with punctures. As mentionned above in Remark \ref{remark}, if $T$ and $T'$ are admissible triangulations (that are triangulations where every puncture is in a self-folded triangle), then a similar result has been stated in \cite{QZ}. But an analogue of the canonical reduction needs to be defined in the case where $T$ and $T'$ are linked by a $\diamond$-flip. 
 
More generally, if $T$ is any triangulation, $T'=\mathfrak{f}_{\textbf i}(T)$ and $\gamma$ is an arc (thus without selfcrossings) which is not in $T$ and $T'$, then Labardini defined in \cite{Lab2} an indecomposable module $M^T(\gamma)\in \mod \Jac(Q,W)$ and shows that $\mu_{\textbf i}(M^T(\gamma))$ is right equivalent to $M^{T'}(\gamma)$. A priori, this does not entirely prove that there is a bijection compatible with any triangulation between arcs and direct summands of cluster-tilting objects in $\mathcal{C}(S,\mathcal M)$ since $\mu_{\textbf i}(M^T(\gamma))$ is only defined up to right equivalence and not up to isomorphism. In this case, the right equivalence $\mu_{\textbf i}(Q,W)\to (Q',W')$ constructed by Labardini in \cite{L1} is much more complicated to describe and so it is not so clear that an analogue of the `canonical' reduction described in the present work does exist. 
\end{Remark}
 
\subsection{The problem with bands}
 
 The aim of this subsection is to show that the situation is not as nice for bands. Before exhibiting counter-examples, let us redefine the bijection $b^T$ in a more precise way than in  \cite{ABCP} and \cite{BZ}. 
 
\subsubsection{Band modules and the bijection $b^T$}
 In a string algebra, a string $b=\alpha_1\ldots\alpha_n$ is called a \emph{band} if $s(\alpha_1)=t(\alpha_n)$, if any power $b^m$ of $b$ is a string and if $b$ is not a power of any string. 
 
 Let $b=\alpha_1\ldots\alpha_n$ be a band in $\Jac(Q,W)$ and $\lambda_1,\ldots,\lambda_n$ be in $k^*$. We define a module $B^T(b:\lambda_1,\ldots,\lambda_n)$ in $\mod \Jac(Q,W)$ : 
  For any $\ell\in \mathbb Z/n\mathbb Z$ let $B_\ell$ be a $1$-dimensional $k$-vector space. Then for any vertex $\textbf i$ in $Q_0$ set 

\[ (B^T)_{\textbf i}  :=\bigoplus _{\ell, t(\alpha_\ell)=\textbf i}B_\ell .\] For any arrow $\alpha:\textbf i\to \textbf j$ in $Q^{\rm op}_1$, if there exists $\ell$ such that $\alpha=\alpha_\ell$  (resp. $\alpha^{-1}=\alpha_\ell$ ) then $B_{\ell-1}$ is a direct summand of  $(B^T)_{\textbf i}$ (resp. $(B^T)_{\textbf j}$) and $B_\ell$ a summand of $(B^T)_{\textbf j}$ (resp. $(B^T)_{\textbf i}$); the restriction of $(B^T)_\alpha$ to $B_{\ell -1}$ (resp. $B_\ell$) is defined to be the multiplication by $\lambda_\ell$ from $B_{\ell-1}$ to $B_\ell$ (resp. the multiplication by $\lambda_{\ell}^{-1}$ from $B_\ell$ to $B_{\ell -1}$) .  

\begin{Definition}
Let $b$ be a band and $\lambda\in k^*$. 
The (regular simple) \emph{band module} associated with $(b;\lambda)$ is defined to be $B^T(b;\lambda):=B^T(b: \lambda,1,\ldots,1)$.
\end{Definition} 

The following is classical and easy to check \cite{BR}.

\begin{Lemma}
Let $b=\alpha_1\ldots\alpha_n$ be a band and denote by $b':=\alpha_{2}\ldots\alpha_n\alpha_1$. Then for any $\lambda_1,\ldots, \lambda_n$ we have isomorphisms
\[ \begin{array}{rcl} B^T(b:\lambda_1,\ldots,\lambda_n) & \simeq & B^T(b': \lambda_{2},\ldots,\lambda_n,\lambda_1)\\ & \simeq & B^T(b^{-1}:\lambda_n^{-1},\ldots,\lambda_1^{-1})\\ & \simeq & B^T(b:1,\lambda_1\lambda_{2},\ldots,\lambda_n).\end{array}\]
\end{Lemma}  
 
\begin{Corollary}\label{lemmabands}
 The module $B^T(b:\lambda_1,\ldots,\lambda_n)$ is isomorphic to $B^T(b;\prod_{\ell}\lambda_\ell)$ and we have 
 \[B^T(b;\lambda)\simeq B^T(b^{-1};\lambda^{-1}) \simeq B^T(b';\lambda)\] for any $b'$ cyclic permutation of $b$.
\end{Corollary} 
 
\medskip
 
Denote by $\pi_1^{\rm free,*}(S)$ the set of non trivial conjugacy classes of the fundamental group $\pi_1(S)$. This set coincides with the set of non contractible oriented closed curves on $S$ up to free homotopy. Consider the subset $\pi_1^{\rm free, irred,*}(S)\subset  \pi_1^{\rm free}(S)$ of non contractible irreducible closed curves $\gamma$, that are conjugacy classes of closed curves which are not conjugate to a power of a closed curve.  
A natural bijection $b^T$ between $\pi_1^{\rm free, irred,*}(S)$ and the set of bands in $(Q,W)$ up to cyclic permutation is described in  \cite{ABCP}, which sends a curve (transversal to $T$) to the sequence of arrows (and inverse arrows) of $Q^{\rm op}_1$ corresponding to the sequence of angles of $T$ intersected by $\gamma$. 
 
Combining this bijection with the definition of $B^T$ and the construction of all the band modules in \cite{BR} (associated with power of bands) we obtain a natural bijection
$$\xymatrix{b^T:(\pi_1^{\rm free}(S)\times k^*)/\sim  \ar[r] & \{ \textrm{band modules in }\Jac(Q,W)\}/ \textrm{iso}},$$
where the equivalence relation in $\pi_1^{\rm free}(S)\times k^*$ is generated by $(\gamma,\lambda)\sim(\gamma^{-1},\lambda^{-1})$. 

Note that this bijection is not exactly the one described in \cite{ABCP, BZ}. The one in \cite{ABCP,BZ} needs to make a choice of an orientation for each element in the set $\pi_1^{\rm free,irred}(S)/\gamma\sim\gamma^{-1}$, choice which is not canonical.

\subsubsection{Example}\label{subsection exemple}
Let $(S,\mathcal M)$ be the annulus with $3$ marked points, and consider the following triangulation $T$ on $(S,\mathcal{M})$ corresponding to the following quiver $Q_T^{\rm op}$.

\[\scalebox{1}{\begin{tikzpicture}[scale=0.6,>=stealth]
\draw (0,0) circle (1);
\draw (0,0) circle (4);
\draw[red] (0,0) circle (3);
\node at (0,4) {$\bullet$};
\node at (0,-1) {$\bullet$};
\node at (0,-4) {$\bullet$};
\draw[red,->] (3,-0.1)--(3,0.1);
\node[red] at (3.3,0) {$\gamma$};
\draw (0,-4) --node [fill=white, inner sep=0pt] (A2){$\textbf 2$} (0,-1);
\draw(0,-1) ..node[fill=white, inner sep=0pt](A1){$\textbf 1$} controls (-2,-1) and (-2,2)..(0,4);
\draw(0,-1) ..node[fill=white, inner sep=0pt](A3){$\textbf 3$} controls (2,-1) and (2,2)..(0,4);
\draw[blue,->,thick] (A1).. controls (-1,2) and (1,2)..(A3);
\draw[blue, ->,thick] (A1).. controls (-2,-2) and (-1,-2.5).. (A2);
\draw[blue, <-,thick] (A3).. controls (2,-2) and (1,-2.5).. (A2);
\end{tikzpicture}}\]

Let $\gamma$ be the following simple generator of $\pi_1(S)$, and $\lambda$ be in $k^*$. The module $B^T(\gamma,\lambda)$ is isomorphic to the following representation: 

\[\scalebox{1}{\begin{tikzpicture}[scale=0.5,>=stealth]
\node (A1) at (0,0) {$k$};
\node (A2) at (3,3) {$k$};
\node (A3) at (6,0) {$k$};

\draw[->] (A1)--node[yshift=5pt,xshift=-5pt]{$1$}(A2);
\draw[->] (A2)--node[yshift=5pt,xshift=5pt]{$\lambda$}(A3);
\draw[->] (A1)--node[yshift=-6pt]{$1$}(A3);
\node at (-2,1) {$B^T(\gamma,\lambda)\simeq$};
\end{tikzpicture}}\]

Define the triangulations $T^{\textbf 1}:=\mathfrak f_{\textbf 1}(T)$ and $T^{\textbf 2}:=\mathfrak{f}_{\textbf 2}(T)$. Then a direct computation gives the isomorphisms:

\[\scalebox{1}{\begin{tikzpicture}[scale=0.5,>=stealth]
\node (A1) at (0,0) {$k$};
\node (A2) at (3,3) {$k$};
\node (A3) at (6,0) {$k$};

\draw[<-] (A1)--node[yshift=6pt,xshift=-6pt]{$-1$}(A2);
\draw[->] (A2)--node[yshift=5pt,xshift=5pt]{$\lambda$}(A3);
\draw[<-] (A1)--node[yshift=-6pt]{$1$}(A3);
\node at (-4,1) {$\mu_{\textbf 1}(B^T(\gamma,\lambda))\simeq$};
\node at (10,1) {$\simeq B^{T^{\textbf 1}}(\gamma,-\lambda).$};
\end{tikzpicture}}\]

Note that here, since $\textbf 1$ is a sink in $Q^{\rm op}$, $\widetilde{\mu}_{\textbf 1}(Q^{\rm op},W)$ is already reduced and so the canonical reduction is the identity morphism.

Other direct computations give the following isomorphisms:
\[\scalebox{1}{\begin{tikzpicture}[scale=0.5,>=stealth]
\node (A1) at (0,0) {$k$};
\node (A2) at (3,3) {$0$};
\node (A3) at (6,0) {$k$};

\draw[<-] (A1)--(A2);
\draw[<-] (A2)--(A3);

\draw[->] (A1)--node[yshift=6pt]{$1$}(A3);
\draw[->] (0.5,-0.5)--node[yshift=-6pt]{$\lambda$}(5.5,-0.5);
\node at (-4,1) {$\mu_{\textbf 2}(B^T(\gamma,\lambda))\simeq$};
\node at (10,1) {$\simeq B^{T^{\textbf 2}}(\gamma,\lambda).$};

\end{tikzpicture}}\]

\[\scalebox{1}{\begin{tikzpicture}[scale=0.5,>=stealth]
\node (A1) at (0,0) {$k$};
\node (A2) at (3,3) {$k$};
\node (A3) at (6,0) {$k$};

\draw[->] (A1)--node[yshift=6pt,xshift=-6pt]{$-\lambda$}(A2);
\draw[->] (A2)--node[yshift=5pt,xshift=5pt]{$1$}(A3);
\draw[->] (A1)--node[yshift=-6pt]{$1$}(A3);
\node at (-4,1) {$\mu_{\textbf 2}(B^{T^{\textbf 2}}(\gamma,\lambda))\simeq$};
\node at (10,1) {$\simeq B^{T}(\gamma,-\lambda)$};
\end{tikzpicture}}\]

This implies that $\mu_{\textbf 2}$ (with the canonical reduction) is not an involution, and therefore the autoequivalence $\Phi_{\textbf 2}^2$ (defined with the canonical reduction) of $\mathcal{C}_{T}$ is not isomorphic to the identity functor. This was already noticed in \cite[Thm 10.13, (10.23)]{DWZ}.

\end{document}